\documentclass[a4paper,reqno]{amsart}
\usepackage{amsfonts}
\theoremstyle{plain}
\begingroup
\newtheorem{theorem}{Theorem}[section]
\newtheorem{lemma}[theorem]{Lemma}
\newtheorem{proposition}[theorem]{Proposition}

\endgroup

\theoremstyle{definition}
\begingroup
\newtheorem{definition}[theorem]{Definition}
\newtheorem{remark}[theorem]{Remark}

\endgroup

\theoremstyle{remark}
\begingroup

\endgroup

\mathchardef\emptyset="001F

\numberwithin{equation}{section}

\newcommand{\e}{\varepsilon}
\newcommand{\Om}{\Omega}
\newcommand{\Ga}{\Gamma}
\newcommand{\C}{{\mathbb C}}
\newcommand{\KK}{{\mathbb K}}
\newcommand{\R}{{\mathbb R}}
\newcommand{\Rn}{{\R}^n}

\newcommand{\N}{\mathbb N}
\newcommand{\Q}{\mathbb Q}
\newcommand{\Mnn}{{\mathbb M}^{n{\times}n}_{sym}}
\newcommand{\MD}{{\mathbb M}^{n{\times}n}_D}
\newcommand{\wto}{\rightharpoonup}

\newcommand{\setmeno}{\!\setminus\!}
\newcommand{\into}{{\int_{\Omega}}}
\renewcommand{\div}{{\rm div}}
\newcommand{\supp}{{\rm supp}}
\newcommand{\tr}{{\rm tr}}
\newcommand{\hn}{{\mathcal H}^{n-1}}
\newcommand{\Ln}{{\mathcal L}^n}
\newcommand{\LL}{{\mathcal L}}

\newcommand{\QQ}{{\mathcal Q}}
\newcommand{\D}{{\mathcal D}}

\newcommand{\HH}{{\mathcal H}}
\newcommand{\K}{{\mathcal K}}
\newcommand{\dist}{{\rm dist}}

\newcommand{\V}{{\mathcal V}}

\newcommand{\ol}{\overline}
\newcommand{\tki}{t_k^i}
\newcommand{\tkim}{t_k^{i-1}}
\newcommand{\wki}{w_k^i}
\newcommand{\uki}{u_k^i}
\newcommand{\pki}{p_k^i}
\newcommand{\uk}{u_k}
\newcommand{\pk}{p_k}

\title[Quasistatic evolution problems in perfect plasticity]
{Quasistatic evolution problems for\\
 linearly elastic - perfectly plastic materials}
\author[Gianni Dal Maso]{Gianni Dal Maso}
\author[Antonio DeSimone]{Antonio DeSimone}
\author[Maria Giovanna Mora]{Maria Giovanna Mora}
\address[G.~Dal Maso, A.~DeSimone, and M.G.~Mora]{SISSA, Via Beirut 2-4, 
34014 Trieste, Italy}
\email[Gianni Dal Maso]{dalmaso@sissa.it}
\email[Antonio DeSimone]{desimone@sissa.it}
\email[Maria Giovanna Mora]{mora@sissa.it}

\begin{document}
\begin{abstract}
The problem of quasistatic evolution in small strain associative elastoplasticity is studied 
in the framework of the variational theory for rate-independent processes. 
Existence of solutions is proved through the use of incremental variational problems 
in spaces of functions with bounded deformation. This provides a new approximation result for the solutions of the quasistatic evolution problem, which are shown to be absolutely continuous in time. Four equivalent formulations of the problem in rate form are derived. A strong formulation of the flow rule is obtained by introducing a precise definition of the stress on the singular set of the plastic strain.
\end{abstract}
\maketitle

{\small

\bigskip
\keywords{\noindent {\bf Keywords:} 
quasistatic evolution, rate-independent processes, perfect plasticity, Prandtl-Reuss plasticity, shear bands, incremental problems, variational problems in $BD$.}

\subjclass{\noindent {\bf 2000 Mathematics Subject Classification:} 74C05 (74G65, 49J45, 47J20, 35Q72)}
}
\tableofcontents
\bigskip
\bigskip

\begin{section}{Introduction}
In this paper we study quasistatic evolution problems in small strain associative elastoplasticity. 
More precisely, we consider the case of a material whose elastic behaviour is linear and isotropic, 
and whose plastic response is governed by the Prandtl-Reuss flow rule, without hardening (perfect plasticity).

This is a classical problem in mechanics and it is usually formulated as follows in a domain $\Om\subset\Rn$. 
The linearized strain $Eu$, defined as the symmetric part of the spatial gradient of the displacement $u$, 
is decomposed as the sum $Eu=e+p$, where $e$ and $p$ are the elastic and plastic strains. The stress $\sigma$ 
is determined only by $e$, through the formula $\sigma=\C e$, where $\C$ is the elasticity tensor. It is constrained 
to lie in a prescribed subset $\KK $ of the space $\Mnn$ of $n{\times}n$ symmetric matrices, whose boundary 
$\partial \KK $ is referred to as the yield surface.  

Given a time-dependent body force $f(t,x)$, the classical formulation of the quasistatic evolution problem in a time interval $[0,T]$ 
consists in finding functions $u(t,x)$, $e(t,x)$, $p(t,x)$, $\sigma(t,x)$ satisfying the following conditions 
for every $t\in[0,T]$ and every $x\in\Om$:
\begin{itemize}
\smallskip
        \item[(cf1)] additive decomposition: $Eu(t,x)=e(t,x)+p(t,x)$,
        \item[(cf2)] constitutive equation: $\sigma(t,x)=\C e(t,x)$,
        \item[(cf3)] equilibrium: $-\div\,\sigma(t,x)=f(t,x)$,
        \item[(cf4)] associative flow rule: $(\xi-\sigma(t,x)){\,:\,}\dot p(t,x)\le 0$ for every $\xi\in \KK $,
        \smallskip
\end{itemize}
where the colon denotes the scalar product between matrices. The problem is supplemented by initial conditions 
at time $t=0$ and by boundary conditions for $t\in [0,T]$, $x\in\partial\Om$, of the form $u(t,x)=w(t,x)$ 
on a portion $\Ga_0$ of the boundary, and $\sigma(t,x)\nu(x)=g(t,x)$ on the complementary portion $\Ga_1$, 
where $\nu(x)$ is the outer unit normal to~$\partial\Om$, $w(t,x)$ is the prescribed displacement on $\Ga_0$, 
and $g(t,x)$ is the prescribed surface force on $\Ga_1$.

For concreteness, we focus on the case where $\KK $ is a cylinder of the form $\KK =K+\R I$, 
where $I$ is the identity matrix and $K$ is a convex compact subset of $\MD$, 
the space of trace free $n{\times}n$ symmetric matrices. 
This corresponds to yield criteria, often used for metals, which are insensitive to pressure, such as the ones of Tresca and von Mises (see, e.g.,~\cite{Lub}). Then condition (cf4) implies that $\dot p(t,x)\in \MD$ and it is not restrictive to assume that $p(t,x)\in \MD$.

Introducing the normal cone $N_K(\xi)$ to $K$ at $\xi$, the support function
$$
H(\xi):=\sup_{\zeta\in K} \xi{\,:\,}\zeta \,,
$$
and the subdifferential $\partial H(\xi)$ of $H$ at $\xi$, 
the flow rule (cf4) can be written in the equivalent forms 
(see, e.g., \cite[Chapter~4]{Han-Red}):
\begin{itemize}
\smallskip
  \item[$\rm(cf4')$] normality: $\dot p(t,x)\in N_K(\sigma_D(t,x))$,
    \item[$\rm(cf4'')$] flow rule in primal formulation: $\sigma_D(t,x)\in \partial H(\dot p(t,x))$,
         \item[$\rm(cf4''')$] maximal dissipation: $H(\dot p(t,x))=\sigma_D(t,x){\,:\,}\dot p(t,x)$,
         \smallskip
 \end{itemize}
where $\sigma_D(t,x)$ denotes the deviator of $\sigma(t,x)$ (see Section~\ref{mathprel}).

In the engineering literature quasistatic evolution problems of the type considered above are approximated numerically  
by solving a finite number of incremental variational problems (see \cite{Mar}, \cite{Ort-Mar}, and, more recently, 
\cite{Car-Hac-Mie}, \cite{Miehe}, \cite{Ort-Sta}). The time interval $[0,T]$ is divided into $k$ subintervals 
by means of points
$$
0=t_k^0<t_k^1<\dots<t_k^{k-1}<t_k^{k}=T\,,
$$
and the approximate solution $u^i_k$, $e^i_k$, $p^i_k$ at time $t^i_k$ is defined, inductively,  
as a minimizer of the incremental problem
\begin{equation}\label{minintro}
\min_{(u,e,p)\in A(w(t^i_k))} \{\QQ(e)+ \HH(p-p^{i-1}_k)-\langle \LL(t_i^k)|u\rangle
 \}\,,
\end{equation}
where 
\begin{eqnarray}
&\displaystyle \QQ(e):={\textstyle\frac12} \int_\Om \C e(x){\,:\,}e(x)\,dx \,,
\qquad
\HH(p):= \int_\Om H(p(x))\,dx\,,\nonumber
\\
\label{langleLL}
&
\displaystyle\langle \LL(t)|u\rangle :=\int_\Om f(t,x)\,u(x)\,dx +
\int_{\Ga_1} g(t,x)\, u(x) \,d\hn(x) \,,
\end{eqnarray}
$\hn$ is the $(n-1)$-dimensional Hausdorff measure, 
and $A(w(t))$ is defined, at this stage of the discussion, as the set of triples $(u,e,p)$, with $Eu(x)=e(x)+p(x)$ for every $x\in\Om$, 
such that $u$ satisfies the prescribed Dirichlet boundary condition at time $t$, i.e., $u(x)=w(t,x)$ for every $x\in\Ga_0$. 
Finally, the stress at time  $t^i_k$ is obtained as $\sigma^i_k(x):=\C e^i_k(x)$.

Since $\HH$ has linear growth, problem (\ref{minintro}) has, in general, no solution in Sobolev spaces. 
This is very natural from the point of view of mechanics, due to the phenomenon of strain localization. In the absence of hardening, 
solutions can develop shear bands, where shear deformation concentrates. 
Seen from a macroscopic perspective, shear bands can be thought of
as sharp discontinuities of the displacement (slip surfaces). They cannot be resolved by Sobolev functions, but they find a natural mathematical representation if plastic deformations are allowed 
to take values in spaces of measures (see~\cite{Suq}).

These remarks lead naturally to a weak formulation of the problem, where the displacement $u$ belongs to the space 
$BD(\Om)$ of functions with bounded deformation, whose theory was developed in 
\cite{Mat}, \cite{Tem-Stra}, \cite{Koh-Tem}, \cite{Tem}, and the plastic strain $p$ belongs to the space 
$M_b(\Om\cup\Ga_0;\MD)$ of $\MD$-valued bounded Borel measures on $\Om\cup\Ga_0$.

In accordance to the theory of convex functions of measures developed in \cite{Gof-Ser} and
\cite[Chapter~II, Section~4]{Tem}, we define the functional $\HH(p)$ in the weak formulation of problem (\ref{minintro}) as
$$
\HH(p):= \int_{\Om\cup\Ga_0} H(p/|p|)\,d|p|\,,
$$
where $p/|p|$ is the Radon-Nikodym derivative of the measure $p$ with respect to its variation $|p|$, while 
$A(w(t^i_k))$ is defined, here and henceforth, as the set of triples $(u,e,p)$, with $u\in BD(\Om)$, 
$e\in L^2(\Om;\Mnn)$, $p\in M_b(\Om\cup\Ga_0;\MD)$, and $Eu=e+p$ on $\Om$, subject to the relaxed boundary condition 
$p=(w(t^i_k)-u){\,\odot\,}\nu\,\hn$ on $\Ga_0$. In the last formula $\odot$  denotes the symmetric tensor product.

Boundary conditions of this kind are typical in the variational theory of functionals with linear growth 
(see, e.g., \cite{Tem} and \cite{Giu}). The mechanical interpretation of our condition on $\Ga_0$ is that, 
if the prescribed boundary displacement is not attained, a plastic slip is developed at the boundary, 
whose strength is proportional to the difference between the prescribed and the attained boundary displacements.

In the case $p^{i-1}_k=0$ the weak formulation of problem (\ref{minintro}) has been studied in detail in 
\cite{Tem-Stra}, \cite{Anz-Gia}, \cite{Koh-Tem}, \cite{Tem}, and \cite{Anz} at the beginning of the 80's. 
With respect to this body of work, it is important to emphasize a change of perspective. 
The model we study (Prandtl-Reuss plasticity) takes explicitly into account the history of plastic deformation. 
Setting instead  $p^{i-1}_k=0$ in (\ref{minintro}) makes the problem oblivious to the accumulation of plastic strain. 
This is the so called Hencky theory of plasticity, in which elastic unloading following plastic loading is not correctly 
resolved (see \cite{Hill} and \cite{Suq}). 

We can rely however on the results of the above mentioned papers to solve problem (\ref{minintro}) in the general case 
(Theorem~\ref{existencemin}), provided a safe-load condition is satisfied. Then we define the piecewise constant interpolations
\begin{equation*}
\uk(t):=\uki\,, \quad e_k(t):=e_k^i\,,  \quad   \pk(t):=\pki\,,  \quad
\sigma_k(t):=\sigma_k^i\,, 
\end{equation*}
where $i$ is the largest integer such that $\tki\le t$. 

The aim of this paper is to introduce a weak definition of continuous-time quasistatic evolution in the functional 
framework $u\in BD(\Om)$, $e\in L^2(\Om;\Mnn)$, $p\in M_b(\Om\cup\Ga_0;\MD)$, $\sigma\in L^2(\Om;\Mnn)$, 
and to prove that, up to a subsequence, the discrete-time solutions $\uk(t)$, $e_k(t)$, $\pk(t)$, $\sigma_k(t)$, 
obtained by solving the weak formulations of problems (\ref{minintro}), converge to a continuous-time solution 
$u(t)$, $e(t)$, $p(t)$, $\sigma(t)$, provided $\max_i (t^i_k-t^{i-1}_k)\to0$ as $k\to\infty$. 

Our definition fits the general scheme of continuous-time energy formulation of 
rate-independent processes developed in \cite{Mie-The}, \cite{Mie-The-Lev}, \cite{Mie-Beijing}, 
\cite{Mie}, \cite{Mie-Rou}, and \cite{Mai-Mie}. Following those papers, for every time interval $[s,t]$ 
contained in $[0,T]$ we introduce the dissipation associated with $\HH$, defined by
\begin{equation*}
\D_\HH(p;s,t):=\sup\Big\{ \sum_{j=1}^N \HH(p(t_j)-p(t_{j-1})):\, 
s=t_0\le t_1\le \dots\le t_N=t, \,
N\in\N \Big\}\, .
\end{equation*}
The general definition proposed in \cite{Mai-Mie} reads in our case as follows: 
a quasistatic evolution is a function $t\mapsto(u(t),e(t),p(t))$ from 
$[0,T]$ into $BD(\Om){\times}L^2(\Om;\Mnn){\times}M_b(\Om\cup\Ga_0;\MD)$
which satisfies the following conditions:
\begin{itemize}
\smallskip
\item[{\rm (qs1)}] global stability: for every $t\in [0,T]$ we have $(u(t),e(t),p(t))\in 
A(w(t))$ and
$$
\QQ(e(t))-\langle \LL(t)|u(t)\rangle 
\le \QQ(\eta)+\HH(q-p(t))
-\langle \LL(t)|v\rangle
$$
for every $(v,\eta,q)\in A(w(t))$;  
\smallskip 
\item[{\rm (qs2)}] energy balance: the function $t\mapsto p(t)$ from $[0,T]$ into 
$M_b(\Om\cup\Ga_0;\MD)$ has bounded variation and for every $t\in[0,T]$
\begin{eqnarray*}
&\displaystyle \QQ(e(t)) + \D_\HH(p; 0,t)-\langle \LL(t)|u(t)\rangle = 
\QQ(e(0))-\langle \LL(0)|u(0)\rangle +{}
\\
&\displaystyle
{}+ \int_0^t \{\langle \sigma(s)|E\dot w(s)\rangle -\langle \LL(s)|\dot w(s)\rangle 
- \langle\dot \LL(s)| u(s)\rangle
 \}\, ds\,,
\end{eqnarray*}
where $\sigma(t):=\C e(t)$, dots denote time derivatives, the first brackets 
$\langle \cdot|\cdot\rangle$ in the integral denote the scalar product in $L^2(\Om;\Mnn)$, while the other brackets $\langle \cdot|\cdot\rangle$ are defined as in~(\ref{langleLL}).
\smallskip
\end{itemize}

The main result of the present paper is the proof of the existence of a quasistatic evolution satisfying prescribed 
initial conditions (Theorem~\ref{main}), provided a uniform safe-load condition is satisfied. 

A different formulation of the problem in rate form was proposed in \cite{Joh} and \cite{Suq}, where an existence result 
is proved by a visco-plastic approximation.
It turns out that our definition is equivalent to the one considered in those papers (Theorem~\ref{varineq} and Remark~\ref{Suquet}). 
Therefore the existence result is not new, but our proof is completely different and leads to a different 
approximation of the solutions (Theorem~\ref{convstress}). Moreover it shows that this problem can be included 
in the general theory developed in \cite{Mie-Beijing} and \cite{Mai-Mie}.

Our proof is obtained by considering the discrete-time solutions $\uk(t)$, $e_k(t)$, $\pk(t)$, $\sigma_k(t)$ and 
by showing that they satisfy an approximate energy inequality (Lemma~\ref{lm:41}), which is similar to 
\cite[Theorem~4.1]{Mai-Mie}. This allows us to apply the generalization  (Lemma~\ref{lm:44}) of the classical 
Helly Theorem proved in \cite[Theorem~3.2]{Mai-Mie}, and to extract a subsequence, independent of $t$ and still 
denoted $p_k$, such that $p_k(t)\wto p(t)$ weakly$^*$ in $M_b(\Om\cup\Ga_0;\MD)$ for every $t\in [0,T]$. 

Extracting a further subsequence, possibly depending on $t$, we may assume that $\uk(t)\wto u(t)$ weakly$^*$ in 
$BD(\Om)$ and $e_k(t)\wto e(t)$ weakly in $L^2(\Om;\Mnn)$. We prove (Theorem~\ref{lm:45}) that
$(u(t),e(t),p(t))$ satisfies the global stability condition~(qs1). Since there exists at most one $(u,e)\in BD(\Om){\times}L^2(\Om;\Mnn)$
such that $(u,e,p(t))$ satisfies~(qs1) (Remark~\ref{unique}), 
we have $\uk(t)\wto u(t)$ and $e_k(t)\wto e(t)$ for the same subsequence (independent of~$t$) 
for which $p_k(t)\wto p(t)$.

One of the inequalities in the energy balance (qs2) is then proved by passing to the limit in the approximate energy 
inequality obtained for the discrete-time solutions, while the opposite inequality follows (Theorem~\ref{inequality}) 
from the global stability, by adapting the proofs of \cite[Theorem~4.4]{Mai-Mie} and \cite[Lemma~7.1]{DM-Fra-Toa}.

The second part of the paper is devoted to the regularity of solutions and to the comparison of our definition of 
quasistatic evolution with other definitions in rate form. We prove (Theorem~\ref{ACesigma}) 
that, if the data of the problem are absolutely continuous functions of time, then for every quasistatic evolution the functions $t\mapsto u(t)$, $t\mapsto e(t)$,  $t\mapsto p(t)$, and 
$t\mapsto\sigma(t)$ are absolutely continuous on $[0,T]$ with values in $BD(\Om)$, $L^2(\Om;\Mnn)$, 
$M_b({\Om\cup\Ga_0};\MD)$, $L^2(\Om;\Mnn)$, respectively. Moreover, we establish a pointwise estimate for the time derivatives of these functions which implies that,
if the data of the problem are Lipschitz continuous on $[0,T]$, then
the same is true for $t\mapsto u(t)$, $t\mapsto e(t)$,  $t\mapsto p(t)$, and $t\mapsto\sigma(t)$ 
(Remark~\ref{Lipschitz}).

Similar arguments prove that $t\mapsto e(t)$ and $t\mapsto \sigma(t)$ are uniquely determined by their initial 
conditions (Theorem~\ref{uniqueness}), while elementary examples in dimension one show that, in general, 
this is not true for $t\mapsto u(t)$ and ${t\mapsto p(t)}$ (see \cite[Section~2.1]{Suq}).

These regularity results allow us (Proposition~\ref{equivalence}) to write the energy balance (qs2) as balance of powers: for a.e.\ $t\in[0,T]$
$$
\langle\sigma(t)|\dot e(t)\rangle +\HH(\dot p(t))=
\langle \LL(t)|\dot u(t)\rangle + 
\langle\sigma(t)|E\dot w(t)\rangle-\langle \LL(t)|\dot w(t)\rangle\,.
$$
We then show that our definition of quasistatic evolution is equivalent to four different sets of conditions, 
expressed in rate form (Theorems~\ref{varineq} and~\ref{strong}). One of them can be interpreted as the weak formulation, 
in the spaces $BD(\Om)$, $L^2(\Om;\Mnn)$, $M_b({\Om\cup\Ga_0};\MD)$, $L^2(\Om;\Mnn)$, of the four conditions 
(cf1)--(cf4), considered in the classical presentation of the problem; another one takes into account 
the weak formulation of maximal dissipation $\rm (cf4''')$; the third one coincides with the definition considered in~\cite{Suq};  the last one (Theorem~\ref{strong} and Remark~\ref{remstrong}) presents 
a strong formulation of the normality rule in both forms~$\rm(cf4')$ and~$\rm(cf4'')$. This requires a precise representative of $\sigma_D(t)$ defined $|\dot p(t)|$-a.e.\ on $\Om\cup\Ga_0$. If $K$ is strictly convex, this representative is obtained as limit of averages of $\sigma_D(t)$ (Theorem~\ref{strong2}).

\end{section}

\begin{section}{Notation and preliminary results}

\subsection{Mathematical preliminaries}\label{mathprel}{\ }

\medskip
\noindent
{\bf Measures.}
The Lebesgue measure on $\Rn$ is denoted by $\Ln$, and the 
$(n-1)$-dimen\-sional Hausdorff measure by $\hn$. 
Given a Borel set $B\subset\Rn$ and a finite dimensional Hilbert space 
$X$, $M_b(B;X)$ denotes the space of 
bounded Borel measures on $B$  with values in $X$, endowed with the norm 
$\|\mu\|_1:=|\mu|(B)$, where 
$|\mu|\in M_b(B;\R)$ is the variation of the measure $\mu$. 
For every $\mu\in M_b(B;X)$ we consider the Lebesgue decomposition 
$\mu=\mu^a+\mu^s$, where $\mu^a$ is 
absolutely continuous and $\mu^s$ is singular with respect to Lebesgue 
measure $\Ln$. 

If $\mu^s=0$, we always 
identify $\mu$ with its density with respect to Lebesgue measure $\Ln$. In 
this way
$L^1(B;X)$ is regarded as a subspace of $M_b(B;X)$, with the induced norm. 
In particular $\mu^a\in L^1(B;X)$ 
for every $\mu\in M_b(B;X)$. The indication of the space $X$ is omitted when $X=\R$.
The $L^p$ norm, $1\le p\le\infty$, is denoted 
by $\|\cdot\|_p$. The brackets $\langle \cdot|\cdot\rangle$ denote the duality product between conjugate $L^p$ spaces, as well as between other pairs of spaces, according to the context.

If the relative topology of $B$ is locally compact, by Riesz representation theorem 
(see, e.g., \cite[Theorem~6.19]{Rud})
$M_b(B;X)$ can be identified with the dual of 
$C_0(B;X)$, the space of continuous functions $\varphi\colon B\to X$ such 
that $\{|\varphi |\ge\e\}$ is compact for every $\e>0$. The 
weak$^*$ topology of $M_b(B;X)$ is defined using this duality.

\medskip

\noindent {\bf Matrices.}
The space of {\it symmetric $n{\times}n$ matrices\/} is denoted by $\Mnn$; it is 
endowed with the euclidean scalar product 
$\xi{\,:\,}\zeta:=\tr(\xi\zeta)=\sum_{ij}\xi_{ij}\zeta_{ij}$ and with the 
corresponding euclidean norm 
$|\xi|:= (\xi{\,:\,}\xi)^{1/2}$. 
The orthogonal complement of the subspace $\R I$ spanned by the identity 
matrix $I$ is the subspace $\MD$ of all matrices of 
$\Mnn$ with trace zero. For every $\xi\in\Mnn$ the orthogonal projection 
of $\xi$ on $\R I$ is $\frac1n \tr(\xi)I$, while the 
orthogonal projection on $\MD$ is the {\it deviator\/} $\xi_D$ of $\xi$, so that 
we have the orthogonal decomposition
$$
\textstyle
\xi=\xi_D+\frac1n (\tr\,\xi) I\,.
$$
The {\it symmetrized tensor product\/} $a{\,\odot\,}b$ of two vectors $a$, 
$b\in\Rn$ is the symmetric matrix with 
entries $(a_ib_j+a_jb_i)/2$. 
It is easy to see that $\tr(a{\,\odot\,}b)=a{\,\cdot\,}b$, the scalar 
product of $a$ and $b$, and that $|a{\,\odot\,}b|^2=\frac12 |a|^2 |b|^2 + \frac12 (a{\,\cdot\,}b)^2$, so that $\frac1{\sqrt2} |a||b|\le |a{\,\odot\,}b| \le |a||b|$.

\medskip

\noindent {\bf Functions with bounded deformation.}
Let $U$ be an open set in $\Rn$. For every $u\in L^1(U;\Rn)$ let $Eu$ be the 
$\Mnn$-valued distribution on $U$, whose components are defined by $E_{ij}u=\frac12 (D_j u_i +D_i u_j)$. The space $BD(U)$ of functions with {\it bounded deformation\/} is the space of all $u\in L^1(U;\Rn)$ such that $Eu\in M_b(U;\Mnn)$. It is easy to see that $BD(U)$ is a Banach space with the norm
$$
\|u\|_1+\|Eu\|_1\,.
$$
It is possible to prove that $BD(U)$ is the dual of a normed space
(see \cite{Mat} and \cite{Tem-Stra}). The weak$^*$ topology of $BD(U)$ is defined using this duality. A sequence $u_k$ converges to $u$ weakly$^*$ in $BD(U)$ if and only if $u_k\wto u$ weakly in $L^1(U;\Rn)$ and $Eu_k\wto Eu$ weakly$^*$ in $M_b(U;\Mnn)$. Every bounded sequence in $BD(U)$ has a weakly$^*$ convergent subsequence. Moreover, if $U$ is bounded and has Lipschitz boundary, every bounded sequence in $BD(U)$ has a subsequence which converges weakly in $L^{n/(n-1)}(U;\Rn)$ and strongly in $L^p(U;\Rn)$ for every ${p<n/(n-1)}$.
For the general properties of $BD(U)$ we refer to \cite{Tem}.

In our problem $u\in BD(U)$ represents the {\it displacement\/} of an elasto-plastic body and $Eu$ is the corresponding linearized {\it strain\/}.

\subsection{Mechanical preliminaries}\label{mech}{\ }

\medskip
\noindent
{\bf The reference configuration.}
Throughout the paper $\Om$ is a {\it bounded connected open set\/} in $\Rn$ with 
$C^2$ {\it boundary\/}. We suppose that the 
boundary $\partial\Om$ 
is partitioned into two disjoint open sets $\Ga_0$, $\Ga_1$ and their 
common boundary $\partial \Ga_0=\partial \Ga_1$ 
(topological notions refer here to the relative topology of 
$\partial\Om$). We assume that $\Ga_0\neq \emptyset$ and that 
for every $x\in \partial \Ga_0=\partial \Ga_1$ there exists a $C^2$ 
diffeomorphism defined in an open neighbourhood 
of $x$ in $\Rn$ 
which maps $\partial\Om$ to an $(n-1)$-dimensional plane and $\partial 
\Ga_0=\partial \Ga_1$ to an $(n-2)$-dimensional plane.

On $\Ga_0$ we will prescribe a Dirichlet boundary condition.
This will be done by assigning a function  $w\in H^{1/2}(\Ga_0;\Rn)$, or, equivalently, a function $w\in H^1(\Rn;\Rn)$, whose trace on $\Ga_0$
(also denoted by $w$) is the prescribed boundary value. 
The set $\Ga_1$ will be the part of the boundary on which the traction is prescribed. 

Every function $u\in BD(\Om)$ has a {\it trace\/} on $\partial\Om$, still denoted by $u$, which belongs to $L^1(\partial\Om;\Rn)$. If $u_k$, $u\in BD(\Om)$, $u_k\to u$ strongly in $L^1(\Om;\Rn)$, and $\|Eu_k\|_1 \to \|Eu\|_1$, then $u_k\to u$ strongly in $L^1(\partial\Om;\Rn)$ (see \cite[Chapter~II, Theorem~3.1]{Tem}). Moreover, there exists a constant $C>0$, 
depending on $\Om$ and $\Ga_0$, such that
\begin{equation}
\label{seminorm}
\|u\|_{1,\Om} \le C\, \|u\|_{1,\Ga_0}+ 
C\, \|Eu\|_{1,\Om}
\end{equation}
(see \cite[Proposition~2.4 and Remark~2.5]{Tem}).

We shall frequently use the space $M_b(\Om\cup\Ga_0;\MD)$, which is the dual of $C_0(\Om\cup\Ga_0;\MD)$. The latter space can be identified with the space of functions in $C(\ol\Om;\MD)$ vanishing on $\ol\Ga_1$. The duality product is defined by
\begin{equation}\label{dualmeas}
\langle \tau|\mu\rangle:= \int_{\Om\cup\Ga_0}\!\!\! \!\!\!\tau{\,:\,}d\mu:= \sum_{ij} \int_{\Om\cup\Ga_0}\!\!\!\!\!\! \tau_{ij}\,d\mu_{ij}
\end{equation}
for every $\tau=(\tau_{ij})\in C(\ol\Om;\MD)$ and every $\mu=(\mu_{ij})\in M_b(\Om\cup\Ga_0;\MD)$.

\medskip

\noindent {\bf The set of admissible stresses.}
Let $K$ be a closed convex set in $\MD$, which will play the role of a constraint on the deviatoric part of the stress. Its boundary is interpreted as the {\it yield surface\/}. We assume that there exist two constants $r_K$ and $R_K$, with $0<r_K\le R_K<\infty$, such that
\begin{equation}\label{rk}
\{\xi\in\MD: |\xi|\le r_K\}\subset K\subset \{\xi\in\MD: |\xi|\le R_K\}\,.
\end{equation}
It is convenient to introduce the convex set 
$$
\K_D(\Om):=\{\tau\in L^2(\Om;\MD): \tau(x)\in K\, \hbox{ for a.e. }\,x\in\Om\}\,.
$$
The {\it set of admissible stresses\/} is defined by 
$$
\K(\Om):=\{\sigma\in L^2(\Om;\Mnn): \sigma_D\in\K_D(\Om)\}\,.
$$

The {\it support function\/} 
$H\colon\MD\to {[0,+\infty[}$ of $K$ is given by
\begin{equation}\label{HD}
H(\xi):=\sup_{\zeta\in K} \xi{\,:\,}\zeta \,.
\end{equation}
It turns out that $H$ is convex and positively homogeneous of degree one. In particular it satisfies the triangle inequality
$$
H(\xi+\zeta)\le H(\xi)+H(\zeta)\,.
$$
{}From (\ref{rk}) it follows that
\begin{equation}\label{boundsH}
r_K|\xi|\le H(\xi)\le R_K|\xi| 
\end{equation}
for every $\xi\in \MD$.

{}For every $\mu\in M_b(\Om\cup\Ga_0;\MD)$ let $\mu/|\mu|$ be the Radon-Nikodym derivative of $\mu$  with respect to its variation~$|\mu|$.
Using the theory of convex functions of measures developed in \cite{Gof-Ser}, we introduce the nonnegative Radon measure $H(\mu)\in M_b(\Om\cup\Ga_0)$ defined by $H(\mu):=H(\mu/|\mu|)|\mu|$, i.e.,
\begin{equation}\label{Hmu}
H(\mu)(B):=\int_{B} H(\mu/|\mu|)\,d|\mu|
\end{equation}
for every Borel set $B\subset \Om\cup\Ga_0$.
Finally, we consider the functional $\HH\colon M_b(\Om\cup\Ga_0;\MD)\to\R$ defined by
\begin{equation}\label{calHD}
\HH(\mu):= H(\mu)(\Om\cup\Ga_0)=\int_{\Om\cup\Ga_0} H(\mu/|\mu|)\,d|\mu|\,.
\end{equation}
Using \cite[Theorem~4]{Gof-Ser} and  \cite[Chapter~II, Lemma~5.2]{Tem} we can see that $H(\mu)$ coincides with the measure studied in \cite[Chapter~II, Section~4]{Tem}, hence
\begin{equation}\label{HHduality}
\HH(\mu) = \sup\{\langle \tau|\mu\rangle: 
\tau\in C_0(\Om\cup\Ga_0;\MD)\cap\K_D(\Om)\}
\end{equation}
and $\HH$ is lower semicontinuous on $M_b(\Om\cup\Ga_0;\MD)$ with respect to weak$^*$ convergence.
It follows from the properties of $H$ that $\HH$ satisfies the triangle inequality, i.e., 
\begin{equation}\label{triangle}
\HH(\lambda+\mu)\le \HH(\lambda)+\HH(\mu)
\end{equation}
for every $\lambda,\mu\in M_b(\Om\cup\Ga_0;\MD)$.

\medskip

\noindent {\bf The elasticity tensor.}
Let $\C$ be the {\it elasticity tensor\/}, considered as a symmetric 
positive definite linear operator $\C\colon\Mnn\to\Mnn$. 
We assume that the orthogonal subspaces $\MD$ and $\R I$ are invariant under $\C$.
This is equivalent to saying that there exist a 
symmetric positive definite linear 
operator $\C_D\colon\MD\to\MD$ and a constant
$\kappa>0$ such that
\begin{equation}
\C\xi:= \C_D\xi_D+\kappa (\tr\,\xi) I 
\end{equation}
for every $\xi\in\Mnn$. Note that when $\C$ is isotropic, we have
$\C\xi=2\mu \xi_D + \kappa (\tr\xi)I$, where $\mu>0$ is the shear modulus and $\kappa$ is the modulus of compression, so that our assumptions are satisfied. 

Let $Q\colon\Mnn\to{[0,+\infty[}$ be the quadratic form associated with $\C$, defined by
\begin{equation}\label{W}
\textstyle Q(\xi):=\frac12 \C\xi{\,:\,}\xi=
\frac12\C_D\xi_D{\,:\,}\xi_D+\frac{\kappa}2(\tr\,\xi)^2\
\,.
\end{equation}
It turns out that there exist two constants $\alpha_\C$ and $\beta_\C$, 
with $0<\alpha_\C\le\beta_\C<+\infty$, such that
\begin{equation}\label{boundsC}
\alpha_\C|\xi|^2\le Q(\xi)\le 
\beta_\C|\xi|^2
\end{equation}
for every $\xi\in\Mnn$.
These inequalities imply 
\begin{equation}\label{normC}
|\C\xi|\le 2\beta_\C| \xi| \,.
\end{equation}

It is convenient to introduce the quadratic form $\QQ\colon L^2(\Om;\Mnn)\to\R$ defined by
\begin{equation}\label{calW}
\QQ(e):=\into Q(e)\,dx\
\end{equation}
for every $e\in L^2(\Om;\Mnn)$. It is well known that $\QQ$ is lower semicontinuous on 
$L^2(\Om;\Mnn)$ with respect to weak convergence.

\medskip

\noindent {\bf The prescribed boundary displacements.}
For every $t\in[0,T]$ we prescribe a {\it boundary displacement\/} $w(t)$ in the space $H^1(\Rn;\Rn)$. This choice is motivated by the 
fact that we do not want to impose ``discontinuous'' boundary data, so that, if the displacement develops sharp 
discontinuities, this is due to energy minimization.

We assume also that the function $t\mapsto w(t)$ 
is absolutely continuous from $[0,T]$ into $H^1(\Rn;\Rn)$, so that the 
time derivative
$t\mapsto\dot w(t)$ belongs to $L^1([0,T]; H^1(\Rn;\Rn))$ and 
its strain $t\mapsto E\dot w(t)$ belongs to
$L^1([0,T];L^2(\Rn;\Mnn))$. For the main properties of absolutely continuous functions with values in reflexive Banach spaces we refer to \cite[Appendix]{Bre}.

\medskip

\noindent {\bf Body and surface forces.}
For every $t\in[0,T]$ the {\it body force\/} $f(t)$ belongs to the space $L^n(\Om;\Rn)$ and the {\it surface force\/}
$g(t)$ acting on $\Ga_1$ belongs to $L^{\infty}(\Ga_1;\Rn)$. 
We assume that the functions $t\mapsto f(t)$ and $t\mapsto g(t)$
are absolutely continuous from $[0,T]$ into $L^n(\Om;\Rn)$ and $L^{\infty}(\Ga_1;\Rn)$, respectively,
so that the time derivative $t\mapsto\dot f(t)$
belongs to $L^1([0,T]; L^n(\Om;\Rn))$, the weak$^*$ limit
$$
\dot g(t):= w^*\hbox{-}\lim_{s\to t}\frac{g(s)-g(t)}{s-t}\,,
$$
exists for a.e.\ $t\in[0,T]$, 
and $t\mapsto\|\dot g(t)\|_\infty$ belongs to $L^1([0,T])$ (see Theorem~\ref{thm:a01}).

Throughout the paper we will assume also the following {\it uniform safe-load condition\/}:
there exist a function $t\mapsto \varrho(t)$ from $[0,T]$ into $L^2(\Om;\Mnn)$ and a constant $\alpha >0$
such that for every $t\in[0,T]$
\begin{equation}\label{safeload}
- \div\varrho(t)=f(t)\, \hbox{ a.e.\ on } \Om\,, \qquad [\varrho(t)\nu]=g(t) \hbox{ on } \Ga_1\,,
\end{equation}
and 
\begin{equation}\label{safeload1}
\varrho_D(t,x)+\xi\in K
\end{equation}
for a.e.\  $x\in \Om$ and for every $\xi\in\MD$ with $|\xi|\le\alpha$. 
In these formulas $\varrho_D(t,x)$
denotes the value of  $\varrho_D(t)$ at $x\in\Om$, and the trace $[\varrho(t)\nu]$ of $\varrho(t)\nu$ on $\Ga_1$ is interpreted in the sense of (\ref{sigmanu}) below.
We assume also that the functions $t\mapsto \varrho(t)$ and $t\mapsto \varrho_D(t)$ are
absolutely continuous from $[0,T]$ into $L^2(\Om;\Mnn)$ and $L^\infty(\Om;\MD)$, respectively, so that the time derivative $t\mapsto \dot\varrho(t)$ belongs to $L^1([0,T];L^2(\Om;\Mnn))$,
\begin{equation}\label{dotsigmaD}
\frac{\varrho_D(s)- \varrho_D(t)}{s-t}\ \wto\  \dot\varrho_D(t)
\end{equation}
weakly$^*$ in~$L^\infty(\Om;\MD)$ for a.e.\ $t\in[0,T]$, and $t\mapsto\|\dot\varrho_D(t)\|_\infty$ belongs to $L^1([0,T])$ (see Theorem~\ref{thm:a01}).

\subsection{Stress and strain}
Given a displacement $u\in BD(\Om)$ and a 
boundary datum $w\in H^1(\Rn;\Rn)$, the {\it elastic\/} and {\it plastic strains\/} $e\in 
L^2(\Om;\Mnn)$ and $p\in M_b(\Om\cup\Ga_0;\MD)$ satisfy the equalities 
\begin{eqnarray}
& Eu=e+p \quad \hbox{in }\Om\,, \label{900}
\\
& p=(w-u){\,\odot\,}\nu\,\hn \quad \hbox{on } \Ga_0\,.  \label{901}
\end{eqnarray}
Therefore we have $e=E^au-p^a$ a.e.\ on $\Om$ and $p^s=E^su$ on $\Om$. Since $\tr\,p=0$, 
it follows from (\ref{900}) that $\div\,u=\tr\, e\in L^2(\Om)$
and from (\ref{901}) that $(w-u){\,\cdot\,}\nu=0$ $\hn$-a.e.\ on $\Ga_0$. 
The {\it stress\/} $\sigma\in L^2(\Om;\Mnn)$ is defined by
\begin{equation}\label{sigma}
\sigma:=\C e= \C_D e_D+\kappa\,\tr\, e\,.
\end{equation}
The {\it stored elastic energy\/} is given by
\begin{equation}\label{elenergy}
\QQ(e)=\int_{\Om}Q(e)\, dx={\textstyle\frac12}\langle \sigma | e\rangle \,.
\end{equation}

Given $w\in H^1(\Rn;\Rn)$, the {\it set of admissible displacements and strains\/} for the boundary datum $w$ on $\Ga_0$ is denoted by $A(w)$: it is defined as the set of all triples $(u,e,p)$, with $u\in BD(\Om)$, $e\in L^2(\Om;\Mnn)$, $p\in M_b(\Om\cup\Ga_0;\MD)$, satisfying (\ref{900}) and (\ref{901}).

We shall also use {\it the space $\Pi_{\Ga_0}(\Om)$ of admissible plastic strains\/}, defined as the set of all
$p\in M_b(\Om\cup\Ga_0;\MD)$ for which there exist $u\in BD(\Om)$, $w\in H^1(\Rn;\Rn)$, and $e\in L^2(\Om;\Mnn)$ satisfying (\ref{900}) and (\ref{901}), i.e., $(u,e,p)\in A(w)$.

We now prove a closure property for the multi-valued map $w\mapsto A(w)$.

\begin{lemma}\label{lm:145}
Let $w_k$ be a sequence in $H^1(\Rn;\Rn)$ and let
$(u_k,e_k,p_k)\in A(w_k)$. Assume that $u_k\wto u_\infty$ weakly$^*$ in 
$BD(\Om)$, $e_k\wto e_\infty$ weakly in $L^2(\Om;\Mnn)$,
$p_k\wto p_\infty$ weakly$^*$ in $M_b(\Om\cup\Ga_0;\MD)$, $w_k\wto w_\infty$ 
weakly in $H^1(\Rn;\Rn)$.
Then $(u_\infty,e_\infty,p_\infty)\in A(w_\infty)$. 
\end{lemma}

\begin{proof} 
Since $\Ga_0$ is open in $\partial\Om$, there exists a bounded open set 
$U$ in $\Rn$ such that $\Ga_0=U\cap\partial\Om$, and we define 
$\tilde\Om:=\Om\cup U$. 

For $k=1, 2, \dots, \infty$ let $\tilde u_k\in BD(\tilde\Om)$ be defined by 
$\tilde u_k= u_k$ a.e.\ on $\Om$ and $\tilde  u_k= w_k$ a.e.\ on 
$U\setmeno\Om$. Then 
\begin{equation}\label{EGak}
\begin{array}{c}
E\tilde u_k=Eu_k\quad  \hbox{ on }\Om\,, \\
E\tilde u_k=(w_k-u_k){\,\odot\,}\nu\,\hn \quad \hbox{ on }\Ga_0\,, \\
E\tilde u_k=Ew_k \quad \hbox{ on } U\setmeno\ol\Om\,,
\end{array}
\end{equation}
(see, e.g., \cite[Theorem~2.1 and Remark~2.3]{Tem}).
Since $w_k-u_k$ is bounded in $L^1(\Ga_0;\Rn)$ by the continuity of the trace operator, the sequence  $E\tilde u_k$ is bounded in $M_b(\tilde\Om;\Mnn)$. 
As $\tilde u_k\to\tilde u_\infty$ weakly in $L^1(\tilde\Om;\Rn)$, we conclude that $\tilde u_k\wto\tilde u_\infty$ weakly$^*$ in $BD(\tilde\Om)$.

For $k=1, 2, \dots, \infty$ let $\tilde e_k\in L^2(\tilde\Om;\Mnn)$ be defined by 
$\tilde e_k= e_k$ a.e.\ on $\Om$ and $\tilde  e_k= Ew_k$ 
a.e.\ on $U\setmeno\Om$, and let $\tilde p_k\in M_b(\tilde\Om;\MD)$ 
be defined by $\tilde  p_k= p_k$ 
on $\Om\cup\Ga_0$
and $\tilde  p_k=0$ on $U\setmeno\overline\Om$.
Then $\tilde e_k$ converges to $\tilde e_\infty$ weakly in $L^2(\tilde\Om;\Mnn)$. 
Since the restrictions to $\Om\cup\Ga_0$ of functions in $C_0(\tilde\Om;\MD)$ belong to $C_0(\Om\cup\Ga_0;\MD)$, we obtain also that $\tilde p_k$ converges to $\tilde p_\infty$ weakly$^*$ in $M_b(\tilde\Om;\MD)$.

As $(u_k,e_k,p_k)\in A(w_k)$ for $k<\infty$, using again (\ref{EGak}) we obtain $E\tilde u_k=\tilde e_k+\tilde p_k$ in $\tilde\Om$.
The convergence properties already proved for $(\tilde u_k,\tilde e_k,\tilde p_k)$ show that $E\tilde u_\infty=\tilde e_\infty+\tilde p_\infty$ in $\tilde\Om$.
Consequently, (\ref{EGak}) for $k=\infty$ implies that $(u_\infty,e_\infty,p_\infty)\in A(w_\infty)$. 
\end{proof} 

\medskip

\noindent{\bf The traces of the stress.}
If $\sigma\in L^2(\Om;\Mnn)$ and $\div\,\sigma\in L^2(\Om;\Rn)$, then we 
can define a distribution 
$[\sigma\nu]$ on $\partial\Om$ by
\begin{equation}\label{sigmanu}
\langle [\sigma\nu]|\psi\rangle_{\partial\Om}:=\langle \div\, \sigma|\psi 
\rangle + \langle \sigma|E\psi\rangle \qquad 
\end{equation}
for every $\psi\in H^1(\Om;\Rn)$.
It turns out that $[\sigma\nu]\in H^{-1/2}(\partial\Om;\Rn)$ (see, e.g., 
\cite[Theorem~1.2, Chapter~I]{Tem}).
We will consider the normal and tangential parts of $[\sigma\nu]$, defined 
by
\begin{equation}\label{sigmataunu}
[\sigma\nu]_\nu:=([\sigma\nu]{\,\cdot\,}\nu)\nu\,,\qquad
[\sigma\nu]_\nu^\perp:=[\sigma\nu]-([\sigma\nu]{\,\cdot\,}\nu)\nu\,.
\end{equation}
Since $\nu\in C^1(\partial\Om;\Rn)$, we have that $[\sigma\nu]_\nu$, 
$[\sigma\nu]_\nu^\perp\in H^{-1/2}(\partial\Om;\Rn)$. 
If, in addition, 
$\sigma_D\in L^\infty(\Om;\MD)$, then $[\sigma\nu]_\nu^\perp\in 
L^\infty(\partial\Om;\Rn)$ and
\begin{equation}\label{sigmatau}
\|[\sigma\nu]_\nu^\perp\|_{\infty,\partial\Om} \le \textstyle \frac{1}{\sqrt2} 
\|\sigma_D\|_\infty
\end{equation}
(see \cite[Lemma~2.4]{Koh-Tem}).

\medskip

\noindent {\bf Stress-strain duality.}
Let 
$$
\Sigma(\Om):=\{ \sigma\in L^2(\Om;\Mnn):
\div\,\sigma\in L^n(\Om;\Rn), \, \sigma_D\in L^\infty(\Om;\MD) \}\,.
$$
If $\sigma\in \Sigma(\Om)$,
then $\sigma\in L^r(\Om;\Mnn)$ for every $r<\infty$ by 
\cite[Proposition~2.5]{Koh-Tem}. For every $u\in BD(\Om)$ with 
$\div\,u\in L^{n/(n-1)}(\Om)$ we define the distribution 
$[\sigma_D{\,:\,}E_Du]$ on $\Om$ by
\begin{equation}\label{sigmaDED}
\langle [\sigma_D{\,:\,}E_Du] |\varphi\rangle:= -\langle 
\div\,\sigma|\varphi\, u\rangle
-{\textstyle \frac1n}\langle \tr\,\sigma|\varphi\, \div\,u\rangle -\langle\sigma| u 
{\,\odot\,} \nabla\varphi\rangle
\end{equation}
for every $\varphi\in C^\infty_c(\Om)$.
It is proved in \cite[Theorem~3.2]{Koh-Tem} that $[\sigma_D{\,:\,}E_Du]$ 
is a bounded measure on $\Om$ whose variation satisfies
\begin{equation}\label{varsigmaE}
|[\sigma_D{\,:\,}E_Du]|\le \|\sigma_D\|_\infty |E_D u| \qquad\hbox{in 
}\Om\,.
\end{equation}
Moreover 
\begin{equation}\label{psisigma}
[\psi\sigma_D{\,:\,}E_Du]= \psi[\sigma_D{\,:\,}E_Du] \qquad\hbox{in 
}\Om
\end{equation}
for every $\psi\in C^1(\ol\Om)$, and
\begin{equation}\label{sigmaEa}
[\sigma_D{\,:\,}E_Du]^a=\sigma_D{\,:\,}E^a_Du \qquad\hbox{a.e.\ in }
\Om
\end{equation}
(see \cite[Corollary~3.2]{Anz}). We define the measure 
$[\sigma_D{\,:\,}E^s_Du]$ on $\Om$ by
\begin{equation}\label{sigmaEsu}
[\sigma_D{\,:\,}E^s_Du]:= 
[\sigma_D{\,:\,}E_Du]^s=[\sigma_D{\,:\,}E_Du]-\sigma_D{\,:\,}E^a_Du\,.
\end{equation}
By (\ref{varsigmaE}) we have
\begin{equation}\label{varsigmaEsu}
|[\sigma_D{\,:\,}E^s_Du]|\le \|\sigma_D\|_\infty |E^s_D u| \qquad\hbox{in 
}\Om\,.
\end{equation}
This shows, in particular, that if $\hat\sigma$, $\hat u$ satisfy the same 
properties as $\sigma$, $u$,
and $\sigma_D=\hat \sigma_D$ a.e.\ on $\Om$, $E^s_D u=E^s_D \hat u$ in 
$\Om$, then 
$[\sigma_D{\,:\,}E^s_Du]=[\hat\sigma_D{\,:\,}E^s_D\hat u]$ in $\Om$. 

We define
\begin{equation} \label{<sE>}
\langle \sigma_D|E_Du \rangle:= [\sigma_D{\,:\,}E_Du](\Om)\,,\qquad
\langle \sigma_D|E^s_Du \rangle:= [\sigma_D{\,:\,}E^s_Du](\Om)\,,
\end{equation}
so that $\langle \sigma_D|E_Du \rangle=\langle \sigma_D|E^a_Du 
\rangle+\langle \sigma_D|E^s_Du \rangle$. 
If $\sigma_k\wto\sigma$ weakly in $L^2(\Om;\Mnn)$, 
$\div\,\sigma_k\wto\div\,\sigma$ weakly in $L^n(\Om;\Rn)$, 
and $(\sigma_k)_D$ is bounded in $L^\infty(\Om;\MD)$, then 
$\sigma_k\wto\sigma$ weakly in $L^r(\Om;\Mnn)$
for every $r<+\infty$ (see  \cite[Proposition~2.5]{Koh-Tem}) and 
\begin{equation}\label{convsE}
\langle [(\sigma_k)_D{\,:\,} E_Du] | \varphi \rangle \to 
\langle [\sigma_D{\,:\,} E_Du] | \varphi \rangle\,,
\qquad
\langle [(\sigma_k)_D{\,:\,} E^s_Du] | \varphi \rangle \to 
\langle [\sigma_D{\,:\,} E^s_Du] | \varphi \rangle
\end{equation}
for every $\varphi\in C(\overline\Om)$
(see \cite[Theorem~3.2]{Koh-Tem}, whose proof gives the result also in the case of weak convergence).

We define now a duality between $\Sigma(\Om)$ and $\Pi_{\Ga_0}(\Om)$.
Given $\sigma\in\Sigma(\Om)$ and $p\in\Pi_{\Ga_0}(\Om)$, we fix
$u\in BD(\Om)$, $e\in L^2(\Om;\Mnn)$, and $w\in H^1(\Rn;\Rn)$
satisfying (\ref{900}) and (\ref{901}).
Then we define a measure $[\sigma_D{\,:\,}p]\in M_b(\Om\cup\Ga_0)$ by setting
\begin{eqnarray*}
& [\sigma_D{\,:\,}p]:=\sigma_D{\,:\,}p^a+[\sigma_D{\,:\,}E_D^su]=
[\sigma_D{\,:\,}E_Du] - \sigma_D{\,:\,}e_D
 \quad\hbox{ on } \Om\,,
\\
& [\sigma_D{\,:\,}p]:= [\sigma \nu]_\nu^\perp {\,\cdot\,}(w-u)\, \hn\quad \hbox{ on } \Ga_0\,,
\end{eqnarray*}
so that
\begin{equation}\label{[sigmap]}
\langle [\sigma_D{\,:\,}p] |\varphi \rangle=
\langle [\sigma_D{\,:\,}E_Du] |\varphi \rangle -
\langle \sigma_D{\,:\,}e_D |\varphi \rangle +
\langle  [\sigma \nu]_\nu^\perp | \varphi (w-u)\rangle_{\Ga_0}
\end{equation}
for every $\varphi\in C(\overline\Om)$, where $\langle \cdot | \cdot \rangle_{\Ga_0}$ denotes the duality pairing between $L^\infty(\Ga_0;\Rn)$ and $L^1(\Ga_0;\Rn)$.
Using the previous remarks, it is easy to see that the measure $[\sigma_D{\,:\,}p]$ 
does not depend on the choice of $u$, $e$, and $w$. It follows from the definition and from (\ref{sigmatau}) and (\ref{varsigmaEsu}) that 
\begin{equation}\label{varsigmap}
\begin{array}{cc}
 [\sigma_D{\,:\,}p]^a= \sigma_D{\,:\,}p^a \quad\hbox{a.e.\ on }\Om \,,
 &
  [\sigma_D{\,:\,}p]^s= [\sigma_D{\,:\,}E_D^su]
\qquad\hbox{on }\Om\cup\Ga_0\,,
\vspace{1mm}
\\
 | [\sigma_D{\,:\,}p]| \le \|\sigma_D\|_\infty |p|
\quad\hbox{on }\Om\cup\Ga_0\,,
&
| [\sigma_D{\,:\,}p]^s| \le \|\sigma_D\|_\infty |p^s|
\quad\hbox{on }\Om\cup\Ga_0\,.
\end{array}
\end{equation}
Moreover (\ref{psisigma}) implies that
\begin{equation}\label{psisigmap}
 [\psi \sigma_D{\,:\,}p] =\psi [ \sigma_D{\,:\,}p]\qquad\hbox{in 
}\Om\cup\Ga_0
\end{equation}
for every $\psi\in C^1(\ol\Om)$.
Using the definitions we can deduce that
\begin{equation}\label{int}
\langle[\sigma_D{\,:\,}p]|\varphi \rangle=
\langle \varphi\, \sigma_D|p \rangle
\end{equation}
for every $\sigma\in C^1(\ol\Om;\Mnn)$ and every $\varphi\in C^1(\ol\Om)$, where the duality used in the right-hand side is defined in~(\ref{dualmeas}). Using the continuity properties given by (\ref{varsigmap}) we can prove by approximation that (\ref{int}) holds 
also for every  $\sigma\in C(\ol\Om;\Mnn)$ and every $\varphi\in C(\ol\Om)$. Therefore, for every $\sigma\in C(\ol\Om;\Mnn)$ and every $p\in\Pi_{\Ga_0}(\Om)$ we have
\begin{equation}\label{sigmaDp}
[\sigma_D {\,:\,}p]=\sigma_D {\,:\,}p \qquad \hbox{on }\,  \Om\cup\Ga_0\,,
\end{equation}
where the right-hand side denotes the measure defined by
\begin{equation}\label{sigmaDp2}
(\sigma_D {\,:\,}p)(B):=\int_B \sigma_D {\,:\,}dp:= \sum_{ij}\int_B  \sigma_{ij}\, dp_{ij}
\end{equation}
for every Borel set $B\subset \Om\cup\Ga_0$.

If  $\sigma_k\wto \sigma$ weakly in 
$L^2(\Om;\Mnn)$, $\div\, \sigma_k\wto \div\, \sigma$ weakly in 
$L^n(\Om;\Rn)$, and $(\sigma_k)_D$ is bounded in $L^\infty(\Om;\MD)$, then, 
using (\ref{sigmanu})--(\ref{sigmatau}) and (\ref{convsE}), we obtain
\begin{equation}\label{sigmakp}
\langle[(\sigma_k)_D {\,:\,}p]|\varphi \rangle\to \langle [\sigma_D {\,:\,}p]|\varphi\rangle
\end{equation}
for every $\varphi\in C(\ol\Om)$.

{}Finally, for every  $\sigma\in\Sigma(\Om)$ and $p\in\Pi_{\Ga_0}(\Om)$, 
we define
\begin{equation}\label{sigmap}
\begin{array}{c}
\langle\sigma_D |p\rangle:=[\sigma_D{\,:\,}p](\Om\cup\Ga_0)
=\langle\sigma_D|p^a\rangle +\langle\sigma_D|E_D^su\rangle 
+\langle [\sigma \nu]_\nu^\perp| w-u \rangle_{\Ga_0}=
\vspace{1mm}
\\
=\langle\sigma_D|E_D u\rangle - \langle\sigma_D|e_D \rangle
+\langle [\sigma \nu]_\nu^\perp| w-u \rangle_{\Ga_0}
\,,
\end{array}
\end{equation}
where
$u\in BD(\Om)$, $e\in L^2(\Om;\Mnn)$, and $w\in H^1(\Rn;\Rn)$
satisfy (\ref{900}) and (\ref{901}).

We are now in a position to prove an integration by parts formula for stresses $\sigma\in \Sigma(\Om)$ and displacements $u\in BD(\Om)$, involving the elastic and plastic strains $e$ and~$p$.

\begin{proposition}[Integration by parts]\label{Intparts}
Let $\sigma\in \Sigma(\Om)$, $f\in L^n(\Om;\Rn)$,  $g\in
L^{\infty}(\Ga_1;\Rn)$, and let $(u,e,p)\in A(w)$, with $w\in H^1(\Rn;\Rn)$. Assume that $-\div\,\sigma=f$ a.e.\ on $\Om$ and $[\sigma \nu]=g$ on $\Ga_1$. Then
\begin{equation}\label{intpartsp2}
\langle \sigma_D | p \rangle + 
 \langle \sigma | e-Ew \rangle
=\langle f | u-w \rangle +
\langle g|u-w\rangle_{\Ga_1} \,,
\vspace{1mm}
\end{equation}
where $\langle \cdot | \cdot \rangle_{\Ga_1}$ denotes the duality pairing between $L^\infty(\Ga_1;\Rn)$ and $L^1(\Ga_1;\Rn)$. Moreover
\begin{equation}\label{intpartsp}
\begin{array}{c}  
\langle [\sigma_D{\,:\,}p] |\varphi \rangle + 
 \langle \sigma{\,:\,} (e-Ew)|\varphi \rangle +
\langle\sigma|(u-w){\,\odot\,}\nabla\varphi \rangle
= \vspace{1mm}
\\
= \langle f | \varphi (u-w) \rangle +
\langle g|\varphi (u-w)\rangle_{\Ga_1}
\end{array}
\end{equation}
for every $\varphi\in C^1(\overline\Om)$.
\end{proposition}

\begin{proof}
By \cite[Theorem~3.2 and Propositions~3.3 and~3.4]{Koh-Tem} we have
\begin{equation}\label{intparts}
\begin{array}{c}
\langle \div\,\sigma|\varphi\, v\rangle + \langle 
[\sigma_D{\,:\,}E_Dv] |\varphi \rangle + 
\frac{1}{n} \langle\tr\,\sigma|\varphi \,\div\,v\rangle +
\langle\sigma|v{\,\odot\,}\nabla\varphi \rangle=
\vspace{1mm}
\\
= \langle[\sigma\nu]_\nu^\perp| 
\varphi \,v\rangle_{\Ga_0} +
\langle g|\varphi \,v\rangle_{\Ga_1} 
\end{array}
\end{equation}
for every $\varphi\in C^1(\overline\Om)$ and every  $v\in BD(\Om)$ with 
$\div\, v\in L^2(\Om)$ and 
$v{\,\cdot\,}\nu=0$ $\hn$-a.e.\ on~$\Ga_0$. By (\ref{[sigmap]})  we have
\begin{equation}\label{intpartsp3}
\begin{array}{c}  
\langle [\sigma_D{\,:\,}p] |\varphi \rangle + 
 \langle \sigma{\,:\,} (e-Ew)|\varphi \rangle +
\langle\sigma|(u-w){\,\odot\,}\nabla\varphi \rangle
= 
\vspace{1mm}
\\
= \langle [\sigma_D{\,:\,}E_D(u-w)] |\varphi \rangle +
\frac{1}{n} \langle\tr\,\sigma|\varphi \,\div(u-w)\rangle +
\langle\sigma|(u-w){\,\odot\,}\nabla\varphi \rangle-{}
\vspace{1mm}
\\
{}-\langle  [\sigma \nu]_\nu^\perp | \varphi (u-w)\rangle_{\Ga_0}\,.
\end{array}
\end{equation}
If we apply (\ref{intparts}) with $v=u-w$ we obtain
\begin{equation}\label{intpartsp4}
\begin{array}{c}  
\langle [\sigma_D{\,:\,}E_D(u-w)] |\varphi \rangle +
\frac{1}{n} \langle\tr\,\sigma|\varphi \,\div\,(u-w)\rangle +
\langle\sigma|(u-w){\,\odot\,}\nabla\varphi \rangle-{}
\vspace{1mm}
\\
{}-\langle  [\sigma \nu]_\nu^\perp | \varphi (u-w)\rangle_{\Ga_0}
=
 \langle f | \varphi (u-w) \rangle +
\langle g|\varphi (u-w)\rangle_{\Ga_1}
\,.
\end{array}
\end{equation}
Equality (\ref{intpartsp}) follows now from  (\ref{intpartsp3}) and  (\ref{intpartsp4}). To obtain  (\ref{intpartsp2}) it is enough to take $\varphi=1$ in~(\ref{intpartsp}).
\end{proof}

In order to show the connection between the duality (\ref{sigmap}) and the functional $\HH$ defined in (\ref{calHD}), 
we need the following approximation result.

\begin{lemma}\label{approxsigma}
Let $U$ be a bounded open set in $\Rn$ with the segment property, let 
$\KK $ be a closed convex subset of $\Mnn$, and 
let $\sigma\in L^r(U;\Mnn)$, $1\le r<+\infty$, with $\div\, \sigma\in 
L^r(U;\Rn)$ and $\sigma(x)\in \KK $ for a.e.\ $x\in U$. 
Then there exists a sequence $\sigma_k\in C^\infty(\overline U;\Mnn)$ 
such that $\sigma_k\to \sigma$ strongly in 
$L^r( U;\Mnn)$, $\div\, \sigma_k\to \div\, \sigma$ strongly in 
$L^r( U;\Rn)$, and $\sigma_k(x)\in \KK $ for every $x\in\overline U$.
\end{lemma}

\begin{proof}
Since $ U$ is bounded and has the segment property, there exists a finite 
open cover $(U_i)$, $i=1, \ldots , m$, 
of $\partial U$ and a corresponding sequence of nonzero vectors $y_i$ 
such that, if $x\in \overline U\cap U_i$ for some $i$, 
then $x+ty_i\in  U$ for $0<t<1$. We set $U_0:= U$ and $y_0:=0$. For 
$i=0, \ldots , m$ and $k=1, 2, \ldots$ the open set 
$U^i_k:=\{x\in U_i: x+(1/k)y_i\in  U\}$ contains $\overline U\cap U_i$.
We define $\sigma^i_k(x):=\sigma(x+(1/k)y_i)$ for every $x\in U^i_k$. Let 
$(V_i)$, $i=0, \ldots , m$, be an open cover 
of $\overline U$ such that $V_i\subset\subset U_i$ for every~$i$. Since $ 
\overline U\cap \overline V_i \subset U^i_k$,
for every $i$ and $k$ we can find a mollifier $\psi^i_k$ of class 
$C^\infty_c(\Rn)$ such that the convolution 
$\sigma^i_k\star \psi^i_k$ is well defined in a neighbourhood of 
$ \overline U\cap \overline V_i$ and
\begin{equation}\label{1/k}
\|\sigma^i_k\star \psi^i_k-\sigma^i_k\|_{r, U\cap V_i}\le \frac1k\qquad 
\hbox{and}
\qquad \|\div\, \sigma^i_k\star \psi^i_k- \div\,\sigma^i_k\|_{r, U\cap 
V_i}\le \frac1k \,.
\end{equation}
As $\KK $ is closed and convex, we have $\sigma^i_k\star \psi^i_k(x)\in \KK $ 
for every $x$ in a neighbourhood of 
$ \overline U\cap \overline V_i$.

Let $(\varphi_i)$, $i=0, \ldots , m$, be a $C^\infty$ partition of unity 
for $\overline  U$ subordinate to $(V_i)$ and let
$$
\sigma_k:=\sum_{i=0}^m \varphi_i (\sigma^i_k\star \psi^i_k)\,.
$$
Then $\sigma_k$ is of class $C^\infty$ in a neighbourhood of 
$\overline U$ and $\sigma_k(x)\in \KK $ for every $x$ in a 
neighbourhood of $\overline U$. Since $\sigma^i_k\to \sigma$ strongly in 
$L^r( U\cap V_i;\Mnn)$ and  
$\div\, \sigma^i_k\to \div\, \sigma$ strongly in $L^q( U\cap V_i;\Rn)$, 
from (\ref{1/k}) and from the identity
$$
\div\,\sigma:=\sum_{i=0}^m (\varphi_i \,\div\,\sigma + \sigma\,\nabla 
\varphi_i)
$$
we deduce that $\sigma_k\to \sigma$ strongly in 
$L^r( U;\Mnn)$ and $\div\, \sigma_k\to \div\, \sigma$ strongly in 
$L^r( U;\Rn)$.
\end{proof}

The following proposition provides a variant of (\ref{HHduality}) expressed by using the 
duality~(\ref{sigmap}).

\begin{proposition}\label{prp:star}
Let $p\in \Pi_{\Ga_0}(\Om)$. Then
\begin{equation}\label{Hpge}
H(p)\ge [\sigma_D : p] \quad\hbox{on }\,\Om\cup\Ga_0
\end{equation}
for every $\sigma\in\Sigma(\Om)\cap \K(\Om)$, and
\begin{equation}\label{Hpsup}
\HH(p)=\sup \{ \langle\sigma_D | p\rangle:  \sigma\in\Sigma(\Om)\cap \K(\Om)\} \,.
\end{equation}
Moreover, if $g\in L^\infty(\Ga_1;\Rn)$ and there exists $\varrho\in \Sigma(\Om)\cap \K(\Om)$ such that $[\varrho\nu]=g$ on $\Ga_1$, then
\begin{equation}\label{Hpsupg}
\HH(p)=\sup \{ \langle\sigma_D | p\rangle:  \sigma\in\Sigma(\Om)\cap \K(\Om),\ 
[\sigma\nu]=g \, \hbox{ on }\,\Ga_1 \}\,.
\end{equation}
\end{proposition}

\begin{proof}
Let $\sigma\in\Sigma(\Om)\cap\K(\Om)$. To prove (\ref{Hpge}) it is enough to show that
\begin{equation}\label{HHp}
\langle H(p)|\varphi\rangle \ge \langle [\sigma_D {\,:\,} p] | \varphi \rangle 
\end{equation}
for every $\varphi\in C(\ol\Om)$ with $\varphi\ge 0$ on $\ol\Om$.
By Lemma~\ref{approxsigma} there exists a sequence
$(\sigma_k)$ in $C^{\infty}(\ol\Om;\Mnn)\cap\K(\Om)$ such that $\sigma_k\to \sigma$ strongly in 
$L^n(\Om;\Mnn)$ and $\div\, \sigma_k\to \div\, \sigma$ strongly in 
$L^n(\Om;\Rn)$.
By (\ref{HD}), (\ref{Hmu}), and (\ref{int}) we have
$$
\langle H(p)|\varphi\rangle \ge \langle [(\sigma_k)_D {\,:\,} p] | \varphi \rangle \,,
$$
and (\ref{HHp}) follows from (\ref{sigmakp}). This concludes the proof of~(\ref{Hpge}).

By \cite[Chapter~II, Section~4]{Tem} we have
$$
\HH(p)=\,\sup \{\langle \sigma_D|p\rangle :
\sigma\in C^\infty(\Rn;\Mnn)\cap\K(\Om),\ \supp\,\sigma\cap\Ga_1=\emptyset\}
\,.
$$
This equality, together with (\ref{int}) and~(\ref{Hpge}), implies (\ref{Hpsup}) and (\ref{Hpsupg}) with $g=0$. 

Let $\phi\in C^\infty(\R)$ be such that $0\le \phi\le 1$, $\phi(s)=0$ for $s\le 1$, and $\phi(s)=1$
for $s\geq 2$. For $\delta>0$ we consider the function $\psi_\delta(x):=\phi(\frac{1}{\delta}\dist(x,\Ga_1))$ 
defined for every $x\in\ol\Om$.
Let $\sigma\in\Sigma(\Om)\cap\K(\Om)$ be such that $[\sigma\nu]=0$ on $\Ga_1$. Then $\sigma_\delta:=\psi_\delta\sigma+(1-\psi_\delta)
\varrho\in\Sigma(\Om)\cap\K(\Om)$ and $[\sigma_\delta\nu]=g$ on $\Ga_1$. 
Moreover, by (\ref{psisigmap}) we have
$$
\langle(\sigma_\delta)_D|p\rangle=\langle[\sigma_D{\,:\,}p]|\psi_\delta\rangle
+\langle[\varrho_D{\,:\,}p]|1-\psi_\delta\rangle\,.
$$
Since the right-hand side converges to $\langle\sigma_D|p\rangle$ as $\delta\to 0$,
equality (\ref{Hpsupg}) follows from the equality already proved for $g=0$ and from~(\ref{Hpge}).
\end{proof}

\end{section}

\begin{section}{The minimum problem}

In this section we study in detail the minimum problem used in the incremental formulation of the quasistatic evolution. 
The data are the current value $p_0\in \Pi_{\Ga_0}(\Om)$ of the plastic strain and the updated values 
$w\in  H^1(\Rn;\Rn)$, $f\in L^n(\Om;\Rn)$, and $g\in L^{\infty}(\Ga_1;\Rn)$ of the boundary displacement and of the body and surface loads. The total load $\LL\in BD(\Om)'$ is defined by
\begin{equation}\label{LL}
\langle \LL|u\rangle:= \langle f|u\rangle+\langle g|u\rangle_{\Ga_1}
\end{equation}
for every $u\in BD(\Om)$.
By solving the minimum problem
\begin{equation}\label{minp0}
\min_{(u,e,p)\in A(w)}  \{\QQ(e)+\HH(p-p_0)-\langle \LL|u\rangle\}
\end{equation}
we get the updated values $u$, $e$, and $p$ of displacement, elastic and plastic strain.

For the existence result we will assume the following safe-load condition:
there exist  $\varrho\in L^2(\Om;\Mnn)$ and $\alpha >0$ such that
\begin{equation}\label{sl10}
- \div\varrho=f\, \hbox{ a.e.\ on } \Om\,, \qquad [\varrho\nu]=g \hbox{ on } \Ga_1\,,
\end{equation}
and
\begin{equation}\label{sl11}
\varrho_D(x)+\xi\in K\, \hbox{  } 
\end{equation}
for a.e.\  $x\in \Om$ and for every $\xi\in\MD$ with $|\xi|\le\alpha$.

\subsection{Existence of a minimizer}
We begin by proving two technical lemmas concerning the safe-load condition.

\begin{lemma}\label{fg}
Let $w\in H^1(\Rn;\Rn)$, $f\in L^n(\Om;\Rn)$, $g\in L^{\infty}(\Ga_1;\Rn)$, and let $\LL$ be defined by~(\ref{LL}).
Assume (\ref{sl10}) and (\ref{sl11}). Then
\begin{equation}\label{fg1}
\langle \LL|u\rangle =  \langle \varrho|e\rangle
+\langle \varrho_D|p\rangle -  \langle \varrho|Ew\rangle + \langle \LL|w\rangle
\end{equation}
for every $(u,e,p)\in A(w)$.
\end{lemma}

\begin{proof} 
The result follows from the definition (\ref{sigmap}) of the duality product $\langle \varrho_D|p\rangle$ 
and from the integration by parts formula~(\ref{intpartsp2}).
\end{proof}

\begin{lemma}\label{eps*}
Let $f\in L^n(\Om;\Rn)$, $g\in L^{\infty}(\Ga_1;\Rn)$, $\varrho\in  L^2(\Om;\Mnn)$, and $\alpha>0$.
Assume (\ref{sl10}) and (\ref{sl11}). Then
\begin{equation}\label{eps*1}
\HH(p)-\langle \varrho_D|p\rangle \ge \alpha \|p\|_1
\end{equation}
for every $p\in \Pi_{\Ga_0}(\Om)$.
\end{lemma}

\begin{proof} 
By Proposition~\ref{prp:star} we have
\begin{eqnarray*}
&\HH(p) -\langle \varrho_D|p\rangle =\sup
\{\langle\sigma_D -\varrho_D| p\rangle: \sigma \in  \Sigma(\Om) \cap \K(\Om)\} 
\ge
\\
&\ge \sup
\{\langle\tau_D| p\rangle: \tau \in  \Sigma(\Om),\ \|\tau_D\|_\infty\le \alpha\}\,.
\end{eqnarray*}
{}From (\ref{int}) 
it follows that 
$$
\HH(p) -\langle \varrho_D|p\rangle \ge \sup
\{\langle\tau_D| p\rangle : \tau \in   C^\infty(\ol\Om;\Mnn), \  \|\tau_D\|_\infty\le \alpha\}\,,
$$
where the duality product in the right-hand side is defined by~(\ref{dualmeas}).
The conclusion follows now from standard arguments in measure theory.
\end{proof}

We are now in a position to prove the existence of a solution to~(\ref{minp0}).

\begin{theorem}\label{existencemin}
Let $w\in H^1(\Rn;\Rn)$, $p_0\in \Pi_{\Ga_0}(\Om)$, $f\in L^n(\Om;\Rn)$, $g\in L^{\infty}(\Ga_1;\Rn)$, and let $\LL$ be defined by~(\ref{LL}).
Assume (\ref{sl10}) and (\ref{sl11}). Then the minimum problem (\ref{minp0}) has a solution.
\end{theorem}

\begin{proof} 
By Lemma~\ref{fg} the minimum problem (\ref{minp0}) is equivalent to  
\begin{equation}\label{minp1}
\min_{(u,e,p)\in A(w)}  \{\QQ(e) -\langle \varrho|e\rangle +\HH(p-p_0)
-\langle \varrho_D|p-p_0\rangle\}\,,
\end{equation}
in the sense that these problems have the same solutions.
Let $(u_k,e_k,p_k)\in A(w)$ be a minimizing sequence. By Lemma~\ref{eps*} we have
$$
\HH(p_k-p_0) -\langle \varrho_D|p_k-p_0\rangle\ge \alpha \|p_k-p_0\|_1\,,
$$
while (\ref{boundsC}) gives
$$
\QQ(e_k)-\langle\varrho|e_k\rangle\ge \frac{\alpha_\C}{2}\|e_k\|_2^2-\frac{1}{2\alpha_\C}\|\varrho\|_2^2\,.
$$
Therefore, the sequences $e_k$ and $p_k$ are bounded in $L^2(\Om;\Mnn)$
and in $M_b(\Om\cup\Ga_0;\MD)$, respectively.
Since $Eu_k=e_k+p_k$ in $\Om$, it follows that $Eu_k$ is bounded in $M_b(\Om;\Mnn)$. 
Since $(w-u_k){\,\odot\,}\nu\,\hn=p_k$ is bounded in $M_b(\Ga_0;\MD)$, the traces of $u_k$ are bounded 
in $L^1(\Ga_0;\Rn)$. Therefore $u_k$ is bounded in $BD(\Om)$ by (\ref{seminorm}). 
Up to extracting a subsequence, we may assume that $u_k\wto u$ weakly$^*$ in 
$BD(\Om)$, $e_k\wto e$ weakly in $L^2(\Om;\Mnn)$,
$p_k\wto p$ weakly$^*$ in $M_b(\Om\cup\Ga_0;\MD)$. By Lemma~\ref{lm:145} we have
$(u,e,p)\in A(w)$. By lower semicontinuity 
\begin{equation}\label{0part}
\QQ(e)-\langle \varrho|e\rangle \le\liminf_{k\to \infty}
\{\QQ(e_k)-\langle\varrho|e_k\rangle\} \,.
\end{equation}
To conclude we just need to show that
\begin{equation}\label{2part}
\HH(p-p_0) - \langle \varrho_D|p-p_0\rangle\le
\liminf_{k\to\infty}\{\HH(p_k-p_0) -\langle \varrho_D|p_k-p_0\rangle \}
\,.
\end{equation}
To this aim, let $\phi\in C^\infty(\R)$ be such that $0\le \phi\le 1$, $\phi(s)=0$ for $s\le 1$, and
$\phi(s)=1$ for $s\ge 2$. Let $\delta>0$ and $\psi_\delta(x):=\phi(\frac{1}{\delta}\dist(x,\Gamma_1))$
for every $x\in\ol\Om$. Since the measure $H(p_k-p_0)-[\varrho_D{\,:\,}(p_k-p_0)]$ is nonnegative on 
$\Om\cup\Ga_0$ by (\ref{Hpge}), we have
\begin{equation}\label{partlsc}
\HH(\psi_\delta(p_k-p_0)) -\langle [\varrho_D{\,:\,}(p_k-p_0)]|\psi_\delta\rangle
\le \HH(p_k-p_0) -\langle \varrho_D|p_k-p_0\rangle
\end{equation}
for every $\delta>0$.
The integration by parts formula (\ref{intpartsp}) gives
\begin{eqnarray*}
& \langle [\varrho_D{\,:\,}(p_k-p_0)]|\psi_\delta\rangle = 
 -\langle \varrho{\,:\,}(e_k- E w)|\psi_\delta\rangle -
 \langle \varrho| (u_k-w){\,\odot\,}\nabla\psi_\delta \rangle +{}&
\\
& {}+
\langle f|\psi_\delta(u_k-w)\rangle
- \langle [\varrho_D{\,:\,}p_0]|\psi_\delta\rangle\,. &
\end{eqnarray*}
Passing to the limit as $k\to\infty$, and using  (\ref{intpartsp}) again,  we deduce that
\begin{equation}\label{1part}
\langle [\varrho_D{\,:\,}(p-p_0)]|\psi_\delta\rangle=
\lim_{k\to\infty}\langle [\varrho_D{\,:\,}(p_k-p_0)]|\psi_\delta\rangle\,.
\end{equation}
By (\ref{partlsc}), (\ref{1part}), and the lower semicontinuity of $\HH$, we have
$$
\HH(\psi_\delta(p-p_0)) - \langle [\varrho_D{\,:\,}(p-p_0)]|\psi_\delta\rangle
\le \liminf_{k\to\infty}\{\HH(p_k-p_0) -\langle \varrho_D|p_k-p_0\rangle \}\,.
$$
Passing to the limit as $\delta\to 0$ we finally obtain (\ref{2part}).

As $(u_k,e_k,p_k)$ is a minimizing sequence and $(u,e,p)\in A(w)$, by (\ref{0part}) and (\ref{2part})
we conclude that $(u,e,p)$ is a minimizer of~(\ref{minp1}).
\end{proof}

\subsection{The Euler conditions}
We now derive the Euler conditions for a minimizer of (\ref{minp0}) in the 
special case $p=p_0$.

\begin{theorem}\label{Euler}
Let $w\in H^1(\Rn;\Rn)$, $f\in L^n(\Om;\Rn)$, $g\in L^{\infty}(\Ga_1;\Rn)$, and let $\LL$ be defined by~(\ref{LL}).
Suppose that $(u,e,p)$ is a solution of (\ref{minp0}) with $p_0=p$, and let 
$\sigma:=\C e$. Then $\sigma\in L^2(\Om;\Mnn)$ and
\begin{equation}\label{Eulereq}
\textstyle
-\HH(q)\le  \langle \sigma|\eta\rangle  -   
 \langle \LL| v\rangle
= \langle \sigma_D|\eta_D \rangle +
 \frac1n \langle \tr\,\sigma|\div\, v\rangle -   
 \langle \LL| v\rangle
\le \HH(-q)
\end{equation}
for every $(v,\eta,q)\in A(0)$.
\end{theorem}

\begin{proof}
Let us fix $(v,\eta,q)\in A(0)$. For every $\e\in \R$ the triple $(u+\e v, e+\e \eta, p+\e 
q)$ belongs to $A(w)$, and hence
\begin{equation*}
\QQ(e+\e\eta) + \HH(\e q) - \e \langle \LL| v\rangle \ge \QQ(e)\qquad 
\hbox{for every }\e\in\R\,.
\end{equation*}
Using the positive homogeneity of $H$ we obtain
$$
\QQ(e\pm \e\eta) + \e\HH(\pm q)
\mp \e \langle \LL| v\rangle \ge \QQ(e)\qquad\hbox{for every }\e>0 
\,.
$$
Taking the derivative with respect to $\e$ at $\e=0$, we get
$$
\langle \sigma|\eta\rangle +\HH(q) - \langle \LL| v\rangle \ge0\,,
\qquad
-\langle \sigma|\eta\rangle +\HH(-q)+  \langle \LL| v\rangle  \ge0\,,
$$
which implies (\ref{Eulereq}).
\end{proof}

\begin{proposition}\label{eulersigma}
Let $\sigma\in L^2(\Om;\Mnn)$, $f\in L^n(\Om;\Rn)$, $g\in L^{\infty}(\Ga_1;\Rn)$, and let $\LL$ be defined by~(\ref{LL}). The following conditions are equivalent:
\begin{itemize}
\item[(a)]$-\HH(q)\le  \langle \sigma|\eta\rangle  -   
 \langle \LL| v\rangle
\le \HH(-q)$
 for every 
$(v,\eta,q)\in A(0)$;
\smallskip
\item[(b)] $\sigma\in\Sigma(\Om)\cap \K(\Om)$,  $-\div\,\sigma=f$ a.e.\ on $\Om$, and $[\sigma \nu]=g$ on~$\Ga_1$.
\end{itemize}
\end{proposition}

\begin{proof} Assume (a) and let $v\in H^1(\Om;\Rn)$ with $v=0$ $\hn$-a.e.\ on $\Ga_0$. Since the triple
$(v,Ev,0)$ belongs to $A(0)$, from (a) we obtain
\begin{equation}\label{weakdiv}
\langle \sigma|E v \rangle-\langle f|v\rangle-\langle g|v\rangle_{\Ga_1}=0\,.
\end{equation}
Since this is true, in particular, for $v\in C^\infty_c(\Om;\Rn)$, we 
conclude that $-\div\,\sigma=f$ on $\Om$, hence $\div\,\sigma\in L^n(\Om;\Rn)$. Using the distributional definition (\ref{sigmanu}) of $[\sigma \nu]$, from (\ref{weakdiv}) we 
obtain also that $[\sigma \nu]=g$ on $\Ga_1$.

Let $\eta\in L^2(\Om;\MD)$. Regarding $-\eta$ as an absolutely continuous measure on $\Om\cup\Ga_0$,  the triple $(0,\eta,-\eta)$ belongs to $A(0)$, thus from 
(a) we obtain
$$
-\HH(-\eta)\le \langle \sigma_D|\eta \rangle\le \HH(\eta)\,.
$$
Let us fix $\xi\in\MD$. Since for every Borel set $B\subset\Om$ we can 
take $\eta(x)=1_B(x)\,\xi$, we deduce that
$$
-H(-\xi)\le \sigma_D(x){\,:\,}\xi \le H(\xi)\qquad
\hbox{for a.e.\ }x\in\Om\,.
$$
Therefore $\sigma_D(x)\in\partial H(0)$ for a.e.\ $x\in\Om$. As 
$\partial H(0)=K$ 
(see, e.g., \cite[Corollary~23.5.3]{Roc}), we obtain that $\sigma_D(x)\in 
K$ for a.e.\ $x\in\Om$, hence $\sigma_D\in L^\infty(\Om;\MD)$ and $\sigma\in \K(\Om)$.

Conversely, assume (b) and  let $(v,\eta,q)\in A(0)$. By Proposition~\ref{prp:star} we have
\begin{equation}\label{wpq}
-\HH(-q)\le \langle\sigma_D |q\rangle\le \HH(q)\,.
\end{equation}
{}From the integration by parts formula~(\ref{intpartsp2})
we get
$$
\langle\sigma_D |q\rangle  =  -\langle\sigma |\eta\rangle 
+\langle f|v\rangle +\langle g|v\rangle_{\Ga_1}
\,,
$$
so that (a) follows now from (\ref{wpq}).
\end{proof}

\begin{theorem}\label{cns}
Let $w\in H^1(\Rn;\Rn)$, $f\in L^n(\Om;\Rn)$, $g\in L^\infty(\Om;\Rn)$,
let $(u,e,p)\in A(w)$, let $\sigma:= \C e$, and let $\LL$ be defined by~(\ref{LL}). Then the following conditions are equivalent:
\begin{itemize}
\item[(a)]$(u,e,p)$ is a solution of (\ref{minp0}) with $p_0=p$;
\smallskip
\item[(b)]$-\HH(q)\le  \langle \sigma|\eta\rangle  -   
 \langle \LL| v\rangle
\le \HH(-q)$
 for every 
$(v,\eta,q)\in A(0)$;
\smallskip
\item[(c)]$\sigma\in\Sigma(\Om)\cap \K(\Om)$,  $-\div\,\sigma=f$ a.e.\ on $\Om$, and $[\sigma \nu]=g$ on~$\Ga_1$.
\end{itemize}
\end{theorem}

\begin{proof} The implication $\rm{(a)}\Rightarrow\rm{(b)}$ is proved in 
Theorem~\ref{Euler}. 
The converse is true by convexity. The equivalence $\rm{(b)}\Leftrightarrow\rm{(c)}$ is proved 
in Proposition~\ref{eulersigma}. 
\end{proof}

Theorem~\ref{cns} gives immediately a stability result with respect to weak convergence of the data.

\begin{theorem}\label{lm:45}
Let $w_k$, $f_k$, $g_k$ be sequences in $H^1(\Rn;\Rn)$, $L^n(\Om;\Rn)$, $L^\infty(\Om;\Rn)$ respectively, let $\LL_k$ be defined by (\ref{LL}) with $f=f_k$ and $g=g_k$, and let
$(u_k,e_k,p_k)\in A(w_k)$. Assume that $u_k\wto u_\infty$ weakly$^*$ in 
$BD(\Om)$, $e_k\wto e_\infty$ weakly in $L^2(\Om;\Mnn)$,
$p_k\wto p_\infty$ weakly$^*$ in $M_b(\Om\cup\Ga_0;\MD)$, $w_k\wto w_\infty$ 
weakly in $H^1(\Rn;\Rn)$, $f_k\wto f_\infty$ weakly in $L^n(\Om;\Rn)$, $g_k\wto g_\infty$  weakly$^*$ in $L^\infty(\Om;\Rn)$, and  let $\LL_\infty$ be defined by (\ref{LL}) with $f=f_\infty$ and $g=g_\infty$. If 
\begin{equation}\label{closedstable}
\QQ(e_k) - \langle\LL_k|u_k\rangle \le \QQ(\eta)+\HH(q-p_k) - \langle\LL_k|v\rangle
\end{equation}
for every~$k$ and every $(v,\eta,q)\in A(w_k)$, then $(u_\infty,e_\infty,p_\infty)\in A(w_\infty)$ and
\begin{equation}\label{closedstable2}
\QQ(e_\infty) - \langle\LL_\infty|u_\infty\rangle\le \QQ(\eta)+\HH(q-p_\infty)- \langle\LL_\infty|v\rangle
\end{equation}
for every $(v,\eta,q)\in A(w_\infty)$.
\end{theorem}

\begin{proof} First we note that $(u_\infty,e_\infty,p_\infty)\in A(w_\infty)$ by 
Lemma~\ref{lm:145}. Let $\sigma_k:=\C e_k$ and  $\sigma_\infty:=\C e_\infty$.
If (\ref{closedstable}) holds, then $u_k$, $e_k$, $p_k$, $w_k$, $f_k$, $g_k$ satisfy condition (a) of Theorem~\ref{cns}. By condition (c) of Theorem~\ref{cns} we have $\sigma_k\in\Sigma(\Om)\cap \K(\Om)$,  $-\div\,\sigma_k=f_k$ a.e.\ on $\Om$, and $[\sigma_k \nu]=g_k$ on~$\Ga_1$.

Since $e_k\wto e_\infty$ weakly in $L^2(\Om;\Mnn)$, we have that $\sigma_k\wto \sigma_\infty$ weakly in $L^2(\Om;\Mnn)$. As $\K(\Om)$ is closed and convex in $L^2(\Om;\Mnn)$, we deduce that $\sigma_\infty\in \K(\Om)$. 
Since $-\div\,\sigma_k=f_k$ a.e.\ on $\Om$ and  $f_k\wto f_\infty$ weakly in $L^n(\Om;\Rn)$, we obtain that $-\div\,\sigma_\infty=f_\infty$ a.e.\ on $\Om$, hence $\sigma_\infty\in\Sigma(\Om)$. Moreover, from (\ref{sigmanu}) it follows that 
$[\sigma_k \nu]\wto [\sigma_\infty \nu]$ weakly in $H^{-1/2}(\partial\Om;\Rn)$.
As $[\sigma_k \nu]=g_k$ on~$\Ga_1$ and $g_k\wto g_\infty$  weakly$^*$ in $L^\infty(\Om;\Rn)$, we conclude that $[\sigma_\infty \nu]=g_\infty$ on~$\Ga_1$. Therefore $u_\infty$, $e_\infty$, $p_\infty$, $w_\infty$, $f_\infty$, $g_\infty$ satisfy condition (c) of Theorem~\ref{cns}. Inequality (\ref{closedstable2}) follows now from condition (a) of Theorem~\ref{cns}.
\end{proof}

\subsection{Continuous dependence on the data}
We complete our study of the solutions $(u,e,p)$ of the minimum problem (\ref{minp0}) in the special case $p=p_0$
by proving the continuous dependence, in the norm topology, of $u$ and $e$ on the data $p_0$, $w$, $f$, and~$g$.

\begin{theorem}\label{thmhoelder}
For $i=1,2$ let $w_i\in H^1(\Rn;\Rn)$, $f_i\in L^n(\Om;\R^n)$, $g_i\in L^\infty(\Ga_1;\R^n)$, and let $\LL_i$ be defined by (\ref{LL}) with $f=f_i$ and $g=g_i$.
Suppose that $(u_i,e_i,p_i)$ is a solution of (\ref{minp0}) with $p_0=p_i$, $w=w_i$,
$\LL=\LL_i$, and let 
$$
\omega_{12}:=\|p_2-p_1\|_1 +   \|p_2-p_1\|_1^{1/2} + \|f_2-f_1\|_n +
\|g_2-g_1\|_{\infty,\Ga_1} +  \|Ew_2-Ew_1\|_2\,.
$$
Then
\begin{eqnarray}
& \|e_2-e_1\|_2\le C_1\, \omega_{12}\,,
\label{eq12}
\\
& \|Eu_2-Eu_1\|_1\le C_2\, \omega_{12} \,,
\label{hoelder}
\\
& \|u_2-u_1\|_1\le C_3 \,(\omega_{12} + \|w_2-w_1\|_2)\,,
\label{hoelderu}
\end{eqnarray}
where $C_1$, $C_2$, and $C_3$ are positive constants depending only on
$R_K$, $\alpha_\C$, $\beta_\C$, $\Om$, and $\Ga_0$.
\end{theorem}

\begin{proof}
Let $v:=(u_2-w_2)-(u_1-w_1)$, $\eta:= (e_2-Ew_2)-(e_1-Ew_1)$, and 
$q:=p_2-p_1$. Since $(v,\eta,q)\in AP(0)$, by Theorem~\ref{Euler} we obtain
\begin{eqnarray*}
&-\HH(p_2-p_1)\le \langle \C e_1|\eta  \rangle - \langle f_1|v \rangle
-\langle g_1|v \rangle_{\Ga_1}\,,
\\
& \langle \C e_2|\eta  \rangle - \langle f_2|v \rangle
-\langle g_2|v \rangle_{\Ga_1}
\le \HH(p_1-p_2)\,.
\end{eqnarray*}
Adding term by term and using (\ref{boundsH}) we obtain
\begin{equation}\nonumber
\begin{array}{c}
 \langle \C (e_2-e_1)|e_2-e_1 \rangle \le  
 \langle \C  (e_2-e_1)|Ew_2-Ew_1 \rangle + {}
 \\ 
{}+\langle f_2-f_1|v\rangle + \langle  g_2-g_1|v \rangle_{\Ga_1}
+  2\,R_K \|p_2-p_1\|_1\,. 
\end{array}
\end{equation}
By (\ref{boundsC}) and (\ref{normC}) this implies
\begin{equation}\label{ineq0}
\begin{array}{c}
2\,\alpha_\C \|e_2-e_1\|_2^2\le 2\,\beta_\C\|e_2-e_1\|_2\, \|Ew_2-Ew_1\|_2 
+\|f_2-f_1\|_n\, \|v\|_{n/(n-1)} +{} \\
{}+\|g_2-g_1\|_{\infty ,\Ga_1} \|v\|_{1,\Ga_1}+ 2\,R_K \|p_2-p_1\|_1\,.
\end{array}
\end{equation}

Since the embedding of $BD(\Om)$ into $L^{n/(n-1)}(\Om;\Rn)$ is continuous, there exists a constant $A_1$, depending only on $\Om$, such that
\begin{equation}\label{n/n-1}
 \|v\|_{n/(n-1)} \le A_1\,  \|v\|_1 + A_1\, \|Ev\|_1\,.
 \end{equation}
By (\ref{seminorm}) there exists a constant  $C>0$, 
depending only on $\Om$ and $\Ga_0$, such that
\begin{equation}
\label{trace}
\|v\|_1 \le C\, \|v\|_{1,\Ga_0}+ 
C\, \|Ev\|_1 \,.
\end{equation}
As $p_2-p_1=-v{\,\odot\,}\nu\, \hn$ on $\Ga_0$, we have
\begin{equation}
\label{tracep}
\|v\|_{1,\Ga_0}\le \sqrt 2\, \|p_2-p_1\|_1\,.
\end{equation}
Since $Ev=(e_2-e_1) +(p_2-p_1) - (Ew_2-Ew_1)$, by the H\"older inequality we have also
\begin{equation}
\label{Ev}
\|Ev\|_1\le   \Ln(\Om)^{1/2} \|e_2-e_1\|_2 + \|p_2-p_1\|_1 +  \Ln(\Om)^{1/2} 
 \|Ew_2-Ew_1\|_2\,.
\end{equation}
By (\ref{n/n-1})--(\ref{Ev}) there exists a constant $A_2$, depending only on $\Om$ and $\Ga_0$, such that
\begin{equation}\label{A2}
 \|v\|_{n/(n-1)} \le A_2\,  \|e_2-e_1\|_2 + A_2\, \|p_2-p_1\|_1 + A_2 \,
 \|Ew_2-Ew_1\|_2\,.
  \end{equation} 
  
Since the trace operator  is continuous from $BD(\Om)$ into $L^1(\partial\Om;\Rn)$, there exists a constant $B_1$, depending only  on $\Om$, such that
\begin{equation}\label{1Ga1}
 \|v\|_{1,\Ga_1} \le B_1\,  \|v\|_1 + B_1\, \|Ev\|_1\,.
\end{equation}
{}From this inequality and from (\ref{trace})--(\ref{Ev}) we deduce that there exists a constant $B_2$, depending only on $\Om$ and $\Ga_0$, such that
\begin{equation}\label{B2}
\|v\|_{1,\Ga_1} \le B_2\,  \|e_2-e_1\|_2 + B_2\, \|p_2-p_1\|_1 + B_2 \,
 \|Ew_2-Ew_1\|_2\,.
\end{equation} 

Therefore (\ref{ineq0}), (\ref{A2}), and (\ref{B2}) imply that
$$
\begin{array}{c}
2\,\alpha_\C\, \|e_2-e_1\|_2^2 \le 2\,\beta_\C\,\|e_2-e_1\|_2\, \|Ew_2-Ew_1\|_2 
+A_2\, \|f_2-f_1\|_n\, \|e_2-e_1\|_2  +{} \\
{}+ A_2\, \|f_2-f_1\|_n\, \|p_2-p_1\|_1 +A_2\, \|f_2-f_1\|_n\, \|Ew_2-Ew_1\|_2
 +{}
\\
{}+B_2\, \|g_2-g_1\|_{\infty ,\Ga_1} \, \|e_2-e_1\|_2+B_2\, \|g_2-g_1\|_{\infty ,\Ga_1} \, \|p_2-p_1\|_1 +{}
\\
{}+B_2\, \|g_2-g_1\|_{\infty ,\Ga_1} \,\|Ew_2-Ew_1\|_2
+ 2\,R_K \|p_2-p_1\|_1\,,
\end{array}
 $$ 
which yields (\ref{eq12}) by the Cauchy inequality.

As $Eu_i=e_i+p_i$ in $\Om$ by (\ref{900}), by the H\"older inequality we obtain
$$
 \|Eu_2-Eu_1\|_1\le \Ln(\Om)^{1/2} 
\|e_2-e_1\|_2 +  \|p_2-p_1\|_1\,,
$$
so that (\ref{eq12}) gives (\ref{hoelder}).

Since $p_2-p_1=[(w_2-w_1)-(u_2-u_1)]{\,\odot\,}\nu\, \hn$ on $\Ga_0$, we have
$$
\|u_2-u_1\|_{1,\Ga_0} \le \|w_2-w_1\|_{1,\Ga_0} + \sqrt 2\, \|p_2-p_1\|_1 \,.
$$
The continuity of the trace operator from $H^1(\Om;\Rn)$ into $L^1(\partial\Om;\Rn)$ implies that there exists a constant $M$, depending only on $\Om$, such that 
$$
\|u_2-u_1\|_{1,\Ga_0} \le M\, \|w_2-w_1\|_2 + M\,\| E w_2- E w_1\|_2  + \sqrt 2\, \|p_2-p_1\|_1\,.
$$
By (\ref{seminorm}) there exists a constant $C$, depending only on $\Om$ and $\Ga_0$, such that
\begin{eqnarray*}
&\|u_2-u_1\|_1\le C\, \|u_2-u_1\|_{1,\Ga_0} + C\, \|Eu_2-Eu_1\|_1
\le
\\
&\le C\,M\, \|w_2-w_1\|_2 +  C\,M\, \| E w_2- E w_1\|_2+
\sqrt 2\,C\, \|p_2-p_1\|_1 + 
 C\, \|Eu_2-Eu_1\|_1 \,.
\end{eqnarray*}
Inequality (\ref{hoelderu}) follows now form~(\ref{hoelder}).
\end{proof}

\begin{remark}\label{unique} 
Theorem~\ref{thmhoelder} 
implies that, if $(u_1, e_1,p_0)$ and $(u_2, e_2,p_0)$ are solutions to problem 
(\ref{minp0}) with the same $w$, $f$, and $g$, then $u_1=u_2$ and $e_1=e_2$ a.e.\ on $\Om$. 
\end{remark}

\end{section}

\begin{section}{Quasistatic evolution}

We now consider time-dependent body and surface forces $f(t)$ and $g(t)$ satisfying the regularity assumptions and the uniform safe-load condition of Section~\ref{mech}. For every $t\in [0,T]$ the total load $\LL(t)\in BD(\Om)'$ applied at time $t$ is defined by
\begin{equation}\label{defLL}
\langle \LL(t)|u\rangle:= \langle f(t)|u\rangle+\langle g(t)|u\rangle_{\Ga_1}
\end{equation}
for every $u\in BD(\Om)$. 

\begin{remark}\label{LLdot}
{}From the hypotheses of  Section~\ref{mech} it follows that the weak$^*$ limit
$$
\dot \LL(t):= w^*\hbox{-}\lim_{s\to t} \frac{\LL(s)-\LL(t)}{s-t}
$$
exists in $BD(\Om)'$ for a.e.\ $t\in [0,T]$, and that
\begin{equation}\label{LLdot2}
\langle\dot \LL(t)|u\rangle:= \langle\dot f(t)|u\rangle+\langle\dot g(t)|u\rangle_{\Ga_1}
\end{equation}
for every $u\in BD(\Om)$.
Therefore the function $t\mapsto \langle\dot \LL(t)|u(t)\rangle$ belongs to $L^1([0,T])$ whenever $t\mapsto u(t)$ belongs to $L^\infty([0,T];BD(\Om))$.

The properties of $\dot\varrho(t)$ and $\dot\varrho_D(t)$ mentioned in Section~\ref{mech} imply that $\dot\varrho(t)\in\Sigma(\Om)$ for a.e.\ $t\in [0,T]$ and
$$
- \div\,\dot\varrho(t)=\dot f(t)\, \hbox{ a.e.\ on } \Om\,, \qquad [\dot\varrho(t)\nu]=\dot g(t) \hbox{ on } \Ga_1\,.
$$
Moreover, thanks to
(\ref{sigmakp}), we can prove that for every $p\in \Pi_{\Ga_0}(\Om)$ the function $s\mapsto \langle \varrho_D(s)|p\rangle$ is differentiable at each $t\in[0,T]$ where $\dot\varrho(t)$ exists and (\ref{dotsigmaD}) holds, with derivative given by
$\langle \dot \varrho_D(t)|p\rangle$. This implies that $t\mapsto \langle \dot \varrho_D(t)|p(t)\rangle$ is measurable for every simple function $t\mapsto p(t)$ from $[0,T]$ into $M_b(\Om\cup\Ga_0;\MD)$ with $p(t)\in  \Pi_{\Ga_0}(\Om)$ for a.e.\ $t\in [0,T]$. By approximation we conclude that $t\mapsto \langle \dot \varrho_D(t)|p(t)\rangle$ belongs to $L^1([0,T])$ whenever $t\mapsto p(t)$ belongs to $L^\infty([0,T];M_b(\Om\cup\Ga_0;\MD))$ and $p(t)\in  \Pi_{\Ga_0}(\Om)$ for a.e.\ $t\in [0,T]$.
\end{remark}

A function $p\colon [0,T]\to M_b(\Om\cup\Ga_0;\MD)$ will be regarded as a function defined on the time interval $[0,T]$ 
with values in the dual of the separable Banach space $C_0(\Om\cup\Ga_0;\MD)$. Therefore for every $s,t\in[0,T]$ with 
$s\le t$ the total variation of $p$ on $[s,t]$ is defined by
$$
\V(p;s,t)=\sup\Big\{\sum_{j=1}^N\|p(t_j)-p(t_{j-1})\|_1:\, s=t_0\le 
t_1\le \dots\le t_N=t, \, N\in\N
\Big\}\,.
$$
By (\ref{HHduality}) we can apply to $\HH$ all results proved in the Appendix with
$X= M_b(\Om\cup\Ga_0;\MD)$, $Y=C_0(\Om\cup\Ga_0;\MD)$, and $\K=\K_D(\Om)\cap C_0(\Om\cup\Ga_0;\MD)$. The $\HH$-variation of $p$ on $[s,t]$, which will 
play the role of the dissipation in the time interval $[s,t]$, is denoted $\D_\HH(p;s,t)$ and is defined by
\begin{equation}\label{diss}
\D_\HH(p;s,t):=\sup\Big\{ \sum_{j=1}^N \HH(p(t_j)-p(t_{j-1})):\, 
s=t_0\le t_1\le \dots\le t_N=t, \,
N\in\N \Big\}\, .
\end{equation}

\subsection{Definition of quasistatic evolution}
We are now in a position to introduce the following definition.

\begin{definition}\label{qstatic}
A {\em quasistatic evolution\/} is a function $t\mapsto(u(t),e(t),p(t))$ from 
$[0,T]$ into $BD(\Om){\times}L^2(\Om;\Mnn){\times}M_b(\Om\cup\Ga_0;\MD)$
which satisfies the following conditions:
\begin{itemize}
\smallskip
\item[{\rm (qs1)}] {\em global stability\/}: for every $t\in [0,T]$ 
we have $(u(t),e(t),p(t))\in 
A(w(t))$ and
\begin{equation}\label{min-ineq}
\QQ(e(t)) -\langle \LL(t)|u(t)\rangle 
\le \QQ(\eta)+\HH(q-p(t))
-\langle \LL(t)|v\rangle
\end{equation}
for every $(v,\eta,q)\in A(w(t))$;  
\smallskip 
\item[{\rm (qs2)}] {\em energy balance\/}: the function $t\mapsto p(t)$ from $[0,T]$ into $M_b(\Om\cup\Ga_0;\MD)$ has bounded variation and for every $t\in[0,T]$
\begin{equation}\label{en-eq}
\begin{array}{c}
\QQ(e(t)) + \D_\HH(p; 0,t) -\langle \LL(t)|u(t)\rangle = 
\QQ(e(0)) -\langle \LL(0)|u(0)\rangle +{}  
\\
\displaystyle
{}+ \int_0^t \{\langle \sigma(s)|E\dot w(s)\rangle -\langle \LL(s)|\dot w(s)\rangle
-\langle \dot \LL(s)|u(s)\rangle\}\, ds\,,
\end{array}
\end{equation}
where $\sigma(t):=\C e(t)$.
\smallskip
\end{itemize}
\end{definition}

\begin{remark}\label{contqo}
Since  the function $t\mapsto p(t)$ from $[0,T]$ into 
$M_b(\Om\cup\Ga_0;\MD)$ has bounded variation, it is bounded and the set of its discontinuity points (in the strong topology) is at most countable (see, e.g., \cite[Lemma~A.1]{Bre}).
By Theorem~\ref{thmhoelder} the same properties hold for the functions 
$t\mapsto e(t)$ and $t\mapsto \sigma(t)$ from $[0,T]$ into $L^2(\Om;\Mnn)$, and for the function
$t\mapsto u(t)$ from $[0,T]$ into $BD(\Om)$. Therefore $t\mapsto e(t)$ and $t\mapsto \sigma(t)$ belong to $L^\infty([0,T];L^2(\Om;\Mnn))$ and $t\mapsto u(t)$ belongs to $L^\infty([0,T];BD(\Om))$.
As $t\mapsto E\dot w(t)$ belongs to $L^1([0,T];L^2(\Om;\Mnn))$ and $t\mapsto \dot w(t)$ belongs to $L^1([0,T];H^1(\Rn;\Rn))$, the integral in the right-hand side of (\ref{en-eq}) is well defined thanks to Remark~\ref{LLdot}
\end{remark}

The following theorem gives an equivalent formulation of conditions (qs1) and (qs2), which uses the function $t\mapsto\varrho(t)$ introduced in the uniform safe-load condition of Section~\ref{mech}.

\begin{theorem}\label{qstatic2}
A function $t\mapsto(u(t),e(t),p(t))$ from 
$[0,T]$ into $BD(\Om){\times} L^2(\Om;\Mnn)\allowbreak{\times} M_b(\Om\cup\Ga_0;\MD)$ is a quasistatic evolution if and only if it
satisfies the following conditions:
\begin{itemize}
\smallskip
\item[$(\rm qs1')$]for every $t\in [0,T]$ we have
$(u(t),e(t),p(t))\in A(w(t))$ and
\begin{equation}\label{min-ineq-v}
\qquad\qquad\QQ(e(t)) -\langle \varrho(t)|e(t)\rangle 
\le \QQ(\eta) -\langle \varrho(t)|\eta\rangle +\HH(q-p(t))
-\langle \varrho_D(t)|q-p(t)\rangle
\end{equation}
for every $(v,\eta,q)\in A(w(t))$;
\smallskip
\item[$(\rm qs2')$]the function $t\mapsto p(t)$ from $[0,T]$ into $M_b(\Om\cup\Ga_0;\MD)$ has bounded variation and for every $t\in[0,T]$
\begin{equation}\label{en-eq-v}
\begin{array}{c}
\QQ(e(t)) + \D_\HH(p; 0,t) -\langle \varrho(t)|e(t)-Ew(t)\rangle -\langle \varrho_D(t)|p(t)\rangle = \vspace{.25cm}
\\
\displaystyle =\QQ(e(0)) -\langle \varrho(0)|e(0)-Ew(0)\rangle -\langle \varrho_D(0)|p(0)\rangle
 - {}
\\
\displaystyle - \int_0^t \{\langle \dot\varrho(s)|e(s)-Ew(s)\rangle +\langle \dot \varrho_D(s)|p(s)\rangle\}\,ds + \int_0^t \langle \sigma(s)|E\dot w(s)\rangle\, ds
\,,
\end{array}
\end{equation}
where $\sigma(t):=\C e(t)$.
\smallskip
\end{itemize}
\end{theorem}

\begin{proof} The equivalence of conditions (qs1) and $(\rm qs1')$ follows from Lemma~\ref{fg}.

As the functions $t\mapsto f(t)$, $t\mapsto g(t)$, and $t\mapsto w(t)$ are absolutely continuous from $[0,T]$ into $L^n(\Om;\Rn)$, $L^\infty(\Ga_1;\Rn)$, and $H^1(\Rn;\Rn)$, respectively, the function $t\mapsto \langle\LL(t)|w(t)\rangle$ is absolutely continuous on $[0,T]$ and its time derivative is given by
$t\mapsto \langle\dot \LL(t)|w(t)\rangle +  \langle\LL(t)|\dot w(t)\rangle$. It follows that
\begin{equation}\label{intproduct}
\int_0^t\{ \langle\dot \LL(s)|w(s)\rangle +  \langle\LL(s)|\dot w(s)\rangle\}\,ds = \langle\LL(t)|w(t)\rangle - \langle\LL(0)|w(0)\rangle\,.
\end{equation}

By Lemma~\ref{fg} we have
\begin{equation}\label{formula0}
\langle\LL(t)|v\rangle= \langle \varrho(t)|\eta-Ez\rangle +\langle \varrho_D(t)|q\rangle +\langle\LL(t)|z\rangle
\end{equation}
for every $t\in[0,T]$,  $z\in H^1(\Rn;\Rn)$, and $(v,\eta,q)\in A(z)$. Taking the derivative with respect to $t$, thanks to Remark~\ref{LLdot} we obtain
$$
\langle\dot \LL(t)|v\rangle= \langle\dot  \varrho(t)|\eta-Ez\rangle +\langle\dot  \varrho_D(t)|q\rangle +\langle\dot \LL(t)|z\rangle
$$
for a.e.\ $t\in[0,T]$, for every  $z\in H^1(\Rn;\Rn)$, and every $(v,\eta,q)\in A(z)$.

If conditions (qs1) or $(\rm qs1')$ hold, then by Remark~\ref{contqo} the function
$t\mapsto(u(t),e(t),p(t))$ belongs to
$L^\infty([0,T];BD(\Om){\times}L^2(\Om;\Mnn){\times}M_b(\Om\cup\Ga_0;\MD))$. As
$(u(t),e(t),p(t))\in A(w(t))$ for every $t\in[0,T]$, we have
$$
\langle\dot \LL(t)|u(t)\rangle= \langle\dot  \varrho(t)|e(t)-Ew(t)\rangle +\langle\dot  \varrho_D(t)|p(t)\rangle +\langle\dot \LL(t)|w(t)\rangle
$$
for a.e.\ $t\in[0,T]$. Therefore (\ref{intproduct}) implies that
\begin{equation}\label{formula1}
\begin{array}{c}
\displaystyle
\int_0^t\{ \langle\dot \LL(s)|u(s)\rangle +  \langle\LL(s)|\dot w(s)\rangle\}\,ds = \langle\LL(t)|w(t)\rangle - \langle\LL(0)|w(0)\rangle+{}
\\
\displaystyle
{}+\int_0^t\{ \langle \dot\varrho(s)|e(s)-Ew(s)\rangle +\langle \dot \varrho_D(s)|p(s)\rangle\}\,ds
\,.
\end{array}
\end{equation}
The equivalence of  conditions (qs2) and $(\rm qs2')$ follows now from~(\ref{formula0}) 
and~(\ref{formula1}).
\end{proof}

\subsection{The existence result}
The following theorem is the main result of the paper.

\begin{theorem}\label{main}
Let $(u_0,e_0,p_0)\in A(w(0))$ satisfy the stability condition
\begin{equation}\label{min-ineq0}
\QQ(e_0) -\langle \LL(0)|u_0 \rangle
\le \QQ(e)+\HH(p-p_0)-\langle \LL(0)|v\rangle 
\end{equation}
for every $(u,e,p)\in A(w(0))$. Then there exists a quasistatic evolution $t\mapsto(u(t),e(t),p(t))$
such that $u(0)=u_0$, $e(0)=e_0$, $p(0)=p_0$.
\end{theorem}

Theorem~\ref{main} will be proved by a time discretization process.
Let us fix a sequence of subdivisions $(\tki)_{0\le i\le k}$ of the 
interval $[0,T]$, with
\begin{eqnarray}
& 0=t_k^0<t_k^1<\dots<t_k^{k-1}<t_k^{k}=T\,,
\label{subdiv5}\\
&\displaystyle
\lim_{k\to\infty}\,
\max_{1\le i\le k} (\tki-\tkim)= 0\,.
\label{fine5}
\end{eqnarray}
{}For $i=0,\ldots,k$ we set $\wki:=w(\tki)$, $f_k^i:=f(\tki)$, $g_k^i:=g(\tki)$, $\LL_k^i:=\LL(\tki)$, and $\varrho_k^i:=\varrho(\tki)$.

{}For every $k$ we define $\uki$, $e_k^i$, and $\pki$ by induction.
We set $(u_k^0,e_k^0, p_k^0):=(u_0,e_0,p_0)$, which, by assumption, belongs to $A(w(0))$, and for $i=1,\ldots,k$ we 
define $(\uki,e_k^i,\pki)$ as a solution to the incremental problem
\begin{equation}\label{min-inc}
\min_{(u,e,p)\in A(\wki)} \{\QQ(e)+\HH(p-p_k^{i-1}) -\langle \LL_k^i|u\rangle \}\,.
\end{equation}
The existence of a solution to this problem is proved in 
Theorem~\ref{existencemin}. We recall that by Lemma~\ref{fg} the minimum problem (\ref{min-inc}) is equivalent to
\begin{equation}\label{min-eqvl}
\min_{(u,e,p)\in A(\wki)} \{\QQ(e)-\langle\varrho_k^i|e\rangle
+\HH(p-p_k^{i-1}) -\langle (\varrho_k^i)_D |p-p_k^{i-1} \rangle 
\}\,.
\end{equation}
Moreover, by the triangle inequality (\ref{triangle}) the triple $(u_k^i,e_k^i,p_k^i)$ is also a solution of the problem
\begin{equation}\label{min-inc2}
\min_{(u,e,p)\in A(\wki)} \{\QQ(e)+\HH(p-p_k^i) -\langle \LL_k^i|u\rangle \}\,.
\end{equation}

{}For $i=0,\ldots,k$ we set $\sigma_k^i:=\C e_k^i$ and for every $t\in[0,T]$ we define the piecewise constant interpolations
\begin{equation}\label{ukt}
\begin{array}{c} 
\uk(t):=\uki\,, \quad e_k(t):=e_k^i\,,  \quad   \pk(t):=\pki\,,  \quad
\sigma_k(t):=\sigma_k^i\,, \vspace{.1cm}\\
w_k(t):=w_k^i \,,
\quad f_k(t):=f_k^i \,, \quad g_k(t):=g_k^i \,,\quad \LL_k(t):=\LL_k^i \,, \quad \varrho_k(t):=\varrho_k^i \,,
\end{array}
\end{equation}
where $i$ is the largest integer such that $\tki\le t$. 
By definition $(u_k(t),e_k(t),p_k(t))\in A(w_k(t))$ and by (\ref{min-inc2}) we have
\begin{equation}\label{min-inc3}
\QQ(e_k(t))-\langle \LL_k(t)|u_k(t)\rangle
 \le \QQ(\eta)+\HH(q-p_k(t))-\langle \LL_k(t)|v\rangle
\end{equation}
for every $(v,\eta,q)\in A(w_k(t))$.

\subsection{The discrete energy inequality}

We now derive an energy estimate for the solutions of the incremental problems. 
Note that a remainder $\delta_k$ is needed because the integral terms which appear in the right-hand side of (\ref{inc-ineq}) 
provide only an approximate value of the work done by the external forces.

\begin{lemma}\label{lm:41}
There exists a sequence $\delta_k\to 0^+$ such that for every $k$ and every $t\in[0,T]$
\begin{equation}\label{inc-ineq}
\begin{array}{l}
\displaystyle
\hskip-1em
\QQ(e_k(t))-\langle\varrho_k(t)|e_k(t)-Ew_k(t)\rangle +{}
\vspace{.2cm}
\\
\displaystyle
\quad{}+\sum_{0<t_k^r\le t}\{\HH(p_k^r-p_k^{r-1})
-\langle \varrho_D(t_k^r)|p_k^r-p_k^{r-1} \rangle\} \le 
\\
\displaystyle
\le \QQ(e_0)-\langle\varrho(0)|e_0-Ew(0)\rangle -{}
\\
\displaystyle
\quad{}-\int_0^{\tki} \langle\dot\varrho(s)|e_k(s)-Ew_k(s)\rangle\,ds 
+\int_0^{\tki}\langle\sigma_k(s)|E\dot w(s)\rangle\, ds +\delta_k\,,
\end{array}
\end{equation}
where $i$ is the largest integer such that $t_k^i\le t$.
\end{lemma}

The integrals in the right-hand side of (\ref{inc-ineq}) can be written as
$$\begin{array}{c}
\displaystyle \int_0^{t_k^i}\langle 
\dot\varrho(s)|e_k(s)-Ew_k(s)\rangle\, ds= \sum_{j=1}^i 
\langle 
\varrho^j_k-\varrho^{j-1}_k|e^{j-1}_k-Ew^{j-1}_k\rangle\,, \\
\displaystyle \int_0^{t_k^i}\langle 
\sigma_k(s)|E\dot w(s)\rangle\, ds= \sum_{j=1}^i 
\langle 
\sigma^{j-1}_k|Ew^j_k- Ew^{j-1}_k\rangle\,,
\end{array}
$$
where the sums involve only the values of $\varrho(t)$ and $w(t)$ at the discretization points~$t^j_k$. 
This is the main difference between inequality  
(\ref{inc-ineq}) and those considered in \cite[Theorem~4.1]{Mai-Mie}.

\begin{proof}[Proof of Lemma~\ref{lm:41}]
We have to prove that there exists a sequence $\delta_k\to 0^+$ such that
\begin{equation}\label{inc-ineq-2}
\begin{array}{l}
\displaystyle
\hskip-1em
\QQ(e_k^i)-\langle\varrho_k^i|e_k^i-Ew_k^i\rangle +{}
\\
\displaystyle
{}\quad +
\sum_{r=1}^i\{\HH(p_k^r-p_k^{r-1})
-\langle (\varrho_k^r)_D|p_k^r-p_k^{r-1} \rangle\} \le
\\ 
\displaystyle
\le \QQ(e_0)-\langle\varrho(0)|e_0-Ew(0)\rangle -{}
\\
\displaystyle
{}\quad - \int_0^{\tki} \langle\dot\varrho(s)|e_k(s)-Ew_k(s)\rangle\,ds 
+\int_0^{\tki}\langle \sigma_k(s)|E\dot 
w(s)\rangle\, ds +\delta_k
\end{array}
\end{equation}
for every $k$ and every $i=1,\ldots,k$.

Let us fix an integer $r$ with $1\le r\le i$ and let $v:=u_k^{r-1}-w_k^{r-1}+w_k^r$ and
$\eta:=e_k^{r-1}-Ew_k^{r-1}+Ew_k^r$.
Since $(v,\eta,p_k^{r-1})\in A(w_k^r)$, by the 
minimality condition (\ref{min-eqvl}) we have
\begin{equation}\label{b01}
\begin{array}{c}
\QQ(e_k^r)-\langle\varrho_k^r|e_k^r\rangle 
+\HH(p_k^r-p_k^{r-1}) -\langle (\varrho_k^r)_D|p_k^r-p_k^{r-1} \rangle \le \vspace{.1cm}\\
\le 
\QQ(e_k^{r-1}+Ew_k^r -Ew_k^{r-1})-\langle\varrho_k^r|e_k^{r-1}+Ew_k^r -Ew_k^{r-1}\rangle
\,,
\end{array}
\end{equation}
where the quadratic form in the right-hand side can be developed as
\begin{equation}\label{b02}
\QQ(e_k^{r-1}+Ew_k^r-Ew_k^{r-1})
=\QQ(e_k^{r-1})+\langle \sigma_k^{r-1} | Ew_k^r- 
Ew_k^{r-1} \rangle 
+ \QQ(Ew_k^r-Ew_k^{r-1})\,.
\end{equation}
{}From the absolute continuity of $w$ with respect to $t$ we obtain
$$
w_k^r-w_k^{r-1}=\int_{t_k^{r-1}}^{t_k^r}\dot w(t)\, dt\,,
$$
where we use a Bochner integral of a function with values in $H^1(\Rn;\Rn)$.
This implies that
\begin{equation}\label{b03}
Ew_k^r- Ew_k^{r-1}=\int_{t_k^{r-1}}^{t_k^r} E\dot w(t)\, dt\, ,
\end{equation}
where we use a Bochner integral of a function with values in $L^2(\Rn;\Mnn)$.
By (\ref{boundsC}) and (\ref{b03}) we get
\begin{equation}\label{b03.5}
\QQ(Ew_k^r-Ew_k^{r-1})\le \beta_{\C}\Big(\int_{t_k^{r-1}}^{t_k^r} 
\|E\dot w(t)\|_2\, dt\Big)^2. 
\end{equation}
{}From the absolute continuity of $\varrho$ with respect to $t$ we have
\begin{equation}\label{zb03}
\langle\varrho_k^r|e_k^{r-1}-Ew_k^{r-1}\rangle =
\langle\varrho_k^{r-1}|e_k^{r-1}-Ew_k^{r-1}\rangle +\int_{t_k^{r-1}}^{t_k^r} 
\langle \dot\varrho(t)|e_k^{r-1}-Ew_k^{r-1}\rangle\,dt\,. 
\end{equation}
By (\ref{b01})--(\ref{zb03}) we obtain
\begin{equation}\label{b04}
\begin{array}{l}
\hskip-1em
\displaystyle
\vphantom{\int_{t_k^{r-1}}^{t_k^r}}
 \QQ(e_k^r)-\langle\varrho_k^r|e_k^r-Ew_k^r\rangle 
+\HH(p_k^r-p_k^{r-1}) -\langle (\varrho_k^r)_D|p_k^r-p_k^{r-1} \rangle \le  
\\
\displaystyle \le \QQ(e_k^{r-1})-\langle\varrho_k^{r-1}|e_k^{r-1}-Ew_k^{r-1}\rangle -
\int_{t_k^{r-1}}^{t_k^r} 
\langle \dot\varrho(t)|e_k^{r-1}-Ew_k^{r-1}\rangle\,dt +{}
 \\
\displaystyle {}\quad+\int_{t_k^{r-1}}^{t_k^r} \langle \sigma_k^{r-1} |E\dot w(t)\rangle \, dt +
\beta_{\C}\Big(\int_{t_k^{r-1}}^{t_k^r}\|E\dot w(t)\|_2\, dt 
\Big)^2 \le
\\
\displaystyle \le \QQ(e_k^{r-1})-\langle\varrho_k^{r-1}|e_k^{r-1}-Ew_k^{r-1}\rangle -
\int_{t_k^{r-1}}^{t_k^r} 
\langle \dot\varrho(t)|e_k^{r-1}-Ew_k^{r-1}\rangle\,dt +{} 
\\
\displaystyle {}\quad +\int_{t_k^{r-1}}^{t_k^r} \langle \sigma_k^{r-1} |E\dot w(t)\rangle \, dt +
\omega_k \int_{t_k^{r-1}}^{t_k^r}\|E\dot w(t)\|_2\, dt \, ,
\end{array}
\end{equation}
where
$$\omega_k:= \beta_{\C}\max_{1\le r\le k}\int_{t_k^{r-1}}^{t_k^r}\|E\dot 
w(t)\|_2\, dt \, \to 0$$
by the absolute continuity of the integral.
Iterating now inequality (\ref{b04}) for $1\le r\le i$, we get 
(\ref{inc-ineq-2}) 
with $\delta_k:=\omega_k\int_0^T \|E\dot w(t)\|_2\, dt$.
\end{proof}

\subsection{Proof of the existence theorem}
We are now in a position to prove Theorem~\ref{main}.

\begin{proof}[Proof of Theorem~\ref{main}]
Let us fix a sequence of subdivisions $(\tki)_{0\le i\le k}$ of the 
interval $[0,T]$ satisfying
(\ref{subdiv5}) and (\ref{fine5}). For every $k$ let $(\uki,e_k^i,\pki)$,
$i=1,\ldots,k$, be defined inductively
 as solutions of the discrete 
problems (\ref{min-inc}), with
$(u_k^0,e_k^0, p_k^0)=(u_0,e_0,p_0)$, and let 
$u_k(t)$, $e_k(t)$, $p_k(t)$, $\sigma_k(t)$, $w_k(t)$, $f_k(t)$, $g_k(t)$, $\LL_k(t)$, 
$\varrho_k(t)$ be defined by (\ref{ukt}).

Let us prove that there exists a constant $C$, depending only on the constants $\alpha_\C$, $\beta_\C$, and $\alpha$, and on the functions $e_0$, $t\mapsto w(t)$, 
and $t\mapsto \varrho(t)$, such that
\begin{equation}\label{b05}
\sup_{t\in[0,T]}\|e_k(t)\|_2\leq C \qquad\hbox{and} 
\qquad \V(\pk;0,T)\leq C
\end{equation}
for every $k$.
As $t\mapsto w(t)$ and $t\mapsto \varrho(t)$ are absolutely continuous with values in 
$H^1(\Rn;\Rn)$ and $L^2(\Om;\Mnn)$, respectively, the functions
$t\mapsto\|E w(t)\|_2$ and $t\mapsto\| \varrho(t)\|_2$ are bounded on $[0,T]$ and
the functions $t\mapsto\|E\dot w(t)\|_2$ and $t\mapsto\|\dot \varrho(t)\|_2$
are integrable on $[0,T]$. This fact, together with (\ref{boundsC}), (\ref{normC}), (\ref{eps*1}), and (\ref{inc-ineq}), implies that
\begin{equation}\label{young}
\begin{array}{c}
\displaystyle
\alpha_{\C}\|e_k(t)\|^2_2 -\sup_{t\in[0,T]}\|\varrho(t)\|_2 
\sup_{t\in[0,T]}\|Ew(t)\|_2+\alpha\sum_{0<t_k^r\le t}\|p_k^r-p_k^{r-1}\|_1\le
\\
\displaystyle
\le 
\beta_{\C}\|e_0\|^2_2+
\|\varrho(0)\|_2
\big(\|e_0\|_2
+\|Ew(0)\|_2\big)
+\sup_{t\in[0,T]}\|Ew(t)\|_2 
\int_0^T  \|\dot \varrho(s)\|_2 \, ds
+{}
\\
\displaystyle
+\sup_{t\in[0,T]}\|e_k(t)\|_2 \Big( 
\int_0^T  \|\dot \varrho(s)\|_2 \, ds
+2\beta_\C\int_0^T  \|E\dot w(s)\|_2 \, ds 
+\sup_{t\in[0,T]}\|\varrho(t)\|_2
\Big)+\delta_k
\end{array}
\end{equation}
for every $k$ and every $t\in[0,T]$. The 
former inequality in (\ref{b05}) can be obtained now by using the Cauchy inequality.
As for the latter, by (\ref{young}) and the first inequality in 
(\ref{b05}) we deduce that 
\begin{equation}\label{c02}
\sum_{0<t_k^r\le t}\|p_k^r-p_k^{r-1}\|_1\le C
\end{equation}
for every $k$ and every $t\in[0,T]$. Since $t\mapsto\pk(t)$ is constant on on the intervals 
${[t^{r-1}_k,t^r_k[}$, the estimate (\ref{c02}) is equivalent to the second inequality in 
(\ref{b05}).

By the generalized version of the classical Helly theorem given in Lemma~\ref{lm:44} there exist a subsequence, still denoted $p_k$, and a function $p\colon [0,T]\to 
M_b(\Om\cup\Ga_0;\MD)$, with bounded 
variation on $[0,T]$, such that $p_k(t)\wto p(t)$ weakly$^*$ in 
$M_b(\Om\cup\Ga_0;\MD)$ for every $t\in [0,T]$.

Since, by (\ref{b05}), $\|e_k(t)\|_2\leq C$ and $\|p_k(t)\|_1\leq C$ for every $k$ and every $t$, arguing as in the proof of The\-orem~\ref{existencemin} we deduce that
$u_k(t)$ is bounded in  $BD(\Om)$ uniformly with respect to $k$ and~$t$. Let us fix $t\in[0,T]$.
There exist an increasing sequence $k_j$  (possibly depending on $t$) and two functions $u(t)\in BD(\Om)$ and $e(t)\in L^2(\Om;\Mnn)$
such that $u_{k_j}(t)\wto u(t)$ weakly$^*$ in $BD(\Om)$ and $e_{k_j}(t)\wto e(t)$ weakly in $L^2(\Om;\Mnn)$.
By (\ref{min-inc3}) we can apply Theorem~\ref{lm:45} and we obtain that the triple $(u(t),e(t),p(t))$ is a solution of the minimum 
problem
\begin{equation}\label{f01}
\min_{(v,\eta,q)\in A(w(t))} \{\QQ(\eta)+\HH(q-p(t))-\langle \LL(t)|v\rangle \}\,.
\end{equation}
By Remark~\ref{unique} there exists a unique $(u,e)\in BD(\Om){\times}L^2(\Om;\Mnn)$
such that $(u,e,p(t))$ is a solution to (\ref{f01}).
Therefore, the convergence result holds for the whole sequence, i.e.,  $u_k(t)\wto u(t)$ weakly$^*$ in $BD(\Om)$ and $e_k(t)\wto e(t)$ weakly in $L^2(\Om;\Mnn)$.

Let us show now that the function $t\mapsto (u(t),e(t),p(t))$ is a quasistatic 
evolution satisfying $(u(0),e(0),p(0))=(u_0,e_0,p_0)$.
The initial condition is fulfilled, since $u_k(0)=u_0$, $e_k(0)=e_0$,
$p_k(0)=p_0$ for every~$k$.
In (\ref{f01}) we have already proved that $(u(t),e(t),p(t))$ satisfies 
(\ref{min-ineq})
for every $t\in[0,T]$.

It remains to prove the energy balance (\ref{en-eq}), or equivalently
(\ref{en-eq-v}). By Theorem~\ref{inequality}, proved below, it is enough to establish the energy inequality
\begin{equation}\label{en-ineq}
\begin{array}{c}
\displaystyle
\QQ(e(t)) -\langle\varrho(t)|e(t)-Ew(t)\rangle+ \D_\HH(p; 0,t)-\langle\varrho_D(t)|p(t)\rangle \le \vspace{.2cm}
\\
\displaystyle\le\QQ(e(0)) -\langle\varrho(0)|e(0)-Ew(0)\rangle - \langle\varrho_D(0)|p(0)\rangle - {}
\\
\displaystyle -\int_0^t\{\langle\dot\varrho(s)|e(s)-Ew(s)\rangle +\langle\dot\varrho_D(s)|p(s)\rangle \}\,ds 
+ \int_0^t \langle \sigma(s)|E\dot w(s)\rangle\, ds\,.
\end{array}
\end{equation}

Let us fix $t\in [0,T]$. As in the proof of Theorem~\ref{existencemin}, let $\delta>0$ and 
$\psi_\delta(x):=\phi(\frac{1}{\delta}\dist(x,\Gamma_1))$
for every $x\in\ol\Om$, where $\phi\in C^\infty(\R)$, $0\le \phi\le 1$, $\phi(s)=0$ for $s\le 1$, and
$\phi(s)=1$ for $s\ge 2$. Since the measure $H(p_k^r-p_k^{r-1})-[\varrho_D(t_k^r)
{\,:\,}(p_k^r-p_k^{r-1})]$ is nonnegative on 
$\Om\cup\Ga_0$ by (\ref{Hpge}),
we have
\begin{equation}\label{dineq}
\HH(\psi_\delta(p_k^r-p_k^{r-1}))
- \langle [\varrho_D(t_k^r){\,:\,}(p_k^r-p_k^{r-1})]| \psi_\delta \rangle \le
\HH(p_k^r-p_k^{r-1})
- \langle \varrho_D(t_k^r)|p_k^r-p_k^{r-1} \rangle
\end{equation}
for every $r=1,\dots,i$.
Since $t\mapsto p_k(t)$ is constant on the intervals ${[t^{r-1}_k,t^r_k[}$, we have
$$
\D_\HH(\psi_\delta p_k;0,t) \le \sum_{0<t_k^r\le t} \HH(\psi_\delta(p_k^r-p_k^{r-1}))\,,
$$
so that the lower semicontinuity of the dissipation (see (\ref{semidiss})) gives
\begin{equation}\label{lscdiss}
\D_\HH(\psi_\delta p;0,t)\le \liminf_{k\to\infty} \sum_{0<t_k^r\le t} 
\HH(\psi_\delta (p_k^r-p_k^{r-1}))
\,.
\end{equation}
It is convenient to write
\begin{equation}\label{devel}
\begin{array}{c}
\displaystyle\sum_{r=1}^i\langle [\varrho_D(t_k^r){\,:\,}(p_k^r-p_k^{r-1})]|\psi_\delta \rangle =
-\sum_{r=1}^i\langle [(\varrho_D(t_k^r)-\varrho_D(t_k^{r-1})){\,:\,}p_k^{r-1}]|\psi_\delta \rangle +{}
\vspace{.1cm} \\
{}+\langle [\varrho_D(t_k^i){\,:\,}p_k^i]|\psi_\delta \rangle
-\langle [\varrho_D(0){\,:\,}p_0]|\psi_\delta \rangle
\,.
\end{array}
\end{equation}
Since $t\mapsto\varrho(t)$ and $t\mapsto f(t)$ are absolutely continuous from $[0,T]$ into 
$L^2(\Om;\Mnn)$ and $L^n(\Om;\Rn)$, respectively, by (\ref{intpartsp}) we have that
$$
\begin{array}{c}
\displaystyle\sum_{r=1}^i\langle [(\varrho_D(t_k^r)-\varrho_D(t_k^{r-1})){\,:\,}p_k^{r-1}]|\psi_\delta \rangle =
-\int_0^{\tki} \langle \dot\varrho(s)|\psi_\delta ( e_k(s)-E w_k(s))\rangle\,ds - {} 
\\
\displaystyle
{}-\int_0^{\tki} \langle \dot\varrho(s)|(u_k(s)-w_k(s)){\,\odot\,}\nabla\psi_\delta \rangle\,ds  +\int_0^{\tki} \langle \dot f(s)|\psi_\delta (u_k(s)-w_k(s))\rangle\,ds
\,.
\end{array}
$$
Passing to the limit as $k\to\infty$ and using  (\ref{intpartsp}) again, we obtain
\begin{equation}\label{devel2}
\lim_{k\to\infty}\sum_{r=1}^i\langle [(\varrho_D(t_k^r)-\varrho_D(t_k^{r-1})){\,:\,}p_k^{r-1}]|\psi_\delta \rangle =
\int_0^t \langle [\dot\varrho_D(s){\,:\,}p(s)]|\psi_\delta \rangle\,ds\,.
\end{equation}
Analogously we can show that
\begin{equation}\label{devel3}
\lim_{k\to\infty} \langle [\varrho_D(t_k^i){\,:\,}p_k^i]|\psi_\delta \rangle =
\langle [\varrho_D(t){\,:\,}p(t)]|\psi_\delta \rangle\,.
\end{equation}
Combining together (\ref{dineq})--(\ref{devel3}),
we obtain that
$$
\begin{array}{c}
\displaystyle \D_\HH(\psi_\delta p;0,t) - \langle [\varrho_D(t){\,:\,}p(t)]|\psi_\delta \rangle
+\langle [\varrho_D(0){\,:\,}p(0)]|\psi_\delta \rangle 
+\int_0^t \langle [\dot\varrho_D(s){\,:\,}p(s)]|\psi_\delta \rangle\,ds \le
\\
\displaystyle\le \liminf_{k\to\infty}\sum_{r=1}^i\{\HH(p_k^r-p_k^{r-1})
- \langle \varrho_D(t_k^r)|p_k^r-p_k^{r-1} \rangle\} 
\end{array}
$$
and passing to the limit as $\delta\to 0^+$, we conclude that
\begin{equation}\label{wdineq}
\begin{array}{c}
\displaystyle \D_\HH(p;0,t) - \langle \varrho_D(t)|p(t) \rangle
+\langle \varrho_D(0)|p(0) \rangle 
+\int_0^t \langle \dot\varrho_D(s)|p(s) \rangle\,ds \le
\\
\displaystyle\le \liminf_{k\to\infty}\sum_{r=1}^i\{\HH(p_k^r-p_k^{r-1})
- \langle \varrho_D(t_k^r)|p_k^r-p_k^{r-1} \rangle\}\,. 
\end{array}
\end{equation}
{}For every $s\in [0,t]$ we have $\sigma_k(s)=\C e_k(s) \wto\C e(s)=\sigma(s)$ weakly in $L^2(\Om;\Mnn)$. 
As $\sigma_k(s)$ is bounded in $L^2(\Om;\Mnn)$ uniformly with respect to $k$ and $s$,
we can pass to the limit  in (\ref{inc-ineq}) as $k\to\infty$ and
we obtain (\ref{en-ineq}) from (\ref{wdineq}) and from the lower semicontinuity of~$\QQ$.
\end{proof}

As in  \cite[Theorem~4.4]{Mai-Mie} and  \cite[Lemma~7.1]{DM-Fra-Toa}, the energy inequality (\ref{en-ineq}) together 
with the global stability $\rm (qs1')$ imply the exact energy balance $\rm (qs2')$.

\begin{theorem}\label{inequality}
Let $t\mapsto(u(t),e(t),p(t))$ be a function from 
$[0,T]$ into $BD(\Om){\times}L^2(\Om;\Mnn)\allowbreak{\times}M_b(\Om\cup\Ga_0;\MD)$
which satisfies the stability condition $\rm (qs1')$ in Theorem~\ref{qstatic2}. Assume that 
$t\mapsto p(t)$ from $[0,T]$ into
$M_b(\Om\cup\Ga_0;\MD)$ has bounded variation. Then for every $t\in[0,T]$we have
\begin{equation}\label{en-ineq2}
\begin{array}{c}
\displaystyle
\QQ(e(t)) -\langle\varrho(t)|e(t)-Ew(t)\rangle+ \D_\HH(p; 0,t)-\langle\varrho_D(t)|p(t)\rangle \ge \vspace{.2cm}
\\
\displaystyle\ge\QQ(e(0)) -\langle\varrho(0)|e(0)-Ew(0)\rangle - \langle\varrho_D(0)|p(0)\rangle - {}
\\
\displaystyle -\int_0^t\{\langle\dot\varrho(s)|e(s)-Ew(s)\rangle +\langle\dot\varrho_D(s)|p(s)\rangle \}\,ds 
+ \int_0^t \langle \sigma(s)|E\dot w(s)\rangle\, ds\,,
\end{array}
\end{equation}
where $\sigma(t):=\C e(t)$. If, in addition, (\ref{en-ineq}) is satisfied, then the exact energy balance $\rm (qs2')$ holds.
\end{theorem}

\begin{proof}
Let us fix $t\in(0,T]$ and let $(s_k^i)_{0\le i\le k}$ be a sequence of 
subdivisions of the interval $[0,t]$
satisfying
\begin{eqnarray}
& 0=s_k^0<s_k^1<\dots<s_k^{k-1}<s_k^k=t\,,
\label{subdiv}\\
&\displaystyle
\lim_{k\to\infty}\,
\max_{1\le i\le k} (s_k^i-s_k^{i-1})= 0\,.
\label{fine}
\end{eqnarray}
{}For every $i=1,\ldots,k$ let $v:=u(s_k^i)-w(s_k^i)+ 
w(s_k^{i-1})$ and $\eta:=e(s_k^i)-Ew(s_k^i)+ 
Ew(s_k^{i-1})$. Since $(v,\eta,p(s_k^i))\in A(w(s_k^{i-1}))$, by the
global stability (\ref{min-ineq-v}) we have
\begin{equation}\label{n01}
\begin{array}{c}
\QQ(e(s_k^{i-1}))-\langle\varrho(s_k^{i-1})|e(s_k^{i-1})\rangle\le \\
\le \QQ(e(s_k^i)-(Ew(s_k^i)-Ew(s_k^{i-1}))) 
-\langle\varrho(s_k^{i-1})|e(s_k^i)-(Ew(s_k^i)-Ew(s_k^{i-1}))\rangle + {}
\\
+\HH(p(s_k^i)-p(s_k^{i-1}))-\langle\varrho_D(s_k^{i-1})|p(s_k^i)-p(s_k^{i-1})\rangle \,.
\end{array}
\end{equation}
The first term in the right-hand side can be written as
$$
\begin{array}{c}
\QQ(e(s_k^i)-(Ew(s_k^i)-Ew(s_k^{i-1})))= \\
= \QQ(e(s_k^i))-\langle 
\sigma(s_k^i)|Ew(s_k^i)-Ew(s_k^{i-1})\rangle 
+\QQ(Ew(s_k^i)-Ew(s_k^{i-1}))\,.
\end{array}
$$
Now, arguing as in (\ref{zb03}) and in the proof of the last inequality in (\ref{b04}), from the
previous equality and from (\ref{n01})
we obtain that there exists a sequence $\omega_k\to 0^+$ such that
$$
\begin{array}{c}
\QQ(e(s_k^{i-1}))-\langle \varrho(s_k^{i-1})| e(s_k^{i-1})-Ew(s_k^{i-1})\rangle 
-\langle\varrho_D(s_k^{i-1})|p(s_k^{i-1}) \rangle\le \vspace{.1cm} \\
\le \QQ(e(s_k^i)) +  \HH(p(s_k^i)-p(s_k^{i-1})) -\langle\varrho(s_k^i)|e(s_k^i)-Ew(s_k^i)\rangle
 -\langle\varrho_D(s_k^i)|p(s_k^i) \rangle + {}
\\ 
\displaystyle
+\int_0^t\langle\dot\varrho(s)|e(s_k^i)-Ew(s_k^i)\rangle\,ds
+\int_0^t\langle\dot\varrho_D(s)|p(s_k^i) \rangle\,ds - {}
\\
\displaystyle
-\int_{s_k^{i-1}}^{s_k^i}\langle 
\sigma(s_k^i)|E\dot w(s)
\rangle\,ds + \omega_k \int_{s_k^{i-1}}^{s_k^i}\| E\dot w(s)\|_2\,ds\,.
\end{array}
$$

On $[0,t]$ we define the piecewise constant functions 
$$
\begin{array}{c}
\ol e_k(s):=e(s_k^i)\,, \quad  E\ol w_k(s):=Ew(s_k^i)\,, \quad 
\ol p_k(s):=p(s_k^i)\,, \quad  \ol \sigma_k(s):=\sigma(s_k^i)\,,
\end{array}
$$ 
where $i$ is the smallest index such that $s\le s_k^i$. Since $\sum_i\HH(p(s_k^i)-p(s_k^{i-1}))\le \D_\HH(p;0,t)$, 
iterating the last inequality for $1\le i\le k$ we obtain
\begin{equation}\label{eq:155}
\begin{array}{c}
\QQ(e(0))-\langle \varrho(0)| e(0)-Ew(0)\rangle -\langle\varrho_D(0)|p(0) \rangle\le 
\vspace{.1cm} \\
\le \QQ(e(t)) + \D_\HH(p;0,t)-\langle \varrho(t)| e(t)-Ew(t)\rangle -\langle\varrho_D(t)|p(t) \rangle + {}
\\
\displaystyle +\int_0^t\langle\dot\varrho(s)|\ol e_k(s)-E\ol w_k(s)\rangle\,ds
+\int_0^t\langle\dot\varrho_D(s)|\ol p_k(s) \rangle\,ds -\int_{0}^{t}\langle \ol \sigma_k(s)|E\dot w(s)
\rangle\,ds +\delta_k\,,
\end{array}
\end{equation}
where $\delta_k:=\omega_k\int_0^T\| E\dot w(s)\|_2\,ds$.
By Remark~\ref{contqo} the set of discontinuity points of the functions
$s\mapsto p(s)$, $s\mapsto e(s)$, and $s\mapsto\sigma(s)$ is at most countable and 
$\|\ol p_k(s)\|_1$, $ \|\ol e_k(s)\|_2$, and $ \|\ol \sigma_k(s)\|_2$ are bounded uniformly 
with respect to $s$ and $k$. Therefore (\ref{fine}) implies that $\ol p_k(s)\to p(s)$ strongly
in $M_b(\Om\cup\Ga_0;\MD)$, $\ol e_k(s)\to e(s)$ and $\ol\sigma_k(s)\to \sigma(s)$ strongly in 
$L^2(\Om;\Mnn)$ for a.e.\ $s\in[0,t]$. Now, (\ref{en-ineq2}) follows from 
(\ref{eq:155}) by the dominated convergence theorem.
\end{proof}

\subsection{Convergence of the approximate solutions}
For every $k$ let $(\uki,e_k^i,\pki)$,
$i=1,\ldots,k$, be defined inductively as solutions of the discrete 
problems (\ref{min-inc}), starting from $(u_k^0,e_k^0, p_k^0) =(u_0,e_0,p_0)$, and let 
$u_k(t)$, $e_k(t)$, $p_k(t)$, $\sigma_k(t)$ be defined by (\ref{ukt}).
Let $t\mapsto(u(t),e(t),p(t))$ be a quasistatic evolution.
Assume that
\begin{equation}\label{convpk}
p_k(t)\wto p(t) \qquad\hbox{weakly}^* \hbox{ in } M_b(\Om\cup\Ga_0;\MD)
\end{equation}
for every $t\in[0,T]$. The following theorem shows, in particular, that stresses 
and elastic strains of the approximate solutions converge strongly in $L^2(\Om;\Mnn)$.

\begin{theorem}\label{convstress}
Assume that the plastic strain of the approximate solutions satisfies (\ref{convpk}). 
Then $e_k(t)\to e(t)$ and $\sigma_k(t)\to\sigma(t)$ strongly in $L^2(\Om;\Mnn)$. Moreover,
\begin{equation}\label{convdiss}
\begin{array}{c}
\displaystyle
\lim_{k\to\infty}\sum_{0<t_k^r\le t}\{\HH(p_k^r-p_k^{r-1}) -\langle\varrho_D(t_k^r)|p_k^r-p_k^{r-1}\rangle\} =
\\
\displaystyle
=\D_\HH(p;0,t) -\langle\varrho_D(t)| p(t)\rangle +\langle\varrho_D(0)| p(0)\rangle
+\int_0^t\langle\dot\varrho_D(s)|p(s)\rangle\,ds
\end{array}
\end{equation}
for every $t\in[0,T]$.
\end{theorem}

\begin{proof}
By the discrete energy inequality (\ref{inc-ineq}) for every $t\in[0,T]$ we have
\begin{equation}\label{ak}
\begin{array}{c}
\displaystyle
\QQ(e_k(t))+\sum_{0<t_k^r\le t}\{\HH(p_k^r-p_k^{r-1})
-\langle \varrho_D(t_k^r)|p_k^r-p_k^{r-1} \rangle\} \le 
\vspace{.2cm}
\\
\le \QQ(e_0)-\langle\varrho(0)|e_0-Ew(0)\rangle +\langle\varrho_k(t)|e_k(t)-Ew_k(t)\rangle -{}
\vspace{.2cm}
\\
\displaystyle
{}-\int_0^{\tki} \langle\dot\varrho(s)|e_k(s)-Ew_k(s)\rangle\,ds 
+\int_0^{\tki}\langle\sigma_k(s)|E\dot w(s)\rangle\, ds +\delta_k\,,
\end{array}
\end{equation}
where $\delta_k\to 0$ and $i$ is the largest integer such that $t_k^i\le t$.
By the energy balance (\ref{en-eq-v}) we have also 
\begin{equation}\label{a}
\begin{array}{c}
\displaystyle
\QQ(e(t)) + \D_\HH(p; 0,t) -\langle\varrho_D(t)| p(t)\rangle +\langle\varrho_D(0)| p(0)\rangle
+\int_0^t\langle\dot\varrho_D(s)|p(s)\rangle\,ds
= \vspace{.25cm}
\\
\displaystyle =\QQ(e_0) -\langle \varrho(0)|e_0-Ew(0)\rangle 
+\langle \varrho(t)|e(t)-Ew(t)\rangle  -{} \vspace{.25cm}
\\
\displaystyle 
- \int_0^t \langle \dot\varrho(s)|e(s)-Ew(s)\rangle\,ds + \int_0^t \langle \sigma(s)|E\dot w(s)\rangle\, ds
\,.
\end{array}
\end{equation}
In the proof of Theorem~\ref{main} we have already seen  that $e_k(t)\wto e(t)$ and $\sigma_k(t)\wto\sigma(t)$ 
weakly in $L^2(\Om;\Mnn)$, and that $\|e_k(t)\|_2$ and  $\|\sigma_k(t)\|_2$ are bounded uniformly with respect to $t$ and~$k$. 
Moreover, $\varrho_k(t)\to \varrho(t)$ and $Ew_k(t)\to Ew(t)$ 
strongly in $L^2(\Om;\Mnn)$.
Therefore the right-hand side of (\ref{ak}) converges to the right-hand side of (\ref{a}). This implies
\begin{eqnarray*}
& \displaystyle
\limsup_{k\to\infty} \Big\{ \QQ(e_k(t))+ \sum_{0<t_k^r\le t}\{\HH(p_k^r-p_k^{r-1})
-\langle \varrho_D(t_k^r)|p_k^r-p_k^{r-1} \rangle\}\Big\} \le 
\\
&\displaystyle
\le \QQ(e(t))+\D_\HH(p;0,t)-\langle\varrho_D(t)| p(t)\rangle 
+\langle\varrho_D(0)| p(0)\rangle
+\int_0^t\langle\dot\varrho_D(s)|p(s)\rangle\,ds\,.
\end{eqnarray*}
By the lower semicontinuity of $\QQ$  and by (\ref{wdineq}) we obtain (\ref{convdiss}) and
$$
\QQ(e_k(t))\to  \QQ(e(t))\,,
$$
which gives the strong convergence of $e_k(t)$, and, consequently, of ${\sigma_k(t)=\C e_k(t)}$.
\end{proof}

\end{section}

\begin{section}{Regularity and uniqueness results}

In this section we prove that every quasistatic evolution $t\mapsto (u(t),e(t),p(t))$ is absolutely continuous with respect to time, and that the functions $t\mapsto e(t)$ and
$t\mapsto \sigma(t)$ are uniquely determined by their initial conditions. 

\subsection{Regularity}
For the general properties of absolutely continuous functions with values in Banach spaces we refer to \cite[Appendix]{Bre} for the reflexive case and 
to the Appendix of the present paper for the case of the dual 
of a separable Banach space.

If $t\mapsto q(t)$ and $t\mapsto v(t)$  are absolutely  continuous from $[0,T]$ into
$M_b(\Om\cup\Ga_0;\MD)$ and $BD(\Om)$, respectively, we define
\begin{equation}\label{dots}
\dot q(t):=w^*\hbox{-}\lim_{s\to t}\frac{q(s)-q(t)}{s-t}\,, \qquad
\dot v(t):=w^*\hbox{-}\lim_{s\to t}\frac{v(s)-v(t)}{s-t}\,.
\end{equation}
By Theorem~\ref{thm:a01} $\dot q(t)$ and $\dot v(t)$ are defined for a.e.\ $t\in[0,T]$, 
the function $t\mapsto \HH(\dot q(t))$ is measurable, and
\begin{equation}\label{Dder01}
\D_\HH(q;0,t)=\int_0^t \HH(\dot q(s))\, ds
\end{equation}
for every $t\in[0,T]$. 

\begin{remark}\label{traceder}
If we apply (\ref{metder}) to the absolutely continuous function $t\mapsto q(t)$, with $X= M_b(\Om\cup\Ga_0;\MD)$, $Y=C_0(\Om\cup\Ga_0;\MD)$, and $\K=\{\varphi\in C_0(\Om\cup\Ga_0;\MD): \|\varphi\|_\infty\le 1\}$, for a.e.\ $t\in[0,T]$ we obtain 
\begin{equation}\label{strictp}
 \|\dot q(t)\|_1 = \lim_{s\to t} \, \Big\|\frac{q(s)-q(t)}{s-t}\Big\|_1 \,.
\end{equation}

By the definition of weak$^*$ convergence in $BD(\Om)$ it follows from (\ref{dots}) that 
for a.e.\ $t\in[0,T]$ we have $(v(s)-v(t))/(s-t)\to \dot v(t)$ strongly in $L^1(\Om;\Rn)$ and  ${(Ev(s)-Ev(t))/(s-t)}\allowbreak\wto E\dot v(t)$ weakly$^*$ in $M_b(\Om;\Mnn)$ as $s\to t$.
If we apply  (\ref{metder}) to the absolutely continuous function $t\mapsto Ev(t)$, with $X= M_b(\Om;\Mnn)$, $Y=C_0(\Om;\Mnn)$, and $\K=\{\varphi\in C_0(\Om;\Mnn): \|\varphi\|_\infty\le 1\}$, for a.e.\ $t\in[0,T]$ we obtain 
$$
\|E\dot v(t)\|_1=  \lim_{s\to t} \, \Big\|\frac{Ev(s)-Ev(t)}{s-t} \Big\|_1  \,.
$$
 This implies that for a.e.\ $t\in[0,T]$ the trace of $\dot v(t)$ is the strong limit in $L^1(\partial\Om;\Rn)$ of the traces of $(v(s)-v(t))/(s-t)$ as $s\to t$ (see \cite[Chapter~II, Theorem~3.1]{Tem}). In other words the time derivative of the trace of $v(t)$ is the trace of the time derivative of $v(t)$. Therefore, using (\ref{defLL}) and (\ref{LLdot2}), we can prove by a  standard argument that
\begin{equation}\label{Leibnitz}
\frac{d}{dt} \langle \LL(t)|v(t)\rangle = \langle \dot \LL(t)|v(t)\rangle + \langle \LL(t)|\dot v(t)\rangle
\end{equation}
for a.e.\ $t\in[0,T]$.
\end{remark}

The next proposition deals with the absolute continuity of the functions 
$t\mapsto e(t)$, $t\mapsto p(t)$, and $t\mapsto u(t)$ from $[0,T]$ into $L^2(\Om;\Mnn)$,
$M_b(\Om\cup\Ga_0;\MD)$, and $BD(\Om)$, respectively.

\begin{theorem}\label{ACesigma}
Let $t\mapsto (u(t),e(t),p(t))$ be a quasistatic evolution. 
Then the functions $t\mapsto e(t)$, $t\mapsto p(t)$, and 
$t\mapsto u(t)$ are absolutely continuous from $[0,T]$ 
into $L^2(\Om;\Mnn)$,  $M_b(\Om\cup\Ga_0;\MD)$, and $BD(\Om)$, 
respectively. Moreover, for a.e.\ $t\in[0,T]$ we have
\begin{eqnarray}\label{edotEwdot}
& \|\dot e(t)\|_2\le C_1(\|\dot\varrho(t)\|_2 +\|\dot\varrho_D(t)\|_\infty +\|E\dot w(t)\|_2)\,,
\\
& \|\dot p(t)\|_1\le C_2(\|\dot\varrho(t)\|_2 +\|\dot\varrho_D(t)\|_\infty +\|E\dot w(t)\|_2)\,,\label{pdot}
\\
& \|E\dot u(t)\|_1\le C_3(\|\dot\varrho(t)\|_2 +\|\dot\varrho_D(t)\|_\infty +\|E\dot w(t)\|_2)\,,\label{Eudot}
\\
& \|\dot u(t)\|_1\le C_4(\|\dot\varrho(t)\|_2 +\|\dot\varrho_D(t)\|_\infty +\|E\dot w(t)\|_2+\|\dot w(t)\|_2)\,,\label{udot}
\end{eqnarray}
where $C_1$ and $C_2$ are  positive constants depending on $R_K$, $\alpha_\C$, $\beta_\C$, $\alpha$, $\sup_{t}\|\varrho(t)\|_2$, 
 $\sup_{t}\|e(t)\|_2$, and $\sup_{t}\|p(t)\|_1$, while $C_3$ depends also on $\Om$ and
$C_4$ also on $\Om$ and $\Ga_0$.
\end{theorem}

\begin{proof}
Since $\HH(p(t_2)-p(t_1))\le \D_\HH(p;t_1,t_2)$,
by the energy equality (\ref{en-eq-v}) we obtain, after an integration by parts,
\begin{equation}\label{t1t2}
\begin{array}{c}
\displaystyle\vphantom{ \int_{t_1}^{t_2} }
{\textstyle\frac12}\langle \sigma(t_2)|e(t_2) \rangle - {\textstyle\frac12}\langle \sigma(t_1)|e(t_1) \rangle
+ \HH(p(t_2)-p(t_1))\le 
\\
\displaystyle\vphantom{ \int_{t_1}^{t_2} }
\le \langle\varrho(t_2)|e(t_2)\rangle -\langle\varrho(t_1)|e(t_1)\rangle
+ \langle\varrho_D(t_2)|p(t_2)\rangle -\langle\varrho_D(t_1)|p(t_1)\rangle - {}
\\
\displaystyle - \int_{t_1}^{t_2}\{\langle\dot\varrho(s)|e(s)\rangle + \langle\dot\varrho_D(s)|p(s)\rangle
-\langle \sigma(s)-\varrho(s)|E\dot w(s) \rangle\}\,ds
\end{array}
\end{equation}
for every $t_1$, $t_2\in [0,T]$ with $t_1<t_2$. Consider now the functions
$v:=u(t_2)-u(t_1)-(w(t_2)-w(t_1))$, $\eta:=e(t_2)-e(t_1)-(Ew(t_2)-Ew(t_1))$, and the measure
$q:=p(t_2)-p(t_1)$. Since $(v,\eta,q)\in A(0)$ and $(u(t_1),e(t_1),p(t_1))$ is a solution of the minimum problem (\ref{minp0}) 
with $p_0=p(t_1)$ and $\LL=\LL(t_1)$, by Theorem~\ref{Euler} and Lemma~\ref{fg} we obtain
\begin{eqnarray*}
& \displaystyle {}-\langle \sigma(t_1)|e(t_2)-e(t_1)\rangle + 
\langle \varrho(t_1)|e(t_2)-e(t_1)\rangle + \langle \varrho_D(t_1)|p(t_2)-p(t_1)\rangle + {} \vspace{.1cm}\\
& {}+\langle \sigma(t_1)-\varrho(t_1)|Ew(t_2)-Ew(t_1) \rangle 
\le \HH(p(t_2)-p(t_1))\,,
\end{eqnarray*}
so that (\ref{t1t2}) implies
\begin{eqnarray*}
& \displaystyle\vphantom{ \int_{t_1}^{t_2} }
{\textstyle\frac12} \langle \sigma(t_2)|e(t_2) \rangle -  {\textstyle\frac12}
\langle \sigma(t_1)|e(t_1) \rangle
-\langle \sigma(t_1)|e(t_2)-e(t_1) \rangle \le \langle\varrho(t_2)-\varrho(t_1)|e(t_2)\rangle + {}
\\ 
& \displaystyle\vphantom{ \int_{t_1}^{t_2} }
{}+ \langle\varrho_D(t_2)-\varrho_D(t_1)|p(t_2)\rangle -
\langle \sigma(t_1)-\varrho(t_1)|Ew(t_2)-Ew(t_1) \rangle - {} \\
& \displaystyle
{} - \int_{t_1}^{t_2}\{\langle\dot\varrho(s)|e(s)\rangle + \langle\dot\varrho_D(s)|p(s)\rangle
-\langle \sigma(s)-\varrho(s)|E\dot w(s) \rangle\}\,ds\,.
\end{eqnarray*}
Therefore,
\begin{eqnarray*}
& \displaystyle
{\textstyle\frac12}
\langle \C(e(t_2)-e(t_1))|e(t_2)- e(t_1)\rangle
\le  \int_{t_1}^{t_2}\langle \sigma(s)-\sigma(t_1)|E\dot w(s) \rangle\,ds + {}
\\
& \displaystyle {}+\int_{t_1}^{t_2}
\{\langle \dot\varrho(s)|e(t_2)-e(s)\rangle + \langle\dot\varrho_D(s)|p(t_2)-p(s)\rangle
-\langle \varrho(s)-\varrho(t_1) |E\dot w(s) \rangle\}
\,ds \,.
\end{eqnarray*}
By (\ref{boundsC}) and (\ref{normC}) we obtain
\begin{equation}\label{squared}
\begin{array}{c}
\displaystyle
\alpha_\C \|e(t_2)-e(t_1)\|^2_2 
\le 2\beta_\C \int_{t_1}^{t_2} \| e(s)- e(t_1)\|_2\, \|E\dot w(s)\|_2 \,ds + {}
\\
\displaystyle {}+\int_{t_1}^{t_2}
\|\dot\varrho(s)\|_2\,\|e(t_2)-e(s)\|_2 \,ds + \int_{t_1}^{t_2}\|\dot\varrho_D(s)\|_\infty\, \|p(t_2)-p(s)\|_1
\,ds +{} 
\\
\displaystyle
{} +\int_{t_1}^{t_2}\|\varrho(s)-\varrho(t_1)\|_2\,\|E\dot w(s)\|_2
\,ds \,.
\end{array}
\end{equation}
By Lemma~\ref{eps*} we have that for every $t_1\le s\le t_2$
$$
\alpha\| p(t_2)-p(s)\|_1 \le \HH(p(t_2)-p(s))- \langle\varrho_D(t_2)|p(t_2)-p(s)\rangle\,,
$$
therefore, inequality (\ref{t1t2}) with $t_1=s$ implies
\begin{eqnarray*}
& \displaystyle \vphantom{\int_{t_1}^{t_2}}
\alpha\| p(t_2)-p(s)\|_1 \le {\textstyle\frac12}\langle \sigma(s)|e(s) \rangle -
{\textstyle\frac12}\langle \sigma(t_2)|e(t_2) \rangle +{} \\
& \displaystyle \vphantom{\int_s^{t_2}}
{} +\langle\varrho(t_2)|e(t_2)-e(s)\rangle +\langle\varrho(t_2)-\varrho(s)|e(s)\rangle 
+\langle\varrho_D(t_2)-\varrho_D(s)|p(s)\rangle -{} \\
& \displaystyle
{} -\int_{s}^{t_2}\{\langle\dot\varrho(t)|e(t)\rangle + \langle\dot\varrho_D(t)|p(t)
\rangle-
\langle \sigma(t)-\varrho(t)|E\dot w(t) \rangle\}\, dt\,.
\end{eqnarray*}
We observe that $\sup_{t}\|\varrho(t)\|_2$, $\sup_{t}\|\varrho_D(t)\|_\infty$, $\sup_{t}\|e(t)\|_2$, and $\sup_{t}\|p(t)\|_1$ are finite (see Remark~\ref{contqo} for $e(t)$).
In the rest of the proof $C$ will denote a positive constant, whose value can change from line to line, depending on these suprema and on  the constants $\alpha_\C$, $\beta_\C$, $\alpha$.
The previous inequality implies that
\begin{eqnarray*}
& \displaystyle\vphantom{\int_s^{t_2}}
\| p(t_2)-p(s)\|_1 \le C( \|e(t_2)-e(s)\|_2 + \|\varrho(t_2)-\varrho(s)\|_2
+ \|\varrho_D(t_2)-\varrho_D(s)\|_\infty) +{}\\
& \displaystyle
{}+ C\int_s^{t_2} \{\|\dot\varrho(t)\|_2+ \|\dot\varrho_D(t)\|_\infty+ \|E\dot w(t)\|_2\}\,dt\,.
\end{eqnarray*}
Therefore, for every $t_1\le s\le t_2$
\begin{equation}\label{pAC}
\| p(t_2)-p(s)\|_1 \le C\, \|e(t_2)-e(s)\|_2 + 
C\int_{t_1}^{t_2} \{\|\dot\varrho(t)\|_2+ \|\dot\varrho_D(t)\|_\infty+ \|E\dot w(t)\|_2\}\,dt\,.
\end{equation}
By (\ref{squared}) and (\ref{pAC}), using $\|e(t_2)-e(s)\|_2\le \|e(t_2)-e(t_1)\|_2 + \|e(s)-e(t_1)\|_2$, we deduce that
\begin{eqnarray*}
& \displaystyle 
\|e(t_2)-e(t_1)\|_2^2 \le 
C\,\|e(t_2)-e(t_1)\|_2
\int_{t_1}^{t_2}\{\|\dot\varrho(s)\|_2+ \|\dot\varrho_D(s)\|_\infty\}\,ds +{} \\
& \displaystyle {}+ C \int_{t_1}^{t_2}\{\|\dot\varrho(s)\|_2+ \|\dot\varrho_D(s)\|_\infty+\|E\dot w(s)\|_2
\}\,\|e(s)-e(t_1)\|_2\,ds +{} \\
& \displaystyle 
{}+ C \Big(\int_{t_1}^{t_2}\{\|\dot\varrho(s)\|_2+ \|\dot\varrho_D(s)\|_\infty+\|E\dot w(s)\|_2
\}\,ds\Big)^2 \,.
\end{eqnarray*}
By the Cauchy inequality we have
$$
\|e(t_2)-e(t_1)\|_2^2 \le \int_{t_1}^{t_2}\psi(s)\, \|e(s)-e(t_1)\|_2\,ds
+ \Big(\int_{t_1}^{t_2}\psi(s)\,ds \Big)^2,
$$
where
$$
\psi(s):= C(\|\dot\varrho(s)\|_2+ \|\dot\varrho_D(s)\|_\infty+ \|E\dot w(s)\|_2)\,.
$$
We can apply now a version of Gronwall inequality, proved in Lemma~\ref{Gronwall} below, which gives
\begin{equation}\label{abcd}
\|e(t_2)-e(t_1)\|_2\le 
\frac32 \int_{t_1}^{t_2}\psi(s)\, ds
\le C \int_{t_1}^{t_2} \{
\|\dot\varrho(s)\|_2+ \|\dot\varrho_D(s)\|_\infty+ \|E\dot w(s)\|_2\}
\,ds\,.
\end{equation}
This implies that $t\mapsto e(t)$ is absolutely continuous from $[0,T]$
into $L^2(\Om;\Mnn)$ and that $\dot e(t)$ satisfies (\ref{edotEwdot}).

Using the absolute continuity of $t\mapsto e(t)$ and (\ref{edotEwdot}), inequality (\ref{pAC}) with $s=t_1$
yields the absolute continuity of ${t\mapsto p(t)}$ and (\ref{pdot}).

{}From the decomposition $Eu(t)=e(t)+p(t)$, it follows that $t\mapsto Eu(t)$  is absolutely continuous from $[0,T]$ into $M_b(\Om;\Mnn)$
and $E\dot u(t)=\dot e(t)+\dot p(t)$ for a.e.\ $t\in[0,T]$. Inequality (\ref{Eudot}) is an easy consequence of this
decomposition.
It remains to prove that $t\mapsto u(t)$ is absolutely continuous from $[0,T]$ into $L^1(\Om;\Rn)$ and satisfies (\ref{udot}).
By (\ref{seminorm}) there exists a constant $C>0$, 
depending on $\Om$ and $\Ga_0$, such that
\begin{equation}
\label{seminorm1}
\|u(t_2)-u(t_1)\|_1 \le C\, \|u(t_2)-u(t_1)\|_{1,\Ga_0}+ 
C\, \|Eu(t_2)- Eu(t_1)\|_1\,.
\end{equation}
Using (\ref{901}) and the continuity of the trace operator from $H^1(\Om;\Rn)$ into $L^1(\partial\Om;\Rn)$, we obtain
that there exists a constant $M$, depending on $\Om$, such that
\begin{equation}\label{Ga0}
\begin{array}{c}
\|u(t_2)-u(t_1)\|_{1,\Ga_0}\le 
\sqrt2\, \|p(t_2)-p(t_1)\|_1 +{}
\\
{}+
M\|w(t_2)-w(t_1)\|_2+M\|Ew(t_2)-Ew(t_1)\|_2\,.
\end{array}
\end{equation}
As $t\mapsto w(t)$, $t\mapsto Eu(t)$, and $t\mapsto p(t)$ are absolutely continuous from $[0,T]$ into 
$H^1(\Om;\Rn)$, $M_b(\Om;\Mnn)$, and $M_b(\Om\cup\Ga_0;\MD)$, respectively, 
inequalities (\ref{seminorm1}) and (\ref{Ga0}) imply that
$t\mapsto u(t)$  is absolutely continuous from $[0,T]$ into $L^1(\Om;\Rn)$ and (\ref{udot}) is satisfied.
\end{proof}

\begin{lemma}\label{Gronwall}
Let $\phi\colon[0,T]\to[0,+\infty[$ be a bounded measurable function and let $\psi\colon[0,T]\to[0,+\infty[$ be an integrable function.
Suppose that
\begin{equation}\label{hypothesis}
\phi(t)^2\le \int_0^t \phi(s)\, \psi(s)\,ds + \Big( \int_0^t\psi(s)\,ds \Big)^2
\end{equation}
for every $t\in[0,T]$.
Then
\begin{equation}\label{thesis}
\phi(t)\le\frac32 \int_0^t\psi(s)\,ds
\end{equation}
for every $t\in[0,T]$.
\end{lemma}

\begin{proof}
Let us fix $t_0\in [0,T]$ and let 
$\gamma_0:= ( \int_0^{t_0}\psi(s)\,ds)^2$.
For every $t\in[0,t_0]$ we define
$V(t):=  \int_0^t \phi(s)\, \psi(s)\,ds$.
Then $V$ is absolutely continuous on $[0,t_0]$,
$\phi(t)^2 \le V(t)+\gamma_0$ for every $t\in [0,t_0]$, and $\dot V(t)\le \psi(t) ({V(t)+\gamma_0})^{1/2}$ for a.e.\ $t\in [0,t_0]$. Integrating between $0$ and $t_0$ we get $2({V(t_0)+\gamma_0})^{1/2}\allowbreak\le 2\gamma_0^{1/2} +  \int_0^{t_0}\psi(s)\,ds=  3  \int_0^{t_0}\psi(s)\,ds$. By (\ref{hypothesis}) we have $\phi(t_0)\le ({V(t_0)+\gamma_0})^{1/2}$, so that the previous inequality gives $2\,\phi(t_0)\le 3  \int_0^{t_0}\psi(s)\,ds$.
\end{proof}

\begin{remark}\label{Lipschitz}
Estimates (\ref{edotEwdot})--(\ref{udot}) imply that, if $t\mapsto w(t)$,  $t\mapsto \varrho(t)$, and $t\mapsto \varrho_D(t)$ are  Lipschitz continuous from $[0,T]$ into $H^1(\Rn;\Rn)$, $L^2(\Om;\Mnn)$, and $L^\infty(\Om;\MD)$, respectively, then the functions $t\mapsto u(t)$, $t\mapsto e(t)$, $t\mapsto p(t)$ are  Lipschitz continuous from $[0,T]$ into $BD(\Om)$, $L^2(\Om;\Mnn)$,
and $M_b(\Om\cup\Ga_0;\MD)$, respectively.
\end{remark}

The following lemma will be crucial in the rest of the paper.

\begin{lemma}\label{dotclosure}
Let $t\mapsto u(t)$, $t\mapsto e(t)$, $t\mapsto p(t)$  be absolutely continuous functions from $[0,T]$ into $BD(\Om)$, $L^2(\Om;\Rn)$,
and $M_b(\Om\cup\Ga_0;\MD)$, respectively.
Assume that $(u(t),e(t),p(t))\in A(w(t))$ for every $t\in[0,T]$.
Then $(\dot u(t),\dot e(t),\dot p(t))\in A(\dot w(t))$ for a.e.\ $t\in[0,T]$.
\end{lemma}

\begin{proof}
It is enough to apply Lemma~\ref{lm:145} to the difference quotients.
\end{proof}

Thanks to the following proposition, we can differentiate the energy balance (\ref{en-eq}) and obtain a balance of powers: 
the rate of change of stored energy plus the rate of plastic dissipation equals the power of external forces.

\begin{proposition}\label{equivalence}
Let $t\mapsto (u(t),e(t),p(t))$ be an absolutely continuous function from $[0,T]$ into 
$BD(\Om){\times}L^2(\Om;\Mnn){\times}M_b(\Om\cup\Ga_0;\MD)$ and let $\sigma(t):=\C e(t)$. 
Then the following conditions are equivalent:
\begin{itemize}
\smallskip
\item[{\rm (a)}] for every $t\in[0,T]$
\begin{eqnarray*}
&\displaystyle \QQ(e(t)) + \D_\HH(p; 0,t) -\langle\LL(t)|u(t)\rangle=
\QQ(e(0))-\langle\LL(0)|u(0)\rangle +{} \\
&\displaystyle {}+\int_0^t \{
\langle \sigma(s)|E\dot w(s)\rangle -\langle\LL(s)|\dot w(s)\rangle-\langle\dot\LL(s)|u(s)\rangle
\}\, ds\,;
\end{eqnarray*}
\item[{\rm (b)}] for a.e.\ $t\in[0,T]$
$$
\langle\sigma(t)|\dot e(t)\rangle +\HH(\dot p(t))=
\langle\sigma(t)|E\dot w(t)\rangle -\langle\LL(t)|\dot w(t)\rangle
+\langle\LL(t)|\dot u(t)\rangle\,;
$$ 
\item[{\rm (c)}] for a.e.\ $t\in[0,T]$
$$
\langle\sigma(t)-\varrho(t)|\dot e(t)\rangle +\HH(\dot p(t))=
\langle\varrho_D(t)|\dot p(t)\rangle+ 
\langle\sigma(t)-\varrho(t)|E\dot w(t)\rangle\,;
$$
\item[{\rm (d)}] for every $t\in[0,T]$
$$
\begin{array}{l}
\displaystyle \QQ(e(t)) +\int_0^t \{ \HH(\dot p(s)) -
\langle\varrho_D(s)|\dot p(s)\rangle\}\,ds =
\\
\displaystyle \qquad=\QQ(e(0))+ \int_0^t \{
\langle \varrho(s)|\dot e(s)\rangle +
\langle \sigma(s)-\varrho(s)|E\dot w(s)\rangle
\}\, ds\,.
\end{array}
$$
\end{itemize}
\end{proposition}

\begin{proof}
Using (\ref{Dder01}) and (\ref{Leibnitz}) we obtain (b) by differentiating (a) and (a) by integrating~(b). Similarly we obtain (d) by integrating~(c) and (c) by differentiating (d).
The equivalence between (b) and (c) follows from Lemmas~\ref{fg} and~\ref{dotclosure}.
\end{proof}

Condition (d) of Proposition~\ref{equivalence} allows us to prove an estimate of $\sup_{t} \|e(t)\|_2$ and $\sup_{t} \|p(t)\|_1$ in terms of the data of the problem.

\begin{proposition}\label{bounded}
Let $t\mapsto (u(t),e(t),p(t))$ be a quasistatic evolution. Then
\begin{equation}\label{bounde}
\begin{array}{c}
\displaystyle
\sup_{t\in[0,T]} \|e(t)\|_2 \,\,\le\,\, C_1\,\Big\{  \|e(0)\|_2+ \sup_{t\in[0,T]} \|\varrho(t)\|_2 +{}
\\
\displaystyle
{}+  \int_0^T \|\dot \varrho(t)\|_2\, dt + 
\int_0^T \|E\dot w(t)\|_2 \,dt\,\Big\} \,,
\end{array}
\end{equation}
and
\begin{equation}\label{boundp}
\begin{array}{c}
\displaystyle
\sup_{t\in[0,T]} \|p(t)\|_1 \,\,\le\,\, \|p(0)\|_1 + C_2\,\Big\{  \|e(0)\|_2^2+  \sup_{t\in[0,T]} \|\varrho(t)\|_2^2 +{}
\\
\displaystyle
 {}+  \Big(\int_0^T \|\dot \varrho(t)\|_2\, dt\Big)^2 + 
\Big( \int_0^T \|E\dot w(t)\|_2 \,dt \Big)^2  \,\Big\}  \,,
 \end{array}
\end{equation}
where $C_1$ is a positive constant depending only on $\alpha_\C$ and $\beta_\C$, while $C_2$ depends also on~$\alpha$.
\end{proposition}

\begin{proof}
By Theorem~\ref{ACesigma} the function $t\mapsto (u(t),e(t),p(t))$ is absolutely continuous from $[0,T]$ into 
$BD(\Om){\times}L^2(\Om;\Mnn){\times}M_b(\Om\cup\Ga_0;\MD)$. As $t\mapsto (u(t),e(t),p(t))$ satisfies (qs2) in Definition~\ref{qstatic}, it satisfies conditions (a) and (d) of Proposition~\ref{equivalence}. After an integration by parts, we obtain from (d)
$$
\begin{array}{l}
\displaystyle \QQ(e(t)) +\int_0^t \{ \HH(\dot p(s)) -
\langle\varrho_D(s)|\dot p(s)\rangle\}\,ds -\langle\varrho(t)| e(t)\rangle=
\\
\displaystyle\qquad = \QQ(e(0))+ \int_0^t \{
\langle \sigma(s)-\varrho(s)|E\dot w(s)\rangle
-\langle \dot\varrho(s)|e(s)\rangle
\}\, ds
-\langle\varrho(0)| e(0)\rangle\,.
\end{array}
$$
By (\ref{boundsC}), (\ref{normC}), and (\ref{eps*1}) for every $t\in [0,T]$ we have
\begin{equation}\label{xyz}
\begin{array}{c}
\displaystyle
\alpha_\C \, \|e(t)\|_2^2 + \alpha \int_0^t \|\dot p(s)\|_1\, ds \,\,\le\,\, \beta_\C \, \|e(0)\|_2^2 +
2\,
\sup_{t\in[0,T]} \|\varrho(t)\|_2 \,
\sup_{t\in[0,T]} \|e(t)\|_2 +{}
\\
\displaystyle
{}+ \sup_{t\in[0,T]} \|e(t)\|_2
\int_0^T \{ 2\,\beta_\C\,  \| E\dot w(s))\|_2 +  \| \dot\varrho(s)\|_2 \}\, ds +{}
\\
\displaystyle
{}+ \ \sup_{t\in[0,T]} \|\varrho(t)\|_2 \int_0^T \| E\dot w(s))\|_2\, ds\,,
\end{array}
\end{equation}
which implies (\ref{bounde}) and  (\ref{boundp}) by the Cauchy inequality.
\end{proof}

\begin{remark}\label{goodestimates}
Let $t\mapsto (u(t),e(t),p(t))$ be a quasistatic evolution. By Proposition~\ref{bounded},
estimates (\ref{edotEwdot})--(\ref{udot}) are satisfied with constants $C_1,\dots,C_4$ depending only on the data of the problem. More precisely, $C_1$ and $C_2$ depend on $R_K$, $\alpha_\C$, $\beta_\C$, $\alpha$, $\sup_{t}\|\varrho(t)\|_2$, $\int_0^T\|\dot \varrho(t)\|_2\,dt$, $\int_0^T\|E\dot w(t)\|_2\,dt$,
 $\|e(0)\|_2$, and $\|p(0)\|_1$, while $C_3$ depends also on $\Om$, and
$C_4$ also on $\Om$ and $\Ga_0$.
\end{remark}

\subsection{Uniqueness of stress and elastic strain}

We now prove that $t\mapsto e(t)$ (and, consequently, $t\mapsto \sigma(t)$) is uniquely determined by its initial condition.

\begin{theorem}\label{uniqueness}
Let $t\mapsto (u(t),e(t),p(t))$ and $t\mapsto (v(t),\eta(t),q(t))$ 
be two quasistatic evolutions and let $\sigma(t):=\C e(t)$ and $\tau(t):=\C \eta(t)$. If $e(0)=\eta(0)$,
then $e(t)=\eta(t)$ for every $t\in[0,T]$. 
Equivalently, if $\sigma(0)=\tau(0)$, then ${\sigma(t)=\tau(t)}$ for every $t\in[0,T]$.
\end{theorem}

\begin{proof} By Theorem~\ref{ACesigma} the functions  $t\mapsto (u(t),e(t),p(t))$ and $t\mapsto (v(t),\eta(t),q(t))$ are absolutely continuous.
By condition (c) of Proposition~\ref{equivalence}
we have
\begin{eqnarray}
& \langle\sigma(t)-\varrho(t)|\dot e(t)-E\dot w(t)\rangle +\HH(\dot p(t))=
\langle\varrho_D(t)|\dot p(t)\rangle\,, \label{prima}\\
& \langle\tau(t)-\varrho(t)|\dot \eta(t)-E\dot w(t)\rangle +\HH(\dot q(t))=
\langle\varrho_D(t)|\dot q(t)\rangle\,.  \label{seconda}
\end{eqnarray}
{}From the global stability condition (\ref{min-ineq}) and from Theorem~\ref{cns} it follows that  for every $t\in[0,T]$ we have
$\tau(t)\in \Sigma(\Om)\cap\K(\Om)$, $-\div\,\tau(t)=f(t)$ a.e.\ on $\Om$, and $[\tau(t)\nu]=g(t)$ on $\Ga_1$.
By Lemma~\ref{dotclosure} we have $(\dot u(t),\dot e(t),\dot p(t))\in A(\dot w(t))$ for a.e.\ $t\in[0,T]$. Therefore 
Proposition~\ref{prp:star} gives $\HH(\dot p(t))\ge \langle \tau_D(t)|\dot p(t)\rangle$.
By (\ref{prima}) this implies
$$
\langle\sigma(t)-\varrho(t)|\dot e(t)-E\dot w(t)\rangle
+ \langle [\tau_D(t)-\varrho_D(t)]|\dot p(t)\rangle
\le 0\,.
$$
As $\div(\tau(t)-\varrho(t))=0$ a.e.\ on $\Om$ and $[(\tau(t)-\varrho(t))\nu]=0$
on $\Ga_1$ by (\ref{safeload}) and Theorem~\ref{cns}, this inequality is equivalent to
$$
\langle\sigma(t)-\tau(t)|\dot e(t)-E\dot w(t)\rangle
\le 0\,.
$$
in view of the integration by parts formula (\ref{intpartsp2}).
Analogously from (\ref{seconda}) we obtain
$$
\langle\tau(t)-\sigma(t)|\dot \eta(t)-E\dot w(t)\rangle
\le 0\,.
$$
Summing these two inequalities we get
$$
\langle\C(e(t)-\eta(t))|\dot e(t)-\dot \eta(t)\rangle\le 0\,,
$$
hence
$$
\frac{d}{dt}\langle\C(e(t)-\eta(t))|e(t)-\eta(t)\rangle
\le 0\,.
$$
If $e(0)=\eta(0)$, we have $\langle\C(e(0)-\eta(0))|e(0)-\eta(0)\rangle=0$, so that  for every $t\in[0,T]$
$\langle\C(e(t)-\eta(t))|e(t)-\eta(t)\rangle\le 0$, which is equivalent to $e(t)=\eta(t)$ by~(\ref{boundsC}).
\end{proof}

\end{section}

\begin{section}{Equivalent formulations in rate form}

Let $t\mapsto (u(t),e(t),p(t))$ be a quasistatic evolution.
Suppose for a moment that $\dot p(t)\in L^2(\Om;\MD)$, and denote the values of $\dot p(t)$ and $\sigma_D(t)$ at $x\in\Om$ 
by $\dot p(t,x)$ and $\sigma_D(t,x)$, respectively. 
We recall that the normal cone $N_K(\xi_0)$ to $K$ at $\xi_0\in \MD$ is defined in the following way: if $\xi_0\in K$, then $N_K(\xi_0)$ is the set of matrices $\zeta\in \MD$ 
such that $\zeta{\,:\,}(\xi-\xi_0)\le 0$ for every $\xi\in K$; if  $\xi_0\notin K$, then $N_K(\xi_0):=\emptyset$. In this section we want to prove that
\begin{equation}\label{normal}
\dot p(t,x)\in N_K(\sigma_D(t,x)) \qquad\hbox{for a.e.\ }x\in\Om\,,
\end{equation}
which represents the classical formulation of the flow rule. 

\subsection{Weak formulation}
By the definition of $N_K$ it is easy to see that (\ref{normal}) is equivalent to saying that
\begin{equation}\label{wnormal}
\langle \sigma_D(t)-\tau_D|\dot p(t)\rangle\ge 0
\end{equation}
for every $\tau\in \Sigma(\Om)\cap \K(\Om)$ with $[\tau\nu]=g(t)$ on $\Ga_1$. Indeed, the fact that (\ref{normal}) implies (\ref{wnormal}) is straightforward, 
while to prove the converse implication it is enough to consider test functions of the form
$\tau=\varphi\,\xi + (1-\varphi)\,\sigma$, with $\varphi\in C^\infty_c(\Om)$, 
$0\le\varphi\le 1$, and $\xi\in K$.

Note that the variational inequality (\ref{wnormal}) makes sense even if $\dot p(t)$ is only a measure, 
since in any case $\dot p(t)\in \Pi_{\Ga_0}(\Om)$ by Theorem~\ref{ACesigma} and Lemma~\ref{dotclosure}, 
so that the duality product $\langle{ \sigma_D(t)-\tau_D}|\dot p(t)\rangle$ is well defined by~(\ref{sigmap}). We will regard (\ref{wnormal}) as the weak formulation of inclusion
(\ref{normal}) when  $\dot p(t)\in M_b(\Om\cup\Ga_0;\MD)$.

The following theorem collects three different sets of conditions, including (\ref{wnormal}) and expressed in terms of the time 
derivatives $\dot p(t)$, $\dot e(t)$, and $\dot u(t)$, which are equivalent to the conditions considered in Definition~\ref{qstatic}.


\begin{theorem}\label{varineq}
Let $t\mapsto (u(t),e(t),p(t))$ be a function from $[0,T]$ into 
$BD(\Om){\times}L^2(\Om;\Mnn)\allowbreak{\times}M_b(\Om\cup\Ga_0;\MD)$ and let $\sigma(t):=\C e(t)$.
Then the following conditions are equivalent:
\begin{itemize}
\smallskip
   \item[(a)] $t\mapsto (u(t),e(t),p(t))$ is a quasistatic evolution;
   \smallskip
   \item[(b)] $t\mapsto (u(t),e(t),p(t))$ is absolutely continuous and
\begin{itemize}
        \item[(b1)] for every $t\in [0,T]$ we have $(u(t),e(t),p(t))\in A(w(t))$,
        $\sigma(t)\in\Sigma(\Om)\cap \K(\Om)$, $-\div\,\sigma(t)=f(t)$ a.e.\ on $\Om$, and $[\sigma(t)\nu]=g(t)$ on~$\Ga_1$,
        \item[(b2)]  for a.e.\ $t\in [0,T]$ we have 
	$$
	\HH(\dot p(t))= \langle\sigma_D(t) | \dot p(t)\rangle\,;
	$$ 
\end{itemize}
 \smallskip
    \item[(c)] $t\mapsto (u(t),e(t),p(t))$ is absolutely continuous and
\begin{itemize}
        \item[(c1)] for every $t\in [0,T]$ we have $(u(t),e(t),p(t))\in A(w(t))$,
        $\sigma(t)\in\Sigma(\Om)\cap \K(\Om)$, $-\div\,\sigma(t)=f(t)$ a.e.\ on $\Om$, and $[\sigma(t)\nu]=g(t)$ on~$\Ga_1$,
        \item[(c2)]  for a.e.\ $t\in [0,T]$ we have 
	$$
	\langle\sigma_D(t)-\tau_D | \dot p(t)\rangle\ge 0
	$$ 
	for every $\tau\in \Sigma(\Om)\cap \K(\Om)$ with $[\tau\nu]=g(t)$ on $\Ga_1$;
\end{itemize}
 \smallskip
   \item[(d)] $t\mapsto (u(t),e(t))$ is absolutely continuous and
\begin{itemize}
        \item[(d1)] for every $t\in [0,T]$ we have 
	$\sigma(t)\in\Sigma(\Om)\cap \K(\Om)$, $-\div\,\sigma(t)=f(t)$ a.e.\ on $\Om$, and $[\sigma(t)\nu]=g(t)$ on~$\Ga_1$,
	\item[(d2)] for a.e.\ $t\in [0,T]$ we have  
	$$
	\langle\tau-\sigma(t) | \dot e(t)\rangle +
	 \langle\div\,\tau -\div\,\sigma(t)| \dot u(t)\rangle\ge \langle [(\tau-\sigma(t))\nu] | \dot w(t)\rangle_{\partial\Om}
	$$
 	for every $\tau\in \Sigma(\Om)\cap \K(\Om)$ with $[\tau\nu]=g(t)$ on $\Ga_1$, where $\langle\cdot | \cdot\rangle_{\partial\Om}$ denotes the duality pairing between $H^{-1/2}(\partial\Om;\Rn)$ and $H^{1/2}(\partial\Om;\Rn)$,
  \item[(d3)] for every $t\in [0,T]$ $p(t)=Eu(t)-e(t)$ on $\Om$
  and  $p(t)= (w(t)-u(t)){\,\odot\,}\nu\,\hn$ on $\Ga_0$.
\end{itemize}
\end{itemize}
\end{theorem}

Note that in conditions (b) and (c) the duality products 
$\langle\sigma_D(t) | \dot p(t)\rangle$ and $\langle\sigma_D(t)-\tau_D | \dot p(t)\rangle$ are well defined by~(\ref{sigmap}), since 
$\dot p(t)\in \Pi_{\Ga_0}(\Om)$ by Lemma~\ref{dotclosure}, and $\sigma(t)$, $\tau \in \Sigma(\Om)$.
 
\begin{proof}[Proof of Theorem~\ref{varineq}] 
We first prove that ${\rm(a)} \Leftrightarrow {\rm(b)}$. 
We already proved in Theorem~\ref{ACesigma} that every quasistatic evolution is absolutely continuous. 
Moreover, Theorem~\ref{cns} shows that (b1) is equivalent to the global stability condition~(qs1) of Definition~\ref{qstatic}.  
By Proposition~\ref{equivalence} it only remains to prove that, for an absolutely continuous function $t\mapsto (u(t),e(t),p(t))$ 
satisfying either (b1) or (qs1), condition (b2) is equivalent to the balance of powers
\begin{equation}\label{balpow}
\langle\sigma(t)|\dot e(t)\rangle +\HH(\dot p(t))=
\langle\sigma(t)|E\dot w(t)\rangle -\langle\LL(t)|\dot w(t)\rangle +\langle\LL(t)|\dot u(t)\rangle
\end{equation}
for a.e.\ $t\in[0,T]$.
Since $ (\dot u(t),\dot e(t),\dot p(t))\in A(\dot w(t))$ for a.e.\ $t\in[0,T]$ by Lemma~\ref{dotclosure}, condition~(b2) is equivalent to~(\ref{balpow}) in view of the integration
by parts formula~(\ref{intpartsp2}).

We now prove that ${\rm(b)} \Leftrightarrow {\rm(c)}$. It is enough to show that, if (b1) is satisfied, 
then ${\rm(b2)} \Leftrightarrow {\rm(c2)}$. Condition (c2) is equivalent to
$$
\langle\sigma_D(t)|\dot p(t)\rangle \ge \sup \{ \langle\tau_D|\dot p(t)\rangle : \tau\in \Sigma(\Om)\cap\K(\Om),\ 
[\tau\nu]=g(t) \, \hbox{ on }\,\Ga_1 \}\,.
$$
Since $\sigma(t)\in \Sigma(\Om)\cap\K(\Om)$ and $[\sigma(t)\nu]=g(t)$ on $\Ga_1$ by (b1), 
the opposite inequality is trivial, so (c2) is equivalent to 
$$
\langle\sigma_D(t)|\dot p(t)\rangle = \sup \{ \langle\tau_D|\dot p(t)\rangle : \tau\in \Sigma(\Om)\cap\K(\Om),\ 
[\tau\nu]=g(t) \, \hbox{ on }\,\Ga_1 \}\,.
$$
This last condition is equivalent to (b2) by Proposition~\ref{prp:star}.

{}Finally, we prove that ${\rm(c)} \Leftrightarrow {\rm(d)}$. 
We observe first that (d3) and the absolute continuity of $t\mapsto (u(t),e(t))$ imply that also $t\mapsto p(t)$ 
is absolutely continuous and $(u(t),e(t),p(t))\in A(w(t))$  for every $t\in [0,T]$. 
It remains to prove that, if (c1) is satisfied, then  ${\rm(c2)} \Leftrightarrow {\rm(d2)}$.

By (\ref{sigmanu}) we have
$$
\langle [(\tau-\sigma(t))\nu] | \dot w(t) \rangle_{\partial\Om} = 
\langle\div\,\tau-\div\,\sigma(t) |\dot w(t)\rangle + \langle \tau-\sigma(t) | E \dot w(t)\rangle\,.
$$
Therefore (d2) is equivalent to
\begin{equation}\label{label10}
\langle \tau-\sigma(t) | \dot e(t)- E \dot w(t)\rangle
+\langle\div\,\tau-\div\,\sigma(t) |\dot u(t)-\dot w(t)\rangle \ge 0\,.
\end{equation}
Since $ (\dot u(t),\dot e(t),\dot p(t))\in A(\dot w(t))$ for a.e.\ $t\in[0,T]$ by Lemma~\ref{dotclosure} and $[(\tau-\sigma(t))\nu]=0$ on $\Ga_1$, condition (c2) is  equivalent to (\ref{label10}) thanks to
the integration
by parts formula~(\ref{intpartsp2}).
\end{proof}

\begin{remark}\label{maxdiss}
By Proposition~\ref{prp:star} the measure $H(\dot p(t))- [\sigma_D(t) {\,:\,} \dot p(t)]$ is nonnegative on $\Om\cup\Ga_0$, so that (b2) of Theorem~\ref{varineq} implies 
\begin{equation}\label{H=sigmap}
H(\dot p(t))= [\sigma_D(t) {\,:\,} \dot p(t)]\qquad\hbox{on }\Om\cup\Ga_0\,. 
\end{equation}
\end{remark}

\begin{remark}\label{Suquet}
Condition (d) of Theorem~\ref{varineq} is the weak formulation of the quasistatic evolution problem 
for perfectly plastic materials, proposed in \cite{Joh} in a slightly different form, and analysed in~\cite{Suq}.
\end{remark}

\subsection{Strong formulation and precise definition of the stress}
Let us return to the classical formulation
(\ref{normal}) of the flow rule, which makes sense if 
$\dot p(t)\in L^2(\Om;\MD)$. It can be written equivalently in the form
\begin{equation}\label{classform}
\frac{\dot p(t,x)}{|\dot p(t,x)|}\in N_K(\sigma_D(t,x)) 
\qquad\hbox{for }\,\Ln\hbox{-a.e.\ }x\in\{|\dot p(t)|>0\}\,.
\end{equation}
When $\dot p(t)\in M_b(\Om\cup\Ga_0;\MD)$, we can consider the Radon-Nikodym derivative $\dot p(t)/|\dot p(t)|$ of  $\dot p(t)$ with respect to its variation $|\dot p(t)|$, which is a function defined $|\dot p(t)|$-a.e.\ on $\Om\cup\Ga_0$. We notice that
$$
\frac{\dot p(t)}{|\dot p(t)|}(x)=\frac{\dot p(t,x)}{|\dot p(t,x)|}
\qquad\hbox{for }\,\Ln\hbox{-a.e.\ }x\in\{|\dot p(t)|>0\}
$$
when $\dot p(t)\in L^2(\Om;\MD)$. It is tempting to consider the inclusion
\begin{equation}\label{pointform}
\frac{\dot p(t)}{|\dot p(t)|}(x)\in N_K(\sigma_D(t,x))
\end{equation}
as a pointwise formulation of the flow rule in the general case 
$\dot p(t)\in M_b(\Om\cup\Ga_0;\MD)$. There is, however, a problem due to the fact that
the left-hand side of (\ref{pointform}) is defined  $|\dot p(t)|$-a.e.\ on $\Om\cup\Ga_0$, while the right-hand side is defined only $\Ln$-a.e.\ on~$\Om$. This difficulty is overcome in Theorem~\ref{strong} below, by introducing a precise representative $\hat\sigma_D(t,x)$ of $\sigma_D(t,x)$, defined almost everywhere with respect to the measure $\mu(t):=\Ln+|\dot p(t)|$. A delicate point in the choice of this representative is the fact that it must also satisfy an integration by parts formula (see Remark~\ref{remstrong}). If $K$ is strictly convex, this representative is essentially unique and can be obtained, in $\Om$, as limit of the averages of~$\sigma_D$ (see Theorem~\ref{strong2}).

\begin{theorem}\label{strong}
Let $t\mapsto (u(t),e(t),p(t))$ be a function from $[0,T]$ into 
$BD(\Om){\times}L^2(\Om;\Mnn)\allowbreak{\times}M_b(\Om\cup\Ga_0;\MD)$, let $\sigma(t):=\C e(t)$, and let $\mu(t):=\Ln+|\dot p(t)|$. Then $t\mapsto (u(t),e(t),p(t))$ is a quasistatic evolution if and only if
\begin{itemize}
\smallskip
        \item[(e)]
$t\mapsto (u(t),e(t),p(t))$ is absolutely continuous and
\begin{itemize}
\smallskip
        \item[(e1)] for every $t\in [0,T]$ we have $(u(t),e(t),p(t))\in A(w(t))$,
        $\sigma(t)\in\Sigma(\Om)\cap \K(\Om)$, $-\div\,\sigma(t)=f(t)$ a.e.\ on $\Om$, and $[\sigma(t)\nu]=g(t)$ on~$\Ga_1$,
        \smallskip
        \item[(e2)]  for a.e.\ $t\in [0,T]$ there exists 
        $\hat\sigma_D(t)\in L^\infty_{\smash{\mu(t)}}(\Om\cup\Ga_0;\MD)$ such that
\begin{eqnarray}
&\displaystyle 
\vphantom{\frac{\dot p(t)}{|\dot p(t)|}}
\hat\sigma_D(t)=\sigma_D(t) \qquad \Ln\hbox{-a.e.\ on }\Om\,,
 \label{a21}
\\
&\displaystyle 
 [\sigma_D(t) {\,:\,} \dot p(t)] =
\Big(\hat\sigma_D(t) {\,:\,} \frac{\dot p(t)}{|\dot p(t)|} \Big)\, |\dot p(t)| 
 \qquad \hbox{on }\Om\cup\Ga_0\,,
 \label{a22}
\\
&\displaystyle 
\frac{\dot p(t)}{|\dot p(t)|}(x)\in N_K(\hat\sigma_D(t,x))\qquad \hbox{for }\,  |\dot p(t)| \hbox{-a.e.\ } x\in\Om\cup\Ga_0\,,
\label{a23}
\end{eqnarray}
where $\hat\sigma_D(t,x)$ denotes the value of $\hat\sigma_D(t)$ at the point~$x$.
\smallskip
\end{itemize}
\end{itemize}
\end{theorem}

\begin{remark}\label{remstrong}
Assume that $t\mapsto (u(t),e(t),p(t))$ is absolutely continuous. If  (e1) holds, then  we can prove, using (\ref{intpartsp}),  that
condition (\ref{a22}) of Theorem~\ref{strong} is equivalent to the following integration by parts formula: for every $\varphi\in C^1(\ol\Om)$ we have
\begin{equation}\label{a22*}
\begin{array}{c}
\langle \varphi\, \hat\sigma_D(t)|\dot p(t)\rangle =
-\langle  \sigma(t)| \varphi\, (\dot e(t)-E\dot w(t))\rangle 
-\langle  \,\sigma(t)| (\dot u(t)- \dot w(t)){\,\odot\,}\nabla  \varphi \rangle +{}
\\
{}+  \langle  f(t)| \varphi\,(\dot u(t)- \dot w(t))\rangle + 
 \langle  g(t)| \varphi\,(\dot u(t)- \dot w(t))\rangle_{\Ga_1}\,,
\end{array}
\end{equation}
where the duality product in the left-hand side is defined by (\ref{dualmeas}).

As $\dot p(t)/|\dot p(t)|=1$ $|\dot p(t)|$-a.e.\ on $\Om\cup\Ga_0$, and $N_K(\xi)=\{0\}$ if $\xi$ is in the interior of $K$, we deduce from (\ref{a23}) that for a.e.\ $t\in[0,T]$
\begin{equation}\label{partialK}
\hat\sigma_D(t,x)\in \partial K
\qquad\hbox{for }|\dot p(t)| \hbox{-a.e.\ }x\in\Om\cup\Ga_0\,.
\end{equation}

Using \cite[Theorem~23.5]{Roc} we can prove that condition (\ref{a23}) is equivalent to
\begin{equation}\label{partialHp}
\hat\sigma_D(t,x)\in \partial H\Big(\frac{\dot p(t)}{|\dot p(t)|}(x)\Big)\qquad \hbox{for }
 |\dot p(t)| \hbox{-a.e.\ } x\in\Om\cup\Ga_0\,.
\end{equation}

Since $\partial H$ is positively homogeneous of degree $0$, this is equivalent to the fact that both the following inclusions are satisfied:
\begin{eqnarray}
&
\hat\sigma_D(t,x)\in
\partial H(\dot p^a(t)(x))
\qquad\hbox{for }  \Ln \hbox{-a.e.\ } x\in\{|\dot p^a(t)|>0\}\,,
\label{hatsigma*1}
\\
&
\displaystyle
\hat\sigma_D(t,x)\in \partial H\Big(\frac{\dot p(t)}{|\dot p(t)|}(x)\Big)
\qquad 
\hbox{for }\, |\dot p^s(t)| \hbox{-a.e.\ } x\in\Om\cup\Ga_0\,.
\label{hatsigma*2}
\end{eqnarray}
\end{remark}

\begin{proof}[Proof of Theorem~\ref{strong}]
Assume that $t\mapsto (u(t),e(t),p(t))$ is a quasistatic evolution. Then $t\mapsto (u(t),e(t),p(t))$ is absolutely continuous by Theorem~\ref{ACesigma} and condition (e1) is satisfied by Theorem~\ref{varineq}.

Let $A(t)\subset\Om$ and $B(t)\subset \Om\cup\Ga_0$ be two disjoint Borel sets such that $A(t)\cup B(t)=\Om\cup\Ga_0$ and
$|\dot p^s(t)|(A(t))=\Ln(B(t))=0$. We define
\begin{eqnarray}
&\hat\sigma_D(t,x):=
\sigma_D(t,x)
\qquad
\text{for  }\Ln\text{-a.e.\ } x\in A(t)\,,
\label{hatsigma1}
\\
&\displaystyle \hat\sigma_D(t,x):=\partial_{\,0} H\Big(\frac{\dot p(t)}{|\dot p(t)|}(x)\Big)
\qquad
\text{for  }\, |\dot p^s(t)| \text{-a.e.\ } x\in B(t)\,,
\label{hatsigma2}
\end{eqnarray}
where $\partial_{\,0} H(\xi)$ denotes the element of $\partial H(\xi)$ with minimum norm.
Then (\ref{a21}) follows from the definition of $\hat\sigma_D(t)$ on $A(t)$ and (\ref{hatsigma*2}) follows from the definition of $\hat\sigma_D(t)$ on $B(t)$. 
To prove (\ref{hatsigma*1}), it is enough to show that
\begin{equation}\label{partialH}
\sigma_D(t,x)\in \partial H(\dot p^a(t)(x))
\qquad \hbox{for }  \Ln \hbox{-a.e.\ } x\in \{|\dot p^a(t)|>0\}\,.
\end{equation}
Taking the absolutely continuous parts in (\ref{H=sigmap}) we obtain 
$H(\dot p^a(t))= \sigma_D(t) {\,:\,} \dot p^a(t)$ $\Ln$-a.e.\ on $\Om$.
Since for $\Ln$-a.e.\ $x\in\Om$ we have $\sigma_D(t,x)\in K=\partial H(0)$ (see, e.g., \cite[Corollary~23.5.3]{Roc}), we obtain $\sigma_D(t,x){\,:\,}\xi\le H(\xi)$ for every $\xi\in \MD$. Therefore for $\Ln$-a.e.\ $x\in \Om$ we have $\sigma_D(t,x){\,:\,}(\xi-  \dot p^a(t))\le H(\xi)-H(\dot p^a(t)(x))$ for every $\xi\in \MD$, which implies~(\ref{partialH}).

To prove (\ref{a22}), we begin by proving the equality on $A(t)$. Since $|\dot p^s(t)|=0$ on $A(t)$, we have $[\sigma_D(t) {\,:\,} \dot p(t)] = \sigma_D(t) {\,:\,} \dot p^a(t)$ on $A(t)$ by (\ref{varsigmap}). As $\hat\sigma_D(t)=\sigma_D(t)$ $\Ln$-a.e.\ on $A(t)$ and $\dot p(t)=\dot p^a(t)$ on $A(t)$, we conclude that
\begin{equation}\label{a22A}
 [\sigma_D(t) {\,:\,} \dot p(t)] = \sigma_D(t) {\,:\,} \dot p^a(t) =
 \Big(\hat\sigma_D(t) {\,:\,} \frac{\dot p(t)}{|\dot p(t)|}\Big)\,  |\dot p(t)| 
 \qquad \hbox{on }A(t)\,.
\end{equation}
To prove the equality on $B(t)$,
we rely on (\ref{H=sigmap}).
Using the definition  (\ref{Hmu}) of $H(\dot p(t))$, the proof of (\ref{a22}) will be complete if we show that
\begin{equation}\label{Hp}
H\Big(\frac{\dot p(t)}{|\dot p(t)|}\Big)=  \hat\sigma_D(t) {\,:\,} \frac{\dot p(t)}{|\dot p(t)|}
\qquad |\dot p(t)| \hbox{-a.e.\ on }B(t)\,.
\end{equation}
But this equality follows from the definition of $\hat\sigma_D(t)$ on $B(t)$, using the Euler identity 
$$
H(\xi)=\zeta{\,:\,}\xi
\qquad\hbox{for every }\xi\in\MD\hbox{ and every }\zeta\in\partial H(\xi)\,. 
$$
This concludes the proof of~(e2).

Conversely, assume~(e). By (\ref{partialHp}), using again the Euler identity, for a.e.\ $t\in [0,T]$ we obtain
$$
H\Big(\frac{\dot p(t)}{|\dot p(t)|}\Big)=  \hat\sigma_D(t) {\,:\,} \frac{\dot p(t)}{|\dot p(t)|}
\qquad |\dot p(t)| \hbox{-a.e.\ on }\Om\cup\Ga_0\,.
$$
{}From the definition (\ref{Hmu}) of the measure $H(\dot p(t))$ and from (\ref{a22}) we deduce that
$\HH(\dot p(t))= \langle\sigma_D(t) | \dot p(t)\rangle$ for a.e.\ $t\in [0,T]$. Therefore 
$t\mapsto (u(t),e(t),p(t))$ is a quasistatic evolution by Theorem~\ref{varineq}.
\end{proof}

For every $r>0$ and every $t\in[0,T]$ we consider the function $\sigma^r(t)\in C(\ol\Om;\Mnn)$ defined by
\begin{equation}\label{sigmar}
\sigma^r(t,x):=\frac{1}{\Ln(B(x,r)\cap \Om)} \int_{B(x,r)\cap \Om}\sigma(t,y)\, dy \,.
\end{equation}
Since $K$ is convex, we have $\sigma^r(t,x)\in K$ for every $x\in \Om$.

If $K$ is strictly convex, i.e., $s\,\xi_1+(1-s)\,\xi_2$ is an interior point of $K$ for every $0<s<1$ and every pair of distinct points $\xi_1$, $\xi_2\in K$, then $H$ is differentiable at all points $\xi\neq 0$ (see, e.g., \cite[Corollary~23.5.4 and Theorem~25.1]{Roc}) and we keep the notation $\partial H(\xi)$ for the gradient. Under this hypothesis, for a.e.\ $t\in [0,T]$ the function $\hat\sigma_D(t)$ is uniquely determined $\mu(t)$-a.e.\ on $\Om\cup\Ga_0$ by (\ref{a21}) and (\ref{partialHp}) as
\begin{eqnarray}
&\hat\sigma_D(t)=\sigma_D(t) \qquad \Ln\hbox{-a.e.\ on }\Om\,,
\label{a21=}
\\
&\displaystyle
\hat\sigma_D(t)= \partial H\Big(\frac{\dot p(t)}{|\dot p(t)|}\Big)\qquad
 |\dot p(t)| \hbox{-a.e.\ on }\Om\cup\Ga_0\,.\label{partialHp=}
\end{eqnarray}
The following theorem shows that, under the same hypothesis, $\hat\sigma_D(t)$ can be obtained in $\Om$ as the limit of $\sigma^r_D(t)$ as ${r\to0}$. This confirms the intrinsic character of the precise representative introduced in Theorem~\ref{strong}.

\begin{theorem}\label{strong2}
Assume that $K$ is strictly convex.
Let $t\mapsto (u(t),e(t),p(t))$ be a quasistatic evolution, let $\mu(t):=\Ln+|\dot p(t)|$, let $\sigma(t):=\C e(t)$, and let 
$\sigma^r(t)$ and $\hat\sigma_D(t)$ be defined by (\ref{sigmar}) and~(\ref{partialHp=}).
Then $\sigma^r_D(t)\to \hat\sigma_D(t)$ strongly in 
$L^1_{\smash{\mu(t)}}(\Om;\MD)$ 
for a.e.\ $t\in[0,T]$. 
\end{theorem}

\begin{proof}
This proof is inspired by the proof of \cite[Theorem~3.7]{Anz}. Since
$\sigma^r_D(t)\to \sigma_D(t)$ strongly in $L^1(\Om;\MD)$ and $\|\sigma^r_D(t)\|_\infty$ is bounded uniformly with respect to $r$, it is enough to prove that $\sigma^r_D(t)\to \hat\sigma_D(t)$ strongly in $L^1_{|\dot p(t)|}(U;\MD)$ for every open set $U\subset\subset\Om$.
Let us fix $U$. Since
$\sigma^r(t)\to \sigma(t)$ strongly in 
$L^2(U;\Mnn)$, $\div\, \sigma^r(t)\to \div\, \sigma(t)$ strongly in 
$L^n(U;\Rn)$, and $\sigma^r_D(t)$ is bounded in $L^\infty(U;\MD)$,
by (\ref{sigmakp}) we have
\begin{equation}\label{sigmarp}
\langle[\sigma^r_D(t) {\,:\,}\dot p(t)]|\varphi \rangle\to \langle [\sigma_D(t) {\,:\,}\dot p(t)]|\varphi\rangle
\end{equation}
for every $\varphi\in C_0(U)$ and for a.e.\ $t\in[0,T]$.
By (\ref{sigmaDp}) we have 
$[\sigma^r_D(t) {\,:\,}\dot p(t)]=\sigma^r_D(t) {\,:\,}\dot p(t)$ on $U$, where the right-hand side is defined by (\ref{sigmaDp2}). By (\ref{H=sigmap}) we have also
$[\sigma_D(t) {\,:\,}\dot p(t)]=H(\dot p(t))$ on $U$. Therefore the definition (\ref{Hmu}) of $H(\dot p(t))$ and (\ref{sigmarp}), together with the boundedness of $\sigma^r_D(t)$, imply that
\begin{equation}\label{sigmarpw*}
\sigma^r_D(t) {\,:\,}\frac{\dot p(t)}{|\dot p(t)|}\ \wto\  H\Big(\frac{\dot p(t)}{|\dot p(t)|}\Big)
\qquad\hbox{weakly}^* \hbox{ in } L^\infty_{|\dot p(t)|}(U)
\end{equation}
for a.e.\ $t\in[0,T]$.

Let us fix $t\in[0,T]$ such that (\ref{partialK}), (\ref{partialHp=}), and (\ref{sigmarpw*}) hold. Since $\sigma^r_D(t)$ is bounded in $L^\infty_{|\dot p(t)|}(U;\MD)$, there exists a sequence $r_j\to 0$ such that $\sigma^{r_j}_D(t)\wto \sigma^*$ for some $\sigma^*\in L^\infty_{|\dot p(t)|}(U;\MD)$. {}From (\ref{sigmarpw*}) we deduce that
\begin{equation}\label{sigma*pH}
\sigma^* {\,:\,}\frac{\dot p(t)}{|\dot p(t)|}\,=\, H\Big(\frac{\dot p(t)}{|\dot p(t)|}\Big)
\qquad |\dot p(t)|\hbox{-a.e.\  on } U\,.
\end{equation}
Let us fix $\xi\in\MD$. Since $\sigma^{r_j}_D(t,x)\in K=\partial H(0)$ for every $x\in U$, we have $\sigma^{r_j}_D(t){\,:\,}\xi\le H(\xi)$ $|\dot p(t)|$-a.e.\ on $U$. As $\sigma^{r_j}_D(t){\,:\,}\xi \wto \sigma^*{\,:\,}\xi$ weakly$^*$ in $L^\infty_{|\dot p(t)|}(U)$, we have also $\sigma^*{\,:\,}\xi\le H(\xi)$ $|\dot p(t)|$-a.e.\ on $U$. Taking (\ref{sigma*pH}) into account, we get
\begin{equation}\label{subdiff}
\sigma^* {\,:\,}\Big(\xi -\frac{\dot p(t)}{|\dot p(t)|}\Big) \,\le\, H(\xi)-H\Big(\frac{\dot p(t)}{|\dot p(t)|}\Big)
\qquad |\dot p(t)|\hbox{-a.e.\  on } U\,. 
\end{equation}
In view of the differentiability properties of $H$, this implies
$$
\sigma^*=\partial H\Big(\frac{\dot p(t)}{|\dot p(t)|}\Big)
\qquad |\dot p(t)|\hbox{-a.e.\  on } U\,. 
$$
By (\ref{partialHp=}) we deduce that $\sigma^*=\hat\sigma_D(t)$ $ |\dot p(t)|$-a.e.\ on $U$. Since the limit does not depend on the sequence $r_j$, we conclude that 
\begin{equation}\label{weak*}
\sigma^r_D(t) \wto \hat\sigma_D(t)\qquad\hbox{weakly}^* \hbox{ in }  L^\infty_{|\dot p(t)|}(U;\MD)\,.
\end{equation}
As $\hat\sigma_D(t,x)\in \partial K$ for  $|\dot p(t)|$-a.e.\ $x\in U$ by Remark~\ref{remstrong} and $\sigma^r_D(t,x)\in K$ for every $x\in U$, the strict convexity of $K$ can be used to improve the weak$^*$ convergence in (\ref{weak*}) and to obtain strong convergence in  $L^1_{|\dot p(t)|}(U;\MD)$ (see, e.g., \cite{Vis}).
\end{proof}

\end{section}

\begin{section}{Appendix}

Let $X$ be the dual of a separable Banach space $Y$. Let $\K$ be a bounded 
closed convex subset
of $Y$ containing the origin as an interior point and let $\HH\colon 
X\to\R$ be its support function, defined by 
$$
\HH(x):=\sup_{y\in \K} \langle x|y \rangle.
$$
Since $\K$ is a bounded neighbourhood of the origin, there exist two 
constants $\alpha_\HH$ and $\beta_\HH$,
with $0<\alpha_\HH\le\beta_\HH<+\infty$, such that
\begin{equation}\label{HHbounds}
\alpha_\HH \|x\|_X\le\HH(x)\le\beta_\HH\|x\|_X \qquad \hbox{for every } x\in X\,.
\end{equation}

Given $f\colon [0,T]\to X$ and $a,b\in[0,T]$ with $a\le b$, we denote the 
total variation of $f$ on $[a,b]$ by
$$
\V(f;a,b):=\sup\Big\{ \sum_{i=1}^N \| f(t_i)-f(t_{i-1})\|_X:\, a=t_0\le 
t_1\le \dots\le t_N=b, \, N\in\N \Big\}\,,
$$
and we define the $\HH$-variation of $f$ on $[a,b]$ as
\begin{equation}\label{pdiss}
\V_\HH(f;a,b):=\sup\Big\{ \sum_{i=1}^N \HH(f(t_i)-f(t_{i-1})):\, a=t_0\le 
t_1\le \dots\le t_N=b, \, N\in\N
\Big\}\,.
\end{equation}
{}From (\ref{HHbounds}) it follows that
$$
\alpha_\HH \V(f;a,b)\le \V_\HH(f;a,b)\le \beta_\HH \V(f;a,b)\,.
$$
Since $\HH$ is weakly$^*$ lower semicontinuous, we have
\begin{equation}\label{semidiss}
\V_\HH(f;a,b)\le \liminf_{k\to\infty} \V_\HH(f_k;a,b)
\end{equation}
whenever $f_k(t)\wto f(t)$  weakly$^*$ for every $t\in[a,b]$.

We now prove a theorem about weak$^*$ derivatives of absolutely continuous functions with values in~$X$ 
and their relationships with the notion of $\HH$-variation.

\begin{theorem}\label{thm:a01}
Let $f\colon [0,T]\to X$ be an absolutely continuous function. Then the weak$^*$-limit 
\begin{equation}\label{a00}
\dot f(t):= w^*\hbox{-}\lim_{s\to t}\frac{f(s)-f(t)}{s-t}
\end{equation}
exists for a.e.~$t\in[0,T]$, and
\begin{equation}\label{metder}
\HH(\dot f(t))= \lim_{s\to t}\HH\Big(\frac{f(s)-f(t)}{s-t}\Big)
\end{equation}
for a.e.~$t\in[0,T]$.
Moreover, the function $t\mapsto \HH(\dot f(t))$ is measurable and
\begin{equation}\label{Dder0}
\V_\HH(f;a,b)=\int_a^b \HH(\dot f(t))\, dt 
\end{equation}
for every $a,b\in[0,T]$ with $a\le b$.
\end{theorem}

\begin{proof}
Let $F$ be the linear span over $\Q$ of a countable dense set in $Y$.
For every $y\in F$ the map $t\mapsto\langle f(t)|y\rangle$ is absolutely continuous on
$[0,T]$; therefore, there exists a set $N_y$ of measure zero such that the limit
$$
D_y(t):= \lim_{s\to t}\frac{\langle f(s)-f(t)|y\rangle}{s-t}
$$
exists for every $t\in [0,T]\setmeno N_y$.
Let $\V(t):=\V(f;0,t)$. Since the function $t\mapsto \V(t)$ is
non-decreasing, it is differentiable for every $t\in[0,T]\setmeno M$, where
$M$ is a set of measure zero.
Let $N$ be the union of $M$ with the sets $N_y$ for $y\in F$. Then, ${\mathcal 
L}^1(N)=0$, the derivative $D_y(t)$ exists for every $y\in F$ and every $t\in
[0,T]\setmeno N$, and
\begin{equation}\label{m01}
|D_y(t)|=\lim_{s\to t} \frac{|\langle f(s)-f(t)|y\rangle|}{|s-t|}\le \dot\V(t)\|y\|_Y
\end{equation}
for every $y\in F$ and every $t\in [0,T]\setmeno N$. 
Now, for $t\in [0,T]\setmeno N$ consider the linear map $y\in F\mapsto 
D_y(t)$. This map is continuous by (\ref{m01}); therefore, there exists a vector in $X$, 
which we call $\dot f(t)$, such that
$$
D_y(t)=\langle \dot f(t)|y\rangle
$$
for every $y\in F$.
Using the density of $F$ and (\ref{m01}) it is easy to show that the 
vector $\dot f(t)$ satisfies
$$
\langle \dot f(t)|y\rangle=\lim_{s\to t}\frac{\langle 
f(s)-f(t)|y\rangle}{s-t}
$$
for every $y\in Y$ and every $t\in [0,T]\setmeno N$, so that (\ref{a00}) 
is satisfied.

We note that the function $t\mapsto \HH(\dot f(t))$ is measurable, since
the map $t\to \langle \dot f(t)|y\rangle$ is measurable for every $y\in Y$ 
and $\HH(\dot f(t))=\sup_{y\in \K_0} \langle \dot f(t)|y\rangle$, where $\K_0$ is a countable
dense subset of $\K$.
Moreover, if $a=t_0\le t_1\le \dots \le t_{N-1}\le  t_N=b$ is a 
subdivision of $[a,b]$, then
$$
\langle f(t_i)-f(t_{i-1})|y\rangle =\int_{t_{i-1}}^{t_i}
\langle \dot f(t)|y\rangle \, dt
\leq \int_{t_{i-1}}^{t_i}\HH(\dot f(t))\, dt
$$
for every $1\leq i \leq N$ and every $y\in \K$, hence
$$
\HH(f(t_i)-f(t_{i-1}))\leq \int_{t_{i-1}}^{t_i}\HH(\dot f(t))\,dt
$$
for every $1\leq i \leq N$. Summing over $i$ and taking the supremum over 
all subdivisions, we obtain
\begin{equation}\label{uz01}
\V_\HH(f;a,b)\leq \int_a^b \HH(\dot f(t))\, dt \,.
\end{equation}

To show the converse inequality, note that the function 
$t\mapsto\V_\HH(f;0,t)$ is non-decreasing;
therefore, it is differentiable for a.e. $t\in[0,T]$ and
\begin{equation}\label{m011}
\int_a^b \frac{d}{dt}\V_\HH(f;0,t)\, dt \le \V_\HH(f;a,b)\,.
\end{equation}
Let $t_0\in[0,T]$ be a point where both $f$ and $\V_\HH(f;0,\cdot)$ are 
differentiable.
Since $\HH$ is positively homogeneous of degree $1$, we have
$$
\HH\Big(\frac{f(t)-f(t_0)}{t-t_0}\Big)\leq 
\frac{\V_\HH(f;0,t)-\V_\HH(f;0,t_0)}{t-t_0}\,
$$
for every $t\neq t_0$. Passing to the limit as $t \to t_0$ and 
using the weak$^*$-lower semicontinuity of $\HH$, we get
$$
\HH(\dot f(t_0))\leq \liminf_{t \to t_0} 
\HH\Big(\frac{f(t)-f(t_0)}{t- t_0}\Big)
\leq \limsup_{t \to t_0} 
\HH\Big(\frac{f(t)-f(t_0)}{t- t_0}\Big)
\leq \frac{d}{dt}\V_\HH(f;0,t)\Big|_{t=t_0}
$$
for a.e.~$t_0 \in[0,T]$. We now integrate the first and the 
last term in the previous inequality from $a$ to $b$ and we obtain
(\ref{Dder0}) and (\ref{metder})  from  (\ref{uz01}) and~(\ref{m011}).
\end{proof}

We conclude this appendix with a lemma which generalizes the classical Helly Theorem for real valued functions with uniformly bounded variation, as well as its extension to reflexive separable Banach spaces (see, e.g., \cite[Chapter~1,Theorem~3.5]{Bar-Pre}).

\begin{lemma}\label{lm:44}
Let $f_k\colon [0,T]\to X$ be a sequence of functions such that $f_k(0)$ and $\V(f_k;0,T)$ are bounded uniformly with respect to~$k$. Then
there exist a subsequence, still denoted $f_k$,
and a function $f\colon [0,T]\to X$  with bounded variation on $[0,T]$, such that $f_k(t)\wto f(t)$ weakly$^*$ for every $t\in [0,T]$.
\end{lemma}

\begin{proof}
It is enough to apply \cite[Theorem~3.2]{Mai-Mie} with ${\mathcal Y}=X$, 
${\mathcal R}(t)= {\mathcal V}(t)$ equal to the corresponding unit ball, and $\mathcal T$ equal to the weak$^*$ topology.
\end{proof}

\end{section}
\bigskip

\noindent {\bf Acknowledgments.} { This work is part of the Project ``Calculus of Variations" 2002, 
supported by the Italian Ministry of Education, University, and Research.}

\bigskip
\bigskip

{
}

\end{document}